\documentclass[a4paper,leqno,draft]{article}
\usepackage{amssymb,amsmath,amsfonts,amsthm,
mathrsfs}
\newtheorem{theorem}{Theorem}[section]
\newtheorem{corollary}[theorem]{Corollary}
\newtheorem{lemma}[theorem]{Lemma}
\newtheorem{proposition}[theorem]{Proposition}
\newtheorem{definition}[theorem]{Definition}
\theoremstyle{definition}
\newtheorem{remark}[theorem]{Remark}
\newtheorem{example}[theorem]{Example}

\newcommand{\wt}[1]{\widetilde{#1}}

\newcommand{\Cinfc}{\ensuremath{\mathcal{C}^\infty_{\text{c}}}}

\renewcommand{\S}{\mathscr{S}}
\newcommand{\E}{\ensuremath{{\cal E}}}

\renewcommand{\L}{\mathcal{L}}

\newcommand{\Ker}{\mathop{{\mathrm{Ker}}}}

\newcommand{\mb}[1]{\ensuremath{\mathbb{#1}}}
\newcommand{\N}{\mb{N}}

\newcommand{\R}{\mb{R}}
\newcommand{\C}{\mb{C}}

\newcommand{\G}{\ensuremath{{\cal G}}}

\newcommand{\Gc}{\ensuremath{{\cal G}_\mathrm{c}}}

\newcommand{\GS}{\G_{{\, }\atop{\hskip-4pt\scriptstyle\S}}\!}
\newcommand{\EM}{\ensuremath{{\cal E}_{M}}}

\newcommand{\Neg}{\mathcal{N}}

\newfont{\bigmath}{cmr12 at 13pt}

\newfont{\grecomath}{cmmi12 at 15pt}

\newcommand{\val}{\mathrm{v}} 
\newcommand{\esp}{\mathrm{e}}

\newfont{\bl}{msbm10 scaled \magstep2}

\newcommand{\beq}{\begin{equation}}
\newcommand{\eeq}{\end{equation}}

\newcommand{\F}{\ensuremath{{\cal F}}}

\newcommand{\notmid}{\mid\kern-0.5em\not\kern0.5em}

\newcommand{\inner}[3][\empty]{\ifx#1\empty\left( #2|#3\right)\else#1( #2|#3 #1)\fi}

\newcommand{\norm}[2]{{\| #1 \|}_{#2}}

\newcommand{\eps}{\varepsilon}

\newcommand{\Om}{\Omega}

\newcommand{\compl}[1]{{#1}^{\mathrm{c}}}

\renewcommand{\Re}{\ensuremath{\mathrm{Re}}}
\renewcommand{\Im}{\ensuremath{\mathrm{Im}}}

\newcommand{\M}{\mathcal{M}}
\newcommand{\mF}{\mathcal{F}}
\newcommand{\mH}{\mathcal{H}}
\newcommand{\mP}{\mathcal{P}}
\newcommand{\mQ}{\mathcal{Q}}

\newcommand{\abs}[2][\empty]{\ifx#1\empty\left|#2\right|%
\else#1\vert #2 #1\vert\fi}
\newcommand{\Ann}{\mathop{\mathrm{Ann}}}
\newcommand{\co}[1]{{#1}^{c}}
\newcommand{\ideal}{\unlhd}

\newcommand{\K}{{\mathbb K}}
\renewcommand{\norm}[2][]{{{\Vert #2 \Vert}_{#1}}}

\newcommand{\caninf}{\ensuremath[(\eps)_\eps]}

\begin{document}

\title{{\bf Hilbert $\wt{\C}$-modules: structural properties and applications to variational problems}}
\author{Claudia Garetto\footnote{Supported by FWF (Austria), grants T305-N13 and  Y237-N13.} \\ 
\texttt{claudia@mat1.uibk.ac.at}\\
\\
Hans Vernaeve\footnote{Supported by FWF (Austria), grants M949-N18 and  Y237-N13.}\\
\texttt{hans.vernaeve@uibk.ac.at}\\
\\
Institut f\"ur Grundlagen der Bauingenieurwissenschaften,\\
Arbeitsbereich f\"ur Technische Mathematik,\\
Leopold-Franzens-Universit\"at Innsbruck, Austria
}
\date{}
\maketitle

\begin{abstract}
We develop a theory of Hilbert $\wt{\C}$-modules by investigating their structural and functional analytic properties. Particular attention is given to finitely generated submodules, projection operators, representation theorems for $\wt{\C}$-linear functionals and $\wt{\C}$-sesquilinear forms. By making use of a generalized Lax-Milgram theorem, we provide some existence and uniqueness theorems for variational problems involving a generalized bilinear or sesquilinear form.
\end{abstract}

\setcounter{section}{-1}
\section{Introduction}
This paper is part of a wide project which aims to introduce functional analytic methods into Colombeau algebras of generalized functions. Our intent is to deal with the general problem of existence and qualitative properties of solutions of partial differential equations in the Colombeau setting, by means of functional analytic tools adapted and generalized to the range of topological $\wt{\C}$-modules. Starting from the topological background given in \cite{Garetto:05a, Garetto:05b}, in this paper we develop a theory of Hilbert $\wt{\C}$-modules. This will be the framework where to investigate variational equalities and inequalities generated by highly singular problems in partial differential equations. 

A first example of a Hilbert $\wt{\C}$-module is the Colombeau space $\G_H$ of generalized functions based on a Hilbert space $(H, \inner{\cdot}{\cdot})$ \cite[Definition 3.1]{Garetto:05a}, where the scalar product is obtained by letting $\inner{\cdot}{\cdot}$ act componentwise at the representatives level. A number of theorems, such as projection theorem, Riesz-representation and Lax-Milgram theorem can be obtained in a direct way when we work on $\G_H$, by applying the corresponding classical results at the level of representatives at first and then by checking the necessary moderateness conditions. This is a sort of transfer method, which has been exclusively employed so far, of producing results in a Colombeau context deeply related to a classical one. The novelty of our work is the introduction of a general notion of a Hilbert $\wt{\C}$-module which has no more the internal structure of $\G_H$, and the completely intrinsic way of developing a topological and functional analytic theory within this abstract setting. As we will see in the course of the paper, the wide generality of our approach on the one hand entails some technicalities in the proofs and on the other hand leads to the introduction of a number of new concepts, such as edged subsets of a $\wt{\C}$-module, normalization property, etc.

We now describe the contents of the sections in more detail.  

The first section serves to collect some basic notions necessary for the comprehension of the paper. We begin in Subsection \ref{subsection_Colombeau} by recalling the definition of the Colombeau space $\G_E$ of generalized functions based on a locally convex topological vector space $E$. In order to view $\G_E$ as a particular example of a  locally convex topological $\wt{\C}$-module, where $\wt{\C}$ is the ring $\G_\C$ of generalized constants, we make use of concepts as valuation and ultra-pseudo-seminorm and of some fundamental ingredients of the theory of topological $\wt{\C}$-modules elaborated in \cite{Garetto:05a, Garetto:05b}. Particular attention is given to $\wt{\C}$-linear maps and $\wt{\C}$-sesquilinear forms acting on locally convex topological $\wt{\C}$-modules and to their basic structure when we work on spaces of $\G_E$-type \cite[Definition 1.1]{Garetto:06a}. We introduce the property of being internal for subsets of $\G_E$ as the analogue of the basic structure for maps. Internal subsets will play a main role in the paper, in the existence and uniqueness theorems for variational equalities and inequalities of Sections \ref{section_varia} and \ref{sec_appl}. We conclude Subsection \ref{subsection_Colombeau} by discussing some issues concerning the ring $\wt{\R}$ of real generalized numbers: definition and properties of the order relation $\ge$, invertibility and negligibility with respect to a subset $S$ of $(0,1]$, characterization of zero divisors and idempotent elements, infimum and close infimum in $\wt{\R}$.

The second part of Section \ref{top_notions} deals with the class of $\wt{\C}$-modules with $\wt{\R}$-seminorms. Making use of the order relation $\ge$ in $\wt{\R}$ and of the classical notion of seminorm as a blueprint, we introduce the concept of $\wt{\R}$-seminorm on a $\wt{\C}$-module $\G$. This induces a topology on $\G$ which turns out to be $\wt{\C}$-locally convex. In other words we find a special class of locally convex topological $\wt{\C}$-modules which contains the spaces of generalized functions based on a locally convex topological vector space as a particular case.

Section \ref{section_Hilbert} is devoted to the definition and the first properties of the family of topological $\wt{\C}$-modules which are the mathematical core of the paper: the Hilbert $\wt{\C}$-modules. They are defined by means of a generalized scalar product $\inner{\cdot}{\cdot}$ with values in $\wt{\C}$ which determines the $\wt{\R}$-norm $\Vert u\Vert=\inner{u}{u}^{\frac{1}{2}}$. This means that they are particular $\wt{\R}$-normed $\wt{\C}$-modules. As first examples of Hilbert $\wt{\C}$-modules we consider the space $\G_H$ based on a vector space $H$ with scalar product and more generally given a net $(H_\eps,\inner{\cdot}{\cdot}_{H_\eps})_\eps$, the quotient of the corresponding moderate nets over negligible nets (Proposition \ref{prop_ex_3}).

In the intent of developing a topological and functional analytic theory of Hilbert $\wt{\C}$-modules, we start in Subsection \ref{subsec_projec} by investigating the notion of projection on a suitable subset $C$ of a Hilbert $\wt{\C}$-module $\G$. This requires some new assumptions on $C$, such as being reachable from a point $u$ of $\G$, the property of being edged, i.e., reachable from any $u$, and a formulation of convexity in terms of $\wt{\R}$-linear combinations which resembles the well-known classical definition for subsets of a vector space but differs from the $\wt{\C}$-convexity introduced in \cite{Garetto:05a}. In detail, we prove that if $C$ is a closed nonempty subset of the Hilbert $\wt{\C}$-module $\G$ such that $C+C\subseteq 2C$ and it is reachable from $u\in\G$, i.e.\ the set $\{\Vert u-w\Vert,\, w\in C\}$ has a close infimum in $\wt{\R}$, then the projection $P_C(u)$ of $u$ on $C$ exists. The operator $P_C$ is globally defined and continuous when $C$ is edged and is $\wt{\C}$-linear when $C=M$ is a closed and edged $\wt{\C}$-submodule of $\G$. We also see, by means of a counterexample, that the condition of being edged is necessary for the existence of $P_M$ and that this operator allows us to extend any continuous $\wt{\C}$-linear map with values in a topological $\wt{\C}$-module from $M$ to the whole of $\G$. In this way we obtain a version of the Hahn-Banach theorem where the fact that $M$ is edged is essential. Moreover, closed and edged submodules of $\G$ can be characterized as those submodules $M$ for which $M + M^\perp=\G$.  

Section \ref{sec_fin} gives a closer look at edged submodules of a Hilbert $\wt{\C}$-module. For the sake of generality we work in the context of $\wt{\K}$-modules, where $\K$ is $\R$ or $\C$, and we state many results in the framework of Banach $\wt{\K}$-modules. In our investigation on submodules we distinguish between cyclic submodules, i.e.\ generated by one element, and submodules generated by $m>1$ elements. In particular, we prove that when a submodule is finitely generated the property of being edged is deeply related to topological closedness and to some structural properties of the generators. We carefully comment our statements providing explanatory examples and counterexamples.

The main topic of Section \ref{section_riesz} is the formulation of a Riesz representation theorem for continuous $\wt{\C}$-linear functionals acting on a Hilbert $\wt{\C}$-module. We prove that a functional $T$ can be written in the form $T(u)=\inner{u}{c}$ if and only if there exists a closed and cyclic $\wt{\C}$-submodule $N$ such that $N^\perp\subseteq\Ker T$. In particular, on $\G_H$, the Riesz representation theorem gives necessary and sufficient conditions for a $\wt\C$-linear functional to be basic. The structural properties of continuous $\wt{\C}$-sesquilinear forms on Hilbert $\wt{\C}$-modules are investigated by making use of the previous representation theorem.

In Section \ref{section_cont} we concentrate on continuous $\wt{\C}$-linear operators acting on a Hilbert $\wt{\C}$-module. In detail, we deal with isometric, unitary, self-adjoint and projection operators obtaining the following characterization: $T$ is a projection operator (i.e.\ self-adjoint and idempotent) if and only if it is the projection $P_M$ on a closed and edged $\wt{\C}$-submodule $M$. 

A version of the Lax-Milgram theorem valid for Hilbert $\wt{\C}$-modules and $\wt{\C}$-sesquilinear forms is proved in Section \ref{section_lax} for forms of the type $a(u,v)=\inner{u}{g(v)}$, when we assume that the range of $g$ is edged and that $a$ satisfies a suitable coercivity condition. This theorem applies to any basic and coercive $\wt{\C}$-sesquilinear form on $\G_H$ and plays a relevant role in the applications of the last section of the paper.

Section \ref{section_varia} concerns variational inequalities involving a continuous $\wt{\R}$-bilinear form in the framework of Hilbert $\wt{\R}$-modules. Under suitable hypotheses on the set $C\subseteq\G$ we prove that the problem
\[
a(u,v-u)\ge \inner{f}{v-u},\qquad\qquad \text{for all $v\in C$}
\]
is uniquely solvable in $C$ if $a$ is a symmetric, coercive and continuous $\wt{\R}$-bilinear form and the functional $I(u)=a(u,u)-2\inner{f}{u}$ has a close infimum on $C$ in $\wt{\R}$. This applies to the case of basic and coercive forms on $\G_H$ when $C$ is internal and can be extended to basic $\wt{\R}$-sesquilinear forms which are non symmetric via some contraction techniques. The theorems of Section \ref{section_varia} are one of the first examples of existence and uniqueness theorems in the Colombeau framework obtained in an intrinsic way via topological and functional analytic methods. 

The paper ends by discussing some concrete problems coming from partial differential operators with highly singular coefficients, which in variational form can be solved by making use of the theorems on variational equalities and inequalities of Section \ref{section_varia}. The generalized framework within which we work allows us to approach problems which are not solvable classically and to get results consistent with the classical ones when the latter exist. 

\section{Basic notions}
\label{top_notions}
This section of preliminary notions provides some topological background necessary for the comprehension of the paper. Particular attention is given to Colombeau spaces of generalized functions, locally convex topological $\wt{\C}$-modules and topological $\wt{\C}$-modules with $\wt{\R}$-seminorms. Main references are \cite{Garetto:05a, Garetto:05b, Garetto:06a}.

\subsection{Colombeau spaces of generalized functions and topological $\wt{\C}$-modules}
\label{subsection_Colombeau}
\subsubsection{First definitions, valuations and ultra-pseudo-seminorms}
Let $E$ be a locally convex topological vector space topologized through the family of seminorms $\{p_i\}_{i\in I}$. The elements of 
\[
\begin{split} 
\M_E &:= \{(u_\eps)_\eps\in E^{(0,1]}:\, \forall i\in I\,\, \exists N\in\N\quad p_i(u_\eps)=O(\eps^{-N})\, \text{as}\, \eps\to 0\},\\
\Neg_E &:= \{(u_\eps)_\eps\in E^{(0,1]}:\, \forall i\in I\,\, \forall q\in\N\quad p_i(u_\eps)=O(\eps^{q})\, \text{as}\, \eps\to 0\},
\end{split}
\]
are called $E$-moderate and $E$-negligible, respectively. The space of \emph{Colombeau generalized functions based on $E$} is defined as the quotient $\G_E := \M_E / \Neg_E$. 

The rings $\wt{\C}=\EM/\Neg$ of \emph{complex generalized numbers} and $\wt{\R}$ of \emph{real generalized numbers} are obtained by taking $E=\C$ and $E=\R$ respectively. One can easily see that for any locally convex topological vector space $E$ (on $\C$) the space $\G_E$ has the structure of a $\wt{\C}$-module. We use the notation $u=[(u_\eps)_\eps]$ for the class $u$ of $(u_\eps)_\eps$ in $\G_E$. This is the usual way adopted in the paper to denote an equivalence class.

$\wt{\C}$ is trivially a module over itself and it can be endowed with a structure of a topological ring. This is done by defining the \emph{valuation} $\val$ of a representative $(r_\eps)_\eps$ of $r\in\wt{\C}$ as $\sup\{b\in \R:\, |r_\eps|=O(\eps^b)\, \text{as $\eps\to 0$}\}$. By observing that $\val((r_\eps)_\eps)=\val((r'_\eps)_\eps)$ for all representatives $(r_\eps)_\eps, (r'_\eps)_\eps$ of $r$, one can let $\val$ act on $\wt{\C}$ and define the map  
\[
\vert\cdot\vert_\esp := \wt{\C}\to [0,+\infty): u\to \vert u\vert_\esp:=\esp^{-\val(u)}.
\]
The properties of the valuation on $\wt{\C}$ make the coarsest topology on $\wt{\C}$ for which the map $\vert \cdot\vert_\esp$ is continuous compatible with the ring structure.
It is common in the already existing literature \cite{NPS:98, Scarpalezos:92, Scarpalezos:98, Scarpalezos:00} to use  the adjective ``sharp'' for such a topology.  

A \emph{topological $\wt{\C}$-module} is a $\wt{\C}$-module $\G$ endowed with a $\wt{\C}$-linear topology, i.e., with a topology such that the addition $\G\times\G\to\G:(u,v)\to u+v$ and the product $\wt{\C}\times\G\to \G:(\lambda,u)\to \lambda u$ are continuous. A \emph{locally convex topological $\wt{\C}$-module} is a topological $\wt{\C}$-module whose topology is determined by a family of \emph{ultra-pseudo-seminorms}. As defined in \cite[Definition 1.8]{Garetto:05a} an ultra-pseudo-seminorm on $\G$ is a map $\mP:\G\to[0,+\infty)$ such that 
\begin{trivlist}
\item[(i)] $\mP(0)=0$,
 \item[(ii)] $\mP(\lambda u)\le \vert\lambda\vert_\esp \mP(u)$ for all $\lambda\in\wt{\C}$, $u\in\G$,
\item[(iii)] $\mP(u+v)\le\max\{\mP(u),\mP(v)\}$.
\end{trivlist}
Note that since $|[(\eps^{-a})_\eps]|_\esp={|[(\eps^a)_\eps]|_\esp}^{-1}$, from (ii) it follows 
\begin{trivlist}
\item[(ii)'] $\mP(\lambda u)= \vert\lambda\vert_\esp \mP(u)$ for all $\lambda=[(c\eps^a)_\eps]$, $c\in\C$, $a\in\R$, $u\in\G$.
\end{trivlist}
The notion of \emph{valuation} can be introduced in the general context of $\wt{\C}$-modules as follows: a valuation on on $\G$ is a function $\val:\G\to(-\infty,+\infty]$ such that
\begin{trivlist}
\item[(i)] $\val(0)=+\infty$,
\item[(ii)] $\val(\lambda u)\ge \val_{\wt{\C}}(\lambda)+\val(u)$ for all $\lambda\in\wt{\C}$, $u\in\G$,
\item[(iii)] $\val(u+v)\ge\min\{\val(u),\val(v)\}$.
\end{trivlist}
As above, from (ii) it follows that 
\begin{trivlist}
\item[(ii)'] $\val(\lambda u)= \val_{\wt{\C}}(\lambda)+\val(u)$ for all $\lambda=[(c\eps^a)_\eps]$, $c\in\C$, $a\in\R$, $u\in\G$.
\end{trivlist}
Any valuation generates an ultra-pseudo-seminorm by setting $\mP(u)=\esp^{-\val(u)}$. An ultra-pseudo-seminorm $\mP$ such that $\mP(u)=0$ if and only if $u=0$ is called \emph{ultra-pseudo-norm}. The topological dual of a topological $\wt{\C}$-module $\G$ is the set $\L(\G,\wt{\C})$ of all continuous and $\wt{\C}$-linear functionals on $\G$. A thorough investigation of $\L(\G,\wt{\C})$ can be found in \cite{Garetto:05a, Garetto:05b} together with interesting examples coming from Colombeau theory.

The family of seminorms $\{p_i\}_{i\in I}$ on $E$ equips $\G_E$ with the structure of a locally convex topological $\wt{\C}$-module by means of the valuations
\[
\val_{p_i}([(u_\eps)_\eps]):=\val_{p_i}((u_\eps)_\eps):=\sup\{b\in\R:\quad p_i(u_\eps)=O(\eps^b)\, \text{as $\eps\to 0$}\}
\] 
and the corresponding \emph{ultra-pseudo-seminorms} $\{\mP_i\}_{i\in I}$, where $\mP_i(u)=\esp^{-\val_{p_i}(u)}$. 

\subsubsection{Basic $\wt{\C}$-linear maps and $\wt{\C}$-sesquilinear forms}
Let $(\G,\{\mP_i\}_{i\in I})$ and $(\mF,\{\mQ_j\}_{j\in J})$ be locally convex topological $\wt{\C}$-modules. Theorem 1.16 and Corollary 1.17 in \cite{Garetto:05a} prove that a $\wt{\C}$-linear map $T:\G\to\mF$ is continuous if and only if it is continuous at the origin, if and only if for all $j\in J$ there exists a constant $C>0$ and a finite subset $I_0$ of $I$ such that the inequality
\beq
\label{cont_T_E}
\mQ_j(Tu)\le C\max_{i\in I_0}\mP_i(u)
\eeq
holds for all $u\in\G$.

In the particular case of $\G=\G_E$ and $\mF=\G_F$, we recall that a $\wt{\C}$-linear map $T:\G_E\to \G_F$ is \emph{basic} if there exists a net $(T_\eps)_\eps$ of continuous linear maps from $E$ to $F$ fulfilling the continuity-property  
\beq
\label{basic_T_E}
\forall j\in J\, \exists I_0\subseteq I\, \text{finite}\ \exists N\in\N\, \exists\eta\in(0,1]\, \forall u\in E\, \forall\eps\in(0,\eta]\qquad
q_j(T_\eps u)\le \eps^{-N}\sum_{i\in I_0}p_i(u),
\eeq
and such that $Tu=[(T_\eps(u_\eps))_\eps]$ for all $u\in\G_E$. It is clear that \eqref{basic_T_E} implies \eqref{cont_T_E} and therefore any basic map is continuous. 

This notion of basic structure can be easily extended to multilinear maps from $\G_{E_1}\times...\times\G_{E_n}\to \G_F$. In this paper we will often work with basic $\wt{\C}$-sesquilinear forms. A \emph{basic} $\wt{\C}$-sesquilinear form $a$ on $\G_E\times\G_F$ is a $\wt{\C}$-sesquilinear map $a$ from $\G_E\times\G_F\to\wt{\C}$ such that there exists a net $(a_\eps)_\eps$ of continuous sesquilinear forms on $E\times F$ fulfilling the continuity-property  
\begin{multline}
\label{basic_a_E}
\exists I_0\subseteq I\, \text{finite}\, \exists J_0\subseteq J\, \text{finite}\ \exists N\in\N\, \exists\eta\in(0,1]\, \forall u\in E\, \forall v\in F\, \forall\eps\in(0,\eta]\\
|a_\eps(u,v)|\le \eps^{-N}\sum_{i\in I_0}p_i(u) \sum_{j\in J_0}q_j(v),
\end{multline}
and such that $a(u,v)=[(a_\eps(u_\eps,v_\eps))_\eps]$ for all $u\in\G_E$ and $v\in\G_F$.

\subsubsection{Internal subsets of $\G_E$}
A subset $A\subseteq\G_E$ is called internal \cite{OV:07} if there exists a net $(A_\eps)_{\eps\in(0,1)}$ of subsets $A_\eps\subseteq E$ such that
\[A=\{u\in\G_E : \exists \textrm{ representative }(u_\eps)_\eps \text{ of }u\, \exists\eps_0\in(0,1)\, \forall \eps\le\eps_0\, u_\eps\in A_\eps\}.\]
If all $A_\eps\neq\emptyset$, then we can take $\eps_0=1$ in the previous definition without loss of generality.
The internal set corresponding with the net $(A_\eps)_\eps$ is denoted by $[(A_\eps)_\eps]$.\\
Let $E$ be a normed vector space and $A$ an internal subset of $\G_E$. Then the following hold \cite{OV:07}:
\begin{itemize}
\item[(i)] $A$ is closed.
\item[(ii)] Let $u\in\G_E$. If $A$ is not empty, then there exists $v\in A$ such that $\norm{u-v}=\min_{w\in A}{\norm{u-w}}$ \cite{OV:07}.
\end{itemize}

\subsubsection{Some properties of the ring of real generalized numbers}
We finally concentrate on the ring $\wt{\R}$ of real generalized numbers. It can be equipped with the order relation $\le$ given by $r\le s$ if and only if there exist $(r_\eps)_\eps$ and $(s_\eps)_\eps$ representatives of $r$ and $s$ respectively such that $r_\eps\le s_\eps$ for all $\eps\in(0,1]$. We say that $r\in\wt{\R}$ is nonnegative iff $0\le r$. We write $r>0$ if and only if $r\ge 0$ and $r\neq 0$. Equipped with this order, $\wt\R$ is a partially ordered ring. One can define the square root of a nonnegative generalized number $r\in\wt{\R}$ by setting $r^{\frac{1}{2}}=[(r_\eps^{\frac{1}{2}})_\eps]$, for any representative $(r_\eps)_\eps$ of $r$ such that $r_\eps\ge 0$ for all $\eps$. We leave it to the reader to check that $|r^2|_\esp=|r|_\esp^2$ for all $r\in\wt{\R}$ and that $|r^{\frac{1}{2}}|_\esp=|r|_\esp^{\frac{1}{2}}$ for all $r\ge 0$. In the sequel we collect some further properties concerning the order relation in $\wt{\R}$ which will be useful in the course of the paper.
\begin{proposition}
Let $a,b,b_n,r$ be real generalized numbers. The following assertions hold:
\begin{itemize}
\item[(i)] $r\ge 0$ if and only if there exists a representative $(r_\eps)_\eps$ of $r$ such that $r_\eps\ge 0$ for all $\eps\in(0,1]$ if and only if for all representatives $(r_\eps)_\eps$ of $r$ and for all $q\in\N$ there exists $\eta\in(0,1]$ such that $r_\eps\ge-\eps^q$ for all $\eps\in(0,\eta]$;
\item[(ii)] if $a\ge 0$, $b\ge 0$ and $a^{2}\le b^{2}$ then $a\le b$;
\item[(iii)] if $a\le b_n$ for all $n\in\N$ and $b_n\to b$ as $n\to\infty$ then $a\le b$.
\end{itemize}
\end{proposition}
 
Let $S\subseteq(0,1]$. We denote by $e_S\in\wt\R$ the generalized number with the characteristic function $(\chi_{S}(\eps))_\eps$ as representative, and $\co S=(0,1]\setminus S$. Then clearly, $e_S\ne 0$ iff $0\in \overline S$ and $e_S\ne 1$ iff $0\in\overline{\co S}$.\\
Let $z\in\wt\C$ and $S\subseteq(0,1]$ with $e_S\ne 0$. Then $z$ is called invertible w.r.t.\ $S$ if there exists $z'\in\wt\C$ such that $z z'=e_S$; $z$ is called zero w.r.t.\ $S$ if $ze_S=0$. The following holds \cite{Vernaeve:07a}:\\
Let $(z_\eps)_\eps$ be a representative of $z$.
\begin{itemize}
\item[(i)] $z$ is zero w.r.t\ $S$ iff $(\forall m\in\N)(\exists \eta >0)(\forall\eps\in S\cap (0,\eta))(\abs{z_\eps}\le\eps^m)$ iff $(\forall T\subseteq S$ with $e_T\ne 0)(z$ is not invertible w.r.t.\ $T)$;
\item[(ii)] $z$ is invertible w.r.t.\ $S$ iff $(\exists m\in\N)(\exists \eta >0)(\forall\eps\in S\cap (0,\eta))(\abs{z_\eps}\ge\eps^m)$ iff  $(\forall T\subseteq S$ with $e_T\ne 0)(z$ is not zero w.r.t.\ $T)$.
\end{itemize}
Finally we have the following characterizations of the zero divisors and the idempotent elements of $\wt{\C}$. 
Let $z,z'\in\wt{\C}$ such that $zz'=0$. Then, there exists $S\subseteq(0,1]$ such that $ze_S=0$ and $z'e_{\co{S}}=0$ \cite{Vernaeve:07a}. If $z=z^2$ then there exists $S\subseteq (0,1]$ such that $z=e_S$ \cite{AJOS:06}.

\subsubsection{Infima in $\wt{\R}$}
Let $A\subseteq\wt\R$. As in any partially ordered set, $\delta\in\wt\R$ is a lower bound for $A$ iff $\delta\le a$, for each $a\in A$; the infimum of $A$, denoted by $\inf A$, if it exists, is the greatest lower bound for $A$. 
As the set of lower bounds of $\overline A$ is equal to the set of lower bounds of $A$, $\inf A$ exists iff $\inf\overline A$ exists and in that case, $\inf A=\inf\overline A$. The following proposition gives a characterization of the infimum.
\begin{proposition}
Let $A\subseteq\wt\R$. Let $\delta\in\wt\R$ be a lower bound for $A$. The following are equivalent.
\begin{itemize}
\item[(i)] $\delta=\inf A$;
\item[(ii)] $\forall m\in\N$\, $\forall S\subseteq(0,1]$ with $e_S\ne 0$\, $\exists a\in A$\, $ae_S\ngeq(\delta + \caninf^m)e_S$;
\item[(iii)] $\forall m\in\N$\, $\forall S\subseteq(0,1]$ with $e_S\ne 0$\, $\exists T\subseteq S$ with $e_T\ne 0$\, $\exists a\in A$\, $ae_T\le(\delta + \caninf^m)e_T$.
\end{itemize}
\end{proposition}
\begin{proof}
$(i)\Rightarrow(ii)$: suppose there exists $m\in\N$ and $S\subseteq(0,1]$ with $e_S\ne 0$ such that for each $a\in A$, $ae_S\ge (\delta + \caninf^m)e_S$. Then also $a\ge a e_{\co S} + (\delta + \caninf^m)e_S\ge \delta + \caninf^m e_S$. So $\delta + \caninf^m e_S$ is a lower bound for $A$. As $e_S\ne 0$, $\delta\ne\inf A$.\\
$(ii)\Rightarrow(iii)$: if $ae_S\ngeq(\delta + \caninf^m)e_S$, then there exists $T\subseteq S$ with $e_T\ne 0$ such that $ae_T\le (\delta + \caninf^m)e_T$.\\
$(iii)\Rightarrow(i)$: let $\rho\in\wt\R$ be a lower bound for $A$. Suppose $\rho\nleq\delta$. Then there exists $S\subseteq(0,1]$ with $e_S\ne 0$ and $m\in\N$ such that $\rho e_S\ge(\delta+\caninf^m)e_S$. By hypothesis, there exists $T\subseteq S$ with $e_T\ne 0$ and $a\in A$ such that $a e_T\le (\delta + \caninf^{m+1})e_T$. Hence $(\delta+\caninf^m)e_T \le\rho e_T\le (\delta + \caninf^{m+1})e_T$, which contradicts the fact that $e_T\ne 0$.
\end{proof}

The infimum of $A$ is called \emph{close} if $\inf A\in\overline A$. In this case we use the notation $\overline{\inf} A$. Unlike in $\R$, an infimum in $\wt\R$ is not automatically close.
\begin{example}
Let $T\subseteq(0,1]$ with $e_T\ne 0$ and $e_{\co T}\ne 0$ and let $A=\{e_T + \caninf^m e_{\co T}:m\in\N\}\cup\{e_{\co T} + \caninf^m e_T: m\in\N\}$. Clearly, $0$ is a lower bound for $A$. Let $\delta\in\wt\R$ be a lower bound for $A$. Then $\delta\le \lim_{n\to\infty} (e_T+\caninf^n e_{\co T})= e_T$, hence $\delta e_{\co T}\le 0$ and similarly, $\delta e_{T}\le 0$. So $\delta=\delta e_{\co T}+\delta e_T\le 0$ and $\inf A=0$. On the other hand, $\abs{e_T + \caninf^m e_{\co T}}_\esp=\abs{e_{\co T} + \caninf^m e_{T}}_\esp=1$, for each $m\in\N$. Hence $0\notin \overline A$.
\end{example}
The close infimum can be easily characterized as follows.
\begin{proposition}
Let $A\subseteq\wt\R$. Let $\delta\in\wt\R$ be a lower bound for $A$. Then $\delta$ is a close infimum iff
\[
\forall m\in\N\ \exists a\in A\qquad a\le \delta +\caninf^m.
\]
\end{proposition}
Clearly, if $A$ reaches a minimum, then $\inf A=\min A$ and the infimum is close.

\subsection{$\wt{\C}$-modules with $\wt{\R}$-seminorms}
We introduce the notion of $\wt{\R}$-seminorm on a $\wt{\C}$-module $\G$. This determines a special kind of topological $\wt{\C}$-modules: the $\wt{\C}$-modules with $\wt{\R}$-seminorms.
\begin{definition}
\label{def_wtr_seminorm}
Let $\G$ be a $\wt{\C}$-module. An $\wt{\R}$-seminorm on $\G$ is a map $p:\G\to\wt{\R}$ such that
\begin{itemize}
\item[(i)] $p(0)=0$ and $p(u)\ge 0$ for all $u\in\G$;
\item[(ii)] $p(\lambda u)=|\lambda|p(u)$ for all $\lambda\in\wt{\C}$ and for all $u\in\G$;
\item[(iii)] $p(u+v)\le p(u)+p(v)$ for all $u,v\in\G$.
\end{itemize}
An $\wt{\R}$-seminorm $p$ such that $p(u)=0$ if and only if $u=0$ is called $\wt{\R}$-norm.
\end{definition}

From the properties which define an $\wt{\R}$-seminorm we easily see that the coarsest topology which makes a family $\{p_i\}_{i\in I}$ of $\wt{\R}$-seminorms on $\G$ continuous equips $\G$ with the structure of a topological $\wt{\C}$-module. Hence, any $\wt{\C}$-module with $\wt{\R}$-seminorms is a topological $\wt{\C}$-module. More precisely we have the following result.

\begin{proposition}
\label{prop_top_wtr_semin}
Any $\wt{\R}$-seminorm $p$ on $\G$ generates an ultra-pseudo-seminorm $\mP$ by setting $\mP(u):=|p(u)|_\esp=\esp^{-\val(p(u))}$. The $\wt{\C}$-linear topology on $\G$ determined by the family of $\wt{\R}$-seminorms $\{p_i\}_{i\in I}$ coincides with the topology of the corresponding ultra-pseudo-seminorms $\{\mP_i\}_{i\in I}$.
\end{proposition}
\begin{proof}
The fact that $\mP$ is an ultra-pseudo-seminorm follows from the properties of $p$ combined with the defining conditions of the ultra-pseudo-norm $|\cdot|_\esp$ of $\wt{\R}$. The families $\{p_i\}_{i\in I}$ and $\{\mP_i\}_{i\in I}$ generate the same topology on $\G$ since for all $\eta>0$, $\delta>0$ and $u\in\G$ we have that
\[
\{u\in\G:\ p_i(u)\le[(\eps^{-\log\eta})_\eps]\}\subseteq\{u\in\G:\ \mP_i(u)\le \eta\}\subseteq \{u\in\G:\ p_i(u)\le [(\eps^{-\log\eta-\delta})_\eps]\}.
\]
\end{proof}
In the particular case of $\G=\G_E$, where $(E,\{p_i\}_{i\in I})$ is a locally convex topological vector space, one can extend any seminorm $p_i$ to an $\wt{\R}$-seminorm on $\G_E$. This is due to the fact that if $(u_\eps)_\eps\in\M_E$ then $(p_i(u_\eps))_\eps\in\EM$ and if $(u_\eps-u_\eps')_\eps\in\Neg_E$ then $|p_i(u_\eps)-p_i(u'_\eps)|\le p_i(u_\eps-u'_\eps)=O(\eps^q)$ for all $q\in\N$. Proposition \ref{prop_top_wtr_semin} says that the sharp topology on $\G_E$ can be regarded as the topology of the $\wt{\R}$-seminorms $p_i(u):=[(p_i(u_\eps))_\eps]$ as well as the topology of the ultra-pseudo-seminorms $\mP_i(u)=|p_i(u)|_\esp$.

\begin{proposition}
\label{prop_cont_bilin}
Let $(\G,\{p_i\}_{i\in I})$, $(\F,\{q_j\}_{j\in J})$ and $(\mH,\{r_k\}_{k\in K})$ be topological $\wt{\C}$-modules with $\wt{\R}$-seminorms.
\begin{itemize}
\item[(i)] A $\wt{\C}$-linear map $T:\G\to\F$ is continuous if and only if the following assertion holds: for all $j\in J$, there exist a finite subset $I_0$ of $I$ and a constant $C\in\wt{\R}$ such that
\beq
\label{ineq_T}
q_j(Tu)\le C\sum_{i\in I_0}p_i(u) 
\eeq
for all $u\in\G$.
\item[(ii)] A $\wt{\C}$-sesquilinear map $a$ from $\G\times\F$ to $\mH$ is continuous if and only if for all $k\in K$, there exist finite subsets $I_0$ and $J_0$ of $I$ and $J$ respectively and a constant $C\in\wt{\R}$ such that 
\beq
\label{ineq_a}
r_k(a(u,v))\le C \sum_{i\in I_0}p_i(u) \sum_{j\in J_0}q_j(v),
\eeq
for all $u\in \G$ and $v\in \F$.
\end{itemize}
\end{proposition}
\begin{proof}
If the inequality \eqref{ineq_T} holds then the $\wt{\C}$-linear map $T$ is continuous, since from \eqref{ineq_T} we have that 
\[
\mathcal{Q}_j(Tu)\le |C|_\esp\max_{i\in I_0}\mP_i(u).
\]
This characterizes the continuity of $T$ as proved by Corollary 1.17 in \cite{Garetto:05a}. Assume now that $T$ is continuous at $0$. Hence, for all $j\in J$ and for all $c\in\R$ there exist $b\in\R$ and a finite subset $I_0$ of $I$  such that $q_j(Tu)\le [(\eps^c)_\eps]$ if $\sum_{i\in I_0}p_i(u)\le [(\eps^b)_\eps]$. Let $q\in \N$. For any $u\in\G$ we have that $[(\eps^b)_\eps]{u}/{(\sum_{i\in I_0}p_i(u)+[(\eps^q)_\eps])}$ belongs to the set of all $v\in\G$ such that $\sum_{i\in I_0}p_i(v)\le [(\eps^b)_\eps]$. Thus, 
\[
q_j(Tu)\le [(\eps^{c-b})_\eps]\biggl(\sum_{i\in I_0}p_i(u)+[(\eps^q)_\eps]\biggr).
\]
Letting $q$ go to $\infty$ we conclude that \eqref{ineq_T} is valid for $C=[(\eps^{c-b})_\eps]$. 

The proof of the second assertion of the proposition is similar and therefore left to the reader.
\end{proof}

We consider now the framework of Colombeau spaces of generalized functions based on a normed space and we provide a characterization for continuous $\wt{\C}$-linear maps given by a representative. We recall that a representative $(T_\eps)_\eps$ of a $\wt{\C}$-linear map $T:\G_E\to \G_F$, if it exists, is a net of linear maps from $E$ to $F$ such that $(T_\eps u_\eps)_\eps\in\M_F$ for all $(u_\eps)_\eps\in\M_E$, $(T_\eps u_\eps)_\eps\in\Neg_F$ for all $(u_\eps)_\eps\in\Neg_E$ and $Tu=[(T_\eps u_\eps)_\eps]$ for all $u=[(u_\eps)_\eps]\in\G_E$.
\begin{proposition}
Let $E$, $F$ be normed spaces and let $(T_\eps)_\eps$ be a net of linear maps from $E$ to $F$ such that $(T_\eps u_\eps)_\eps\in\M_F$ for each $(u_\eps)_\eps\in \Neg_E$. Then $(T_\eps u_\eps)_\eps\in\M_F$ for each $(u_\eps)_\eps\in \M_E$ and $(T_\eps u_\eps)_\eps\in\Neg_F$ for each $(u_\eps)_\eps\in \Neg_E$.
\end{proposition}
\begin{proof}
Let $(u_\eps)_\eps\in \M_E$, i.e., there exists $N\in\N$ such that $\norm{u_\eps}\le\eps^{-N}$, for sufficiently small $\eps$. Suppose that $(T_\eps u_\eps)_\eps\notin\M_F$. Then we can find a decreasing sequence $(\eps_n)_{n\in\N}$ with $\lim_n \eps_n=0$ such that $\norm{T_{\eps_n} u_{\eps_n}}\ge\eps_n^{-n}$. Let $v_{\eps_n}=u_{\eps_n}\eps_n^{n/2}$, $\forall n\in\N$ and let $v_\eps=0$, if $\eps\notin\{\eps_n:n\in\N\}$. Then for any $M\in\N$, $\norm{v_\eps}\le\norm{u_\eps \eps^M}\le\eps^{M-N}$, for sufficiently small $\eps$, but for each $n\in\N$, $\norm{T_{\eps_n}(v_{\eps_n})}=\eps_n^{n/2}\norm{T_{\eps_n}(u_{\eps_n})} \ge \eps_n^{-n/2}$. Hence $v_\eps\in \Neg_E$, but $T_\eps(v_\eps)\notin\M_E$, which contradicts the hypotheses.\\
Similarly, let $(u_\eps)_\eps\in \Neg_E$, i.e., for each $m\in\N$, $\norm{u_\eps}\le\eps^m$ as soon as $\eps\le\eta_m\in(0,1]$. Suppose that $(T_\eps u_\eps)_\eps\notin\Neg_F$. Then we can find $m\in\N$ and a decreasing sequence $(\eps_n)_{n\in\N}$ with $\lim_n \eps_n=0$, such that $\eps_n\le\eta_n$ and $\norm{T_{\eps_n} (u_{\eps_n})}\ge\eps_n^{m}$, $\forall n\in\N$. Let $v_{\eps_n}=u_{\eps_n}\eps_n^{-n/2}$, $\forall n\in\N$ and let $v_\eps=0$, if $\eps\notin\{\eps_n:n\in\N\}$. Then for each $n\in\N$, $\norm{v_{\eps_n}}\le\norm{u_{\eps_n}}\eps_n^{-n/2}\le\eps_n^{n/2}$, but $\norm{T_{\eps_n}(v_{\eps_n})} =\norm{T_{\eps_n}(u_{\eps_n})}\eps_n^{-n/2}\ge\eps_n^{m-n/2}$. Hence $v_\eps\in\Neg_E$, but $T_\eps(v_\eps)\notin\M_E$, which contradicts the hypotheses.
\end{proof}
Inspired by a similar result in \cite{OV:07} we obtain the following proposition.
\begin{proposition}
\label{prop_repres}
Let $E$ and $F$ be normed spaces and $T:\G_E\to\G_F$ a $\wt{\C}$-linear map. If $T$ has a representative $(T_\eps)_\eps$ then it is continuous.
\end{proposition}
\begin{proof}
We prove that if $(u_\eps)_\eps\in\Neg_E$ implies $(T_\eps u_\eps)_\eps\in\Neg_F$ then 
\beq
\label{impl_ind_repr}
\forall n\in\N\, \exists m\in\N\, \exists\eps_0\in(0,1]\, \forall\eps\in(0,\eps_0)\, \forall u\in E\qquad \Vert u\Vert_E\le \eps^m\ \Rightarrow\ \Vert T_\eps u\Vert_F\le \eps^n.
\eeq
Indeed, if we negate \eqref{impl_ind_repr} then we can find some $n'\in\N$, a decreasing sequence $\eps_m$ converging to $0$ and some $u_{\eps_m}\in E$ with $\Vert u_{\eps_m}\Vert_E\le\eps_m^m$ such that $\Vert T_{\eps_m}u_{\eps_m}\Vert_F>\eps_m^{n'}$. Let now $u_\eps=u_{\eps_m}$ for $\eps\in[\eps_m,\eps_{m-1})$ and $u_\eps=0$ for $\eps\in[\eps_0,1]$. By construction we have that $(u_\eps)_\eps\in\Neg_E$ and $\Vert T_{\eps_m} u_{\eps_m}\Vert_F>\eps_m^{n'}$ for all $m$. This is in contradiction with $(T_\eps u_\eps)_\eps\in\Neg_F$.\\
The assertion \eqref{impl_ind_repr} says that for all $n\in\N$ there exists a neighborhood $U=\{u\in\G_E:\, \Vert u\Vert_E\le[(\eps^m)_\eps]\}$ of $0$ which has image $T(U)$ contained in the neighborhood $V=\{v\in\G_F:\, \Vert v\Vert_F\le [(\eps^n)_\eps]\}$. Hence, the map $T$ is continuous at $0$ and thus continuous from $\G_E$ to $\G_F$.
\end{proof}

\begin{proposition}
\label{prop_cont_basic}
Let $E$ and $F$ be normed spaces, $T$ be a continuous $\wt{\C}$-linear map from $\G_E$ to $\G_F$ with a representative and $C\ge 0$ in $\wt{\R}$. Then, the following assertions are equivalent:
\begin{itemize}
\item[(i)] $\Vert Tu\Vert_F\le C\Vert u\Vert_E$ for all $u\in\G_E$;
\item[(ii)] for all representatives $(T_\eps)_\eps$ of $T$, for all representatives $(C_\eps)_\eps$ of $C$ and for all $q\in\N$ there exists $\eta\in(0,1]$ such that 
\beq
\label{ineq_T_eps}
\Vert T_\eps u\Vert_F\le (C_\eps +\eps^q)\Vert u\Vert_E
\eeq
for all $u\in E$ and $\eps\in(0,\eta]$;
\item[(iii)] for all representatives $(T_\eps)_\eps$ of $T$ there exists a representative $(C_\eps)_\eps$ of $C$ and $\eta\in(0,1]$ such that 
\beq
\label{ineq_T_eps_2}
\Vert T_\eps u\Vert_F\le C_\eps\Vert u\Vert_E
\eeq
for all $u\in E$ and $\eps\in(0,\eta]$. 
\end{itemize}
\end{proposition}
\begin{proof}
From Proposition \ref{prop_cont_bilin} we have that the continuity of $T$ is equivalent to $(i)$. In order to prove that $(i)$ implies $(ii)$, we begin by observing that $(i)$ is equivalent to claim that $e_S \Vert Tu\Vert_F\le C e_S\Vert u\Vert_E$ for all $S\subseteq(0,1]$. We want to prove that the negation of $(ii)$ implies that there exists a subset $S$ of $(0,1]$ and some $u\in\G_E$ such that $e_S\Vert Tu\Vert_F> C e_S\Vert u\Vert_E$. From 
\[
\exists (T_\eps)_\eps\, \exists (C_\eps)_\eps\, \exists q\in\N\, \forall\eta\in(0,1]\, \exists\eps\in(0,\eta]\, \exists u\in E\quad \Vert T_\eps u\Vert_F> (C_\eps+\eps^q)\Vert u\Vert_E 
\]
we have that there exists a decreasing sequence $(\eps_k)_k\subseteq(0,1]$ converging to $0$ and a sequence $(u_{\eps_k})_k$ of elements of $E$ with norm $1$ such that
\[
\Vert T_{\eps_k}(u_{\eps_k})\Vert_F>(C_{\eps_k}+\eps_k^q).
\]
Let us fix $x\in E$ with $\Vert x\Vert_E=1$. The net $u_\eps=u_{\eps_k}$ when $\eps=\eps_k$ and $u_\eps=x$ otherwise generates an element $u=[(u_\eps)_\eps]$ of $\G_E$ with $\wt{\R}$-norm $1$. Let now $S=\{\eps_k:\, k\in\N\}$. By construction we have that 
\[
e_S\Vert Tu\Vert_F=[(\chi_S T_{\eps_k}(u_{\eps_k}))_\eps]\ge e_S(C+[(\eps^q)_\eps])> e_S C. 
\]
This contradicts $(i)$. It is easy to prove that $(ii)$ implies $(iii)$. Indeed, by fixing representatives $(T_\eps)_\eps$ and $(C'_\eps)_\eps$ of $T$ and $C$ respectively, we can extract a decreasing sequence $(\eta_q)_{q\in\N}$ tending to $0$ such that $\Vert T_\eps(u)\Vert_F\le (C'_\eps+\eps^q)\Vert u\Vert_E$ for all $u\in E$ and $\eps\in(0,\eta_q]$. The net $n_\eps=\eps^q$ for $\eps\in(\eta^{q+1},\eta^q]$ is negligible and therefore $C_\eps=C'_\eps+n_\eps$ satisfies \eqref{ineq_T_eps_2} on the interval $(0,\eta_0]$. Finally, it is clear that $(iii)$ implies $(i)$.
\end{proof}
Note that from the previous propositions we have that if $T$ is given by a representative $(T_\eps)_\eps$, then it is a basic map.

\section{Hilbert $\wt{\C}$-modules}
\label{section_Hilbert}
\subsection{Definition}
This section is devoted to the definition and the first properties of the class of topological $\wt{\C}$-modules which are the mathematical core of the paper: the \emph{Hilbert $\wt{\C}$-modules}. In the intent of developing a topological and functional analytic theory of Hilbert $\wt{\C}$-modules, we start in Subsection \ref{subsec_projec} by investigating the notion of projection on suitable subsets of a Hilbert $\wt{\C}$-module $\G$. This requires the new concept of \emph{edged} subset of $\G$ and a formulation of convexity, which differently from the $\wt{\C}$-convexity introduced in \cite{Garetto:05a}, resembles the well-known classical definition for subsets of a vector space.

\begin{definition}
\label{def_scalar}
Let $\G$ be a $\wt{\C}$-module. A scalar product $\inner{\cdot}{\cdot}$ is a $\wt{\C}$-sesquilinear form from $\G\times\G$ to $\wt{\C}$ satisfying the following properties:
\begin{itemize}
\item[(i)] $\inner{u}{v}=\overline{\inner{v}{u}}$ for all $u,v\in\G$,
\item[(ii)] $\inner{u}{u}\in\wt{\R}$ and $\inner{u}{u}\ge 0$ for all $u\in\G$,
\item[(iii)] $\inner{u}{u}= 0$ if and only if $u=0$. 
\end{itemize} 
\end{definition}
In the sequel we denote $\sqrt{\inner{u}{u}}$ by $\Vert u\Vert$.

Since $\wt\C$ is not a field, the following proposition is not immediate.
\begin{proposition}
\label{prop_Cauchy}
Let $\G$ be a $\wt{\C}$-module with scalar product $\inner{\cdot}{\cdot}$. Then for all $u,v\in\G$ the Cauchy-Schwarz inequality holds:
\beq
\label{cs_ineq}
|\inner{u}{v}|\le \Vert u\Vert\Vert{v}\Vert.
\eeq
\end{proposition}
\begin{proof}
Let $\alpha\in\wt{\C}$. By definition of a scalar product we know that $\Vert{u+\alpha v}\Vert$ is a positive generalized real number. Hence, the $\wt{\C}$-sesquilinearity of $\inner{\cdot}{\cdot}$ yields
\beq
\label{ineq_alpha}
0\le\Vert{u+\alpha v}\Vert^2=\Vert{u}\Vert^2+\alpha\overline{\inner{u}{v}}+\overline{\alpha}\inner{u}{v}+|\alpha|^2\Vert v\Vert^2.
\eeq
We will derive the Cauchy-Schwarz inequality \eqref{cs_ineq} from \eqref{ineq_alpha} by choosing a suitable sequence of $\alpha\in\wt{\C}$. In detail, let $\alpha_n:=-\inner{u}{v}/(\Vert{v}\Vert^2+[(\eps^n)])$. The equality \eqref{ineq_alpha} combined with $\Vert{v}\Vert^2\le\Vert{v}\Vert^2+[(\eps^n)]$ yields
\begin{multline*}
0\le\Vert{u}\Vert^2-\frac{|\inner{u}{v}|^2}{\Vert{v}\Vert^2+[(\eps^n)]}-\frac{|\inner{u}{v}|^2}{\Vert{v}\Vert^2+[(\eps^n)]}+\frac{|\inner{u}{v}|^2}{(\Vert{v}\Vert^2+[(\eps^n)])^2}\Vert{v}\Vert^2\\
\le \Vert{u}\Vert^2-2\frac{|\inner{u}{v}|^2}{\Vert{v}\Vert^2+[(\eps^n)]}+\frac{|\inner{u}{v}|^2}{\Vert{v}\Vert^2+[(\eps^n)]}.
\end{multline*}
Hence, 
\[
0\le\Vert{u}\Vert^2(\Vert{v}\Vert^2+[(\eps^n)])-|\inner{u}{v}|^2
\]
for all $n$, and since the sequence $(\Vert{v}\Vert^2+[(\eps^n)])_n$ tends to $\Vert{v}\Vert^2$ in $\wt{\R}$ it follows that 
the Cauchy-Schwarz inequality \eqref{cs_ineq} holds. 
\end{proof}
We use the Cauchy-Schwarz inequality in proving the following proposition.
\begin{proposition}
\label{prop_norm}
The map $\Vert\cdot\Vert:\G\to\wt{\R}:u\to \Vert u\Vert:=\inner{u}{u}^{\frac{1}{2}}$ is an $\wt{\R}$-norm on $\G$ and the map $\mathcal{P}:\G\to[0,+\infty):u\to |\inner{u}{u}^{\frac{1}{2}}|_\esp={|\inner{u}{u}|_\esp}^{\frac{1}{2}}$ is an ultra-pseudo-norm on $\G$.
\end{proposition}
\begin{proof}
The third property of Definition \ref{def_scalar} ensures that $\Vert u\Vert=0$ if and only if $u=0$. Let us now take $\lambda\in\wt{\C}$. From the homogeneity of the scalar product we have that
\[
\Vert\lambda u\Vert=\inner{\lambda u}{\lambda u}^{\frac{1}{2}}=\big(|\lambda|^2\Vert u\Vert^2)^{\frac{1}{2}}=|\lambda|\Vert u\Vert.
\]
Finally, we write $\Vert u+v\Vert^2$ as $\Vert u\Vert^2+2\Re\inner{u}{v}+\Vert v\Vert^2$ and since $\Re\inner{u}{v}\le|\inner{u}{v}|$ we obtain from the Cauchy-Schwarz inequality \eqref{cs_ineq} that 
\[
\Vert u+v\Vert^2\le \Vert u\Vert^2+\Vert v\Vert^2+2|\inner{u}{v}|\le (\Vert u\Vert+\Vert v\Vert)^2.
\]
It follows that $\Vert u+v\Vert\le \Vert u\Vert +\Vert v\Vert$ for all $u,v\in\G$. Thus, $\Vert\cdot\Vert$ is an $\wt{\R}$-norm on $\G$. Proposition \ref{prop_top_wtr_semin} combined with the fact that $|\lambda^{\frac{1}{2}}|_\esp=|\lambda|^{\frac{1}{2}}_\esp$ allows us to conclude that $\mP$ is an ultra-pseudo-norm.
\end{proof}
From Proposition \ref{prop_top_wtr_semin} we have that a $\wt{\C}$-module $\G$ with scalar product $\inner{\cdot}{\cdot}$ can be endowed with the topology of the $\wt{\R}$-norm $\Vert\cdot\Vert$ generated by $\inner{\cdot}{\cdot}$ or equivalently with the topology of the ultra-pseudo-norm $\mP(u)=|\inner{u}{u}|^{\frac{1}{2}}_\esp$. This means that any $\wt{\C}$-module with a scalar product is a $\wt{\C}$-module with an $\wt{\R}$-norm and hence a topological $\wt{\C}$-module. Proposition \ref{prop_Cauchy} combined with Proposition \ref{prop_norm} yields the following continuity result.
\begin{proposition}
\label{prop_cont_scal}
Let $\G$ be a $\wt{\C}$-module with scalar product $\inner{\cdot}{\cdot}$, topologized through the ultra-pseudo-norm $\mP(u)={|\inner{u}{u}|_\esp}^{\frac{1}{2}}$. The scalar product is a continuous $\wt{\C}$-sesquilinear map from $\G\times\G$ to $\wt{\C}$.
\end{proposition}
 
\begin{definition}
\label{def_Hilbert}
A Hilbert $\wt{\C}$-module is a $\wt{\C}$-module with scalar product $\inner{\cdot}{\cdot}$ which is complete when endowed with the topology of the corresponding ultra-pseudo-norm $\mP$.
\end{definition}
Since a closed subset of a complete topological $\wt{\C}$-module is complete, we have that a closed $\wt{\C}$-submodule of a Hilbert $\wt{\C}$-module is itself a Hilbert $\wt{\C}$-module.

\begin{example}
\label{ex_Hilbert}
\leavevmode
\begin{trivlist}
\item[(i)] A first example of a Hilbert $\wt{\C}$-module is given by $\G_H$, where $(H,\inner{\cdot}{\cdot})$ is a Hilbert space. The scalar product on $\G_H$ is obtained by letting $\inner{\cdot}{\cdot}$ act componentwise on the representatives of the generalized functions in $\G_H$ as follows: $\inner{u}{v}=[(\inner{u_\eps}{v_\eps})_\eps]$. By Proposition 3.4 in \cite{Garetto:05a} one can omit the assumption of completeness on $H$ and still obtain that $\G_H$ is complete with respect to the sharp topology induced by the scalar product.
\item[(ii)] The topological structure on $\G_H$ determined by the scalar product of $H$ can be equivalently generated by any continuous $\wt{\C}$-sesquilinear form $a$ on $\G_H\times\G_H$ such that $a(u,v)=\overline{a(v,u)}$ for all $u,v\in\G_H$, $a(u,u)\ge 0$ for all $u\in\G_H$ and the following bound from below holds:
\beq
\label{coerc_ex}
\exists C\in\wt{\R}, C\ge 0,\, \text{invertible},\ \forall u\in\G_H\qquad a(u,u)\ge C\Vert u\Vert^2,
\eeq
(see also Definition \ref{def_coercivity}). Since $a$ satisfies the conditions of Definition \ref{def_scalar}, it is a scalar product on $\G_H$ and the corresponding Cauchy-Schwarz inequality is valid. Hence, $\Vert u\Vert_a:=a(u,u)^{\frac{1}{2}}$ is an $\wt{\R}$-norm. Combining the continuity of $a$ with the estimate \eqref{coerc_ex} we have that $\Vert \cdot\Vert_a$ is equivalent to the usual $\wt{\R}$-norm $\Vert\cdot\Vert$. this means that there exist $C_1,C_2\ge 0$ real generalized numbers such that  
\[
C_1\Vert u\Vert \le \Vert u\Vert_a\le C_2\Vert u\Vert
\]
for all $u\in\G_H$.
\end{trivlist}
\end{example}
A further example of a Hilbert $\wt{\C}$-module is provided by the following proposition
\begin{proposition}
\label{prop_ex_3}
Let $(H_\eps,\inner{\cdot}{\cdot}_{H_\eps})_\eps$ be a net of vector spaces with scalar product and let $\G$ be the $\wt{\C}$-module obtained by factorizing the set
\[
\M_{(H_\eps)_\eps}=\{(u_\eps)_\eps:\, \forall\eps\in(0,1]\, u_\eps\in H_\eps\ \text{and}\ \exists N\in\N\, \Vert u\Vert_{H_\eps}=O(\eps^{-N})\}
\]
of moderate nets with respect to the set
\[
\Neg_{(H_\eps)_\eps}=\{(u_\eps)_\eps:\, \forall\eps\in(0,1]\, u_\eps\in H_\eps\ \text{and}\ \forall q\in\N\, \Vert u\Vert_{H_\eps}=O(\eps^{q})\}
\]
of negligible nets. Let $\inner{\cdot}{\cdot}:\G\times\G\to\wt{\C}$ be the $\wt{\C}$-sesquilinear form defined as follows:
\beq
\label{def_scalar_3}
\inner{u}{v}=[(\inner{u_\eps}{v_\eps}_{H_\eps})_\eps].
\eeq
Then, $\inner{\cdot}{\cdot}$ is a scalar product on $\G$ which equips $\G$ with the structure of a Hilbert $\wt{\C}$-module.
\end{proposition}
\begin{proof}
Applying the Cauchy-Schwarz inequality componentwise in any Hilbert space $H_\eps$ we have that \eqref{def_scalar_3} is a well-defined $\wt{\C}$-sesquilinear form on $\G\times\G$ such that the properties $(i),(ii),(iii)$ of Definition \ref{def_scalar} are fulfilled. Let $\G$ be endowed with the topology of this scalar product, i.e., with the topology of the $\wt{\R}$-norm $\Vert u\Vert=[(\Vert u_\eps\Vert_{H_\eps})_\eps]$. We want to prove that any Cauchy sequence in $\G$ is convergent. If $(u_n)_n$ is a Cauchy sequence then we can extract a subsequence $(u_{n_k})_k$ and a corresponding subsequence $((u_{n_k,\eps})_\eps)_k$ of representatives such that $\Vert u_{n_{k+1},\eps}-u_{n_k,\eps}\Vert_{H_\eps}\le \eps^k$ for all $\eps\in(0,\eps_k)$, with $\eps_k\searrow 0$, $\eps_k\le 2^{-k}$ for all $k\in\N$. Arguing as in the proof of \cite[Proposition 3.4]{Garetto:05a} we set $h_{k,\eps}=u_{n_{k+1},\eps}-u_{n_k,\eps}$ for $\eps\in(0,\eps_k)$ and $h_{k,\eps}=0$ for $\eps\in[\eps_k,1]$. Obviously, $(h_{k,\eps})_\eps\in\M_{(H_\eps)_\eps}$ and $\Vert h_{k,\eps}\Vert_{H_\eps}\le\eps^k$ on the whole interval $(0,1]$. Let now
\[
u_\eps:=\sum_{k=0}^\infty h_{k,\eps}+u_{n_0,\eps}.
\]
This sum is locally finite and moderate, since
\[
\Vert u_\eps\Vert_{H_\eps}\le \sum_{k=0}^\infty \Vert h_{k,\eps}\Vert_{H_\eps}+\Vert u_{n_0,\eps}\Vert_{H_\eps}\le \sum_{k=0}^\infty \eps_k^k+\Vert u_{n_0,\eps}\Vert_{H_\eps}\le\sum_{k=0}^\infty 2^{-k}+\Vert u_{n_0,\eps}\Vert_{H_\eps}.
\]
Hence, $(u_\eps)_\eps$ generates an element of $\G$. By construction the sequence $(u_{n_k})_k$ converges to $u$. Indeed, for all $\overline{k}\ge 1$ we have that
\[
\Vert u_{n_{\overline{k}},\eps}-u_\eps\Vert_{H_\eps}=\biggl\Vert u_{n_{\overline{k}},\eps}-u_{n_0,\eps}-\sum_{k=0}^\infty h_{k,\eps}\biggr\Vert_{H_\eps}=\biggl\Vert -\sum_{k=\overline{k}}^\infty h_{k,\eps}\biggr\Vert_{H_\eps}\le  \eps^{\overline{k}-1}\sum_{k=\overline{k}}^\infty\eps_k\le \eps^{\overline{k}-1}\sum_{k=\overline{k}}^\infty 2^{-k}
\] 
and the proof is complete.
\end{proof}

\begin{proposition}
The Hilbert $\wt\C$-module $\G=\M_{(H_\eps)_\eps}/\Neg_{(H_\eps)_\eps}$ defined in the previous proposition is (algebraically and isometrically) isomorphic with an internal submodule of a Hilbert $\wt\C$-module $\G_H$ for some pre-Hilbert space $H$.
\end{proposition}
\begin{proof}
Let $H=\bigoplus_{\lambda\in(0,1]}H_\lambda$ be the direct sum of the pre-Hilbert spaces $H_\lambda$, which is by definition the set of all nets $(u_\lambda)_{\lambda\in(0,1]}$, where $u_\lambda\in H_\lambda$, for each $\lambda$, which satisfy $\sum_{\lambda\in(0,1]}\norm{u_\lambda}_{H_\lambda}^2<+\infty$. This is a pre-Hilbert space \cite[Section 2.6]{KadRin:83} for the componentwise algebraic operations and the inner product
\[\inner{(u_\lambda)_\lambda}{(v_\lambda)_\lambda}= \sum_{\lambda\in(0,1]} \inner{u_\lambda}{v_\lambda}_{H_\lambda}.\] 
(When all $H_\lambda$ are Hilbert spaces, the direct sum is actually a Hilbert space (\cite[Section2.6]{KadRin:83}). Each $H_\lambda$ is canonically (algebraically and isometrically) isomorphic with a submodule $\wt H_\lambda$ of $H$ by the embedding $\iota_\lambda$: $H_\lambda\to H$: $\iota_\lambda(u)=(u_\mu)_\mu$ with $u_\lambda=u$, $u_\mu=0$ if $\mu\ne\lambda$. Hence we can consider the internal subset $[(\wt H_\eps)_\eps]\subseteq \G_H$.
Let now $\iota$: $\G\to \G_H$ be defined on representatives by $\iota((u_\eps)_\eps)=(\iota_\eps(u_\eps))_\eps$. Since $\norm{u_\eps}_{H_\eps}=\norm{\iota_\eps(u_\eps)}_H$, for each $\eps$, $(\iota_\eps(u_\eps))_\eps$ belongs to $\M_H$, resp.\ $\Neg_H$ iff $(u_\eps)_\eps$ belongs to $\M_{(H_\eps)_\eps}$, resp.\ $\Neg_{(H_\eps)_\eps}$. Hence $\iota$ is well-defined and injective.
Clearly, the image of $\iota$ is contained in $[(\wt H_\eps)_\eps]$. Conversely, each $v\in [(\wt H_\eps)_\eps]$ has a representative $(v_\eps)_\eps$ with each $v_\eps\in \wt H_\eps$. So $v_\eps=\iota_\eps(u_\eps)$, for some $u_\eps\in H_\eps$. Again by $\norm{u_\eps}_{H_\eps}=\norm{\iota_\eps(u_\eps)}_H$, the net $(u_\eps)_\eps$ belongs to $\M_{(H_\eps)_\eps}$, so it represents $u\in \G$ with $\iota(u)=v$.
\end{proof}
We see from the previous proposition that there is no loss of generality by considering the $\wt{\C}$-modules $\G_H$ instead of the factors $\M_{(H_\eps)_\eps}/\Neg_{(H_\eps)_\eps}$.

\subsection{Projection on a subset $C$}
\label{subsec_projec}

\begin{definition}
\label{def_edged}
Let $\G$ be a Hilbert $\wt{\C}$-module and $C$ a nonempty subset of $\G$. We say that $C$ is \emph{reachable} from $u\in\G$ if $$\overline{\inf}_{w\in C}\Vert{u-w}\Vert$$ exists in $\wt{\R}$. $C$ is called \emph{edged} if it is reachable from any $u\in\G$.
\end{definition}
From the definition it is clear that if $C$ is edged then $u+C$ is edged too for all $u\in\G$. Since $\overline{\inf}_{w\in\overline{C}}\Vert u-w\Vert=\overline{\inf}_{w\in{C}}\Vert u-w\Vert$ we have that $C$ is edged if and only if $\overline{C}$ is edged.
\begin{theorem}
\label{theo_projec_H}
Let $C$ be a closed nonempty subset of the Hilbert $\wt{\C}$-module $\G$ such that $C+C\subseteq 2C$. If $C$ is reachable from $u\in\G$ then there exists a unique $v\in C$ such that
\[
\Vert{u-v}\Vert=\inf_{w\in C}\Vert{u-w}\Vert.
\]
The element $v$ is called the projection of $u$ on $C$ and denoted by $P_C(u)$.
\end{theorem}
\begin{proof}
Note that when $\overline{\inf}_{w\in C}\Vert u-w\Vert$ exists in $\wt{\R}$ one has that $\overline{\inf}_{w\in C}\Vert u-w\Vert^2=(\overline{\inf}_{w\in C}\Vert u-w\Vert)^2$.
As the properties of $C$ are translation invariant, we can assume $u=0$. We set $\overline{\inf}_{w\in C}\Vert{w}\Vert^2=\delta$ in $\wt{\R}$. By definition of close infimum we can extract a sequence $w_n$ in $C$ such that $\Vert{w_n}\Vert^2\to \delta$. The fact that $C+C\subseteq 2C$ implies that $\frac{w_n+w_m}{2}$ belongs to $C$ for all $n,m\in\N$. So,  
\beq
\label{parallel}
0\le\Vert{w_n-w_m}\Vert^2= -4\Vert{\frac{w_n+w_m}{2}}\Vert^2+2\Vert{w_n}\Vert^2+2\Vert{w_m}\Vert^2\le -4\delta+2\Vert{w_n}\Vert^2+2\Vert{w_m}\Vert^2.
\eeq
From $\Vert{w_n}\Vert^2\to \delta$ it follows that $(w_n)_n$ is a Cauchy sequence in $C$ and therefore it is convergent in $\G$ to an element $v$ of $C$. By continuity of the $\wt{\R}$-norm we have that $\Vert v\Vert^2=\delta$. Finally, if we assume that there exists another $v'\in C$ such that $\Vert{v'}\Vert^2=\delta$, the inequality \eqref{parallel} is valid for $v-v'$ and proves that $v=v'$ in $\G$.
\end{proof}
\begin{corollary}
\label{theo_projec}
Let $C$ be a closed edged subset of the Hilbert $\wt{\C}$-module $\G$ such that $C+C\subseteq 2C$. Then, for all $u\in\G$ there exists a unique $v\in C$ such that
\[
\Vert{u-v}\Vert=\inf_{w\in C}\Vert{u-w}\Vert.
\]
\end{corollary}
The following example shows that the hypothesis of close infimum is necessary in the assumptions of the previous theorem.
\begin{example}
There exists a nonempty closed subset $C$ of $\wt\C$ with $\lambda C + (1-\lambda) C\subseteq C$, for each $\lambda\in\wt{[0,1]}:=\{x\in\wt{\R}:\, 0\le x\le 1\}$ for which $\inf_{c\in C} \abs{c}$ exists, but which is not reachable from $0\in\wt\C$.
\end{example}
\begin{proof}
Let for each $n\in\N$, $S_n\subseteq (0,1]$ with $e_{S_n}\ne 0$ and $S_n\cap S_m=\emptyset$ if $n\ne m$. Let ${\mathcal T}=\{T\subseteq(0,1]: e_T\ne 0 $ and $e_T e_{S_n}=0, \forall n\in\N\}\cup \{S_n: n\in\N\}$. Let $A=\{e_{\co T}: T\in{\mathcal T}\}$. We show that $\inf A=0$. Let $\rho\in\wt\R$, $\rho\le e_{\co T}$ for each $T\in{\mathcal T}$. Suppose that $\rho\nleq 0$. Then there exist $U\subseteq (0,1]$ with $e_U\ne 0$ and $m\in\N$ such that $\rho e_U\ge \caninf^m e_U$. Then also $\caninf^m e_Ue_{S_n}\le\rho e_U e_{S_n}\le e_{\co S_n}e_U e_{S_n}=0$, so $e_U e_{S_n}=0$, $\forall n$. Hence $U\in \mathcal T$, and $\caninf^me_U\le \rho e_U\le e_{\co U} e_U=0$, which contradicts $e_U\ne 0$.\\
Now let $B=\{\lambda_1 a_1 +\cdots + \lambda_m a_m: m\in\N, a_j\in A, \lambda_j\in\wt{[0,1]}, \sum_{j=1}^m\lambda_j=1\}$. Then also $\inf B=0$ and $\lambda B + (1-\lambda) B\subseteq B$, for each $\lambda\in\wt{[0,1]}$. We show that $0$ is not a close infimum for $B$.\\
Let $\lambda_1 a_1 +\cdots + \lambda_m a_m\in B$. Fix representatives $\lambda_{j,\eps}$ of $\lambda_j$ and let \[U_j=\{\eps\in(0,1]: \lambda_{j,\eps}=\max(\lambda_{1,\eps},\dots,\lambda_{m,\eps})\}.\] Then $e_{U_j}=\sum_{i=1}^m\lambda_i e_{U_j}\le m \lambda_j e_{U_j}$, for $j=1,\dots, m$. So $\lambda_1 a_1 + \cdots + \lambda_m a_m \ge \frac{1}{m}(e_{U_1}a_1+\cdots + e_{U_m}a_m)$ with $\bigcup U_j=(0,1]$. Let $a_j=e_{\co T_j}$, $T_j\in\mathcal T$. By the definition of $\mathcal T$, there exists $n\in\N$ such that $e_{T_1} e_{S_n}=\cdots=e_{T_m} e_{S_n}=0$. Then $\lambda_1 a_1 +\cdots + \lambda_m a_m\ge\frac{1}{m}(e_{U_1}e_{\co T_1}+\cdots + e_{U_m}e_{\co T_m})e_{S_n}\ge \frac{1}{m}e_{S_n}$. Hence $\abs{\lambda_1 a_1 +\cdots + \lambda_m a_m}_\esp\ge 1$. Consequently, $0\notin\overline B$.\\
Finally, let $C=\overline B$.
\end{proof}

Under the hypotheses of Corollary \ref{theo_projec} we can define the map $P_C$ as the map which assigns to each $u\in\G$ its projection on $C$. A careful investigation of the properties of the map $P_C$ requires the following lemma, which is obtained by observing the proof of Theorem \ref{theo_projec_H}.
\begin{lemma}
\label{lemma_cont}
Let $C$ be a closed nonempty subset of the Hilbert $\wt{\C}$-module $\G$ such that $C+C\subseteq 2C$, $u$ an element of $\G$ such that $C$ is reachable from $u$ and $(v_n)_n$ a sequence of elements of $C$. If $\Vert{u-v_n}\Vert\to\inf_{w\in C}\Vert{u-w}\Vert=\Vert{u-P_C(u)}\Vert$ in $\wt{\R}$ then $v_n\to P_C(u)$ in $\G$.
\end{lemma}
  
\begin{proposition}
\label{prop_cont}
Let $C$ be a closed edged subset of the Hilbert $\wt{\C}$-module $\G$ such that $C+C\subseteq 2C$. The operator $P_C$ has the following properties:
\begin{trivlist}
\item[(i)] $P_C(u)=u$ if and only if $u\in C$;
\item[(ii)] $P_C(\G)=C$;
\item[(iii)] $P_C^2=P_C$;
\item[(iv)] $P_C$ is a continuous operator on $\G$.
\end{trivlist}
\end{proposition} 
\begin{proof}
$(i)$ It is obvious that $u$ belongs to $C$ if it coincides with its projection. Conversely, if $u\in C$ then $\Vert{u-P_C(u)}\Vert=\inf_{w\in C}\Vert{u-w}\Vert=0$ and therefore $u=P_C(u)$. The assertion $(ii)$ is trivial and from $(i)$ it follows that the operator $P_C$ is idempotent. Let us now prove that $P_C$ is continuous. Since $\G$ is a metric space it is sufficient to prove that $P_C$ is sequentially continuous, i.e., $u_n\to u$ implies $P_C(u_n)\to P_C(u)$. This is guaranteed by Lemma \ref{lemma_cont} if we prove that the sequence $\Vert{u-P_C(u_n)}\Vert$ converges to $\Vert{u-P_C(u)}\Vert$ in $\wt{\R}$. The triangle inequality, valid in $\wt{\R}$ for $\Vert{\cdot}\Vert$, combined with the fact that $\Vert{u_n-P_C(u_n)}\Vert\le\Vert{u_n-P_C(u)}\Vert$ leads to 
\[
\Vert{u-P_C(u_n)}\Vert\le\Vert{u-u_n}\Vert+\Vert{u_n-P_C(u_n)}\Vert\le 2\Vert{u-u_n}\Vert+\Vert{u-P_C(u)}\Vert.
\]
It follows that
\[
0\le \Vert{u-P_C(u_n)}\Vert-\Vert{u-P_C(u)}\Vert\le 2\Vert{u-u_n}\Vert.
\]
Since $u_n\to u$, we conclude that $\Vert{u-P_C(u_n)}\Vert\to\Vert{u-P_C(u)}\Vert$ in $\wt{\R}$.
\end{proof}


\begin{proposition}
\label{prop_charac}
Let $C$ be a closed nonempty subset of the Hilbert $\wt{\C}$-module $\G$ such that $\lambda C+(1-\lambda)C\subseteq C$ for all real generalized numbers $\lambda\in\{[(\eps^q)_\eps]\}_{q\in\N}\cup\{\frac{1}{2}\}$. $C$ is reachable from $u\in\G$ if and only if there exists $v\in C$ such that  
\beq
\label{ineq_Re}
\Re\inner{u-v}{w-v}\le 0
\eeq
for all $w\in C$.  In this case $v=P_C(u)$.
\end{proposition}
\begin{proof}
We begin by assuming that $C$ is reachable from $u$. Then, $P_C(u)\in C$ and $\Vert{u-P_C(u)}\Vert^2=\inf_{w\in C}\Vert{u-w}\Vert^2$. Let $w\in C$. By the hypotheses on $C$ we know that $(1-[(\eps^q)_\eps])P_C(u)+[(\eps^q)_\eps]w$ belongs to $C$. Hence, 
\begin{multline*}
\Vert{u-P_C(u)}\Vert^2\le\Vert{u-P_C(u)-[(\eps^q)_\eps](w-P_C(u))}\Vert^2\\
\le \Vert{u-P_C(u)}\Vert^2-2[(\eps^q)_\eps]\Re\inner{u-P_C(u)}{w-P_C(u)}+[(\eps^q)_\eps]^2\Vert{w-P_C(u)}\Vert^2.
\end{multline*}
By the previous inequality and the invertibility of $[(\eps^q)_\eps]$ we obtain
\[
\Re\inner{u-P_C(u)}{w-P_C(u)}\le \frac{[(\eps^q)_\eps]}{2}\Vert{w-P_C(u)}\Vert^2.
\]
Letting $q$ tend to $\infty$ we conclude that $\Re\inner{u-P_C(u)}{w-P_C(u)}\le 0$.

Assume now that $v\in C$ and $\Re\inner{u-v}{w-v}\le 0$ for all $w\in C$. By the properties of a scalar product we can write
\[
\Vert{u-v}\Vert^2-\Vert{u-w}\Vert^2=\Vert{u-v}\Vert^2-\Vert{u-v}\Vert^2+2\Re\inner{u-v}{w-v}-\Vert{w-v}\Vert^2\le 0
\]
for all $w\in C$. This means that $\Vert{u-v}\Vert^2\le\inf_{w\in C}\Vert{u-w}\Vert^2$ and since $v\in C$ we conclude that $\Vert{u-v}\Vert^2=\min_{w\in C}\Vert{u-w}\Vert^2$. Thus, $v=P_C(u)$.
\end{proof}
In Proposition \ref{prop_new_convex} we will prove that under the assumptions of the previous theorem, in fact $\lambda C+(1-\lambda)C\subseteq C$ for all $\lambda\in \wt{[0,1]}=\{x\in\wt\R: 0\le x\le 1\}$. 
\begin{corollary}
\label{corol_1_new}
Let $M$ be a closed $\wt{\C}$-submodule of $\G$. Then, defining the $\wt{\C}$-submodule $M^\perp:=\{w\in\G:\, \forall v\in M\ \inner{w}{v}=0\}$, we have that $M$ is reachable from $u\in\G$ if and only if there exists $v\in M$ such that $u-v\in M^\perp$ if and only if $u\in\G$ can be uniquely written in the form $u=u_1+u_2$, where $u_1\in M$ and $u_2\in M^\perp$.
\end{corollary}
\begin{proof}
Since $M$ is a closed $\wt{\C}$-submodule, \eqref{ineq_Re} is equivalent to $\inner{u-v}{w}=0$ for all $w\in M$. Indeed, $\Re\inner{u-v}{w}=0$ and from $\inner{u-v}{-iw}=0$ we have $\Im\inner{u-v}{w}=0$. Hence, $u=v+(u-v)$, where $v\in M$ and $u-v\in M^\perp$ are uniquely determined by the scalar product on $\G$ since $M\cap M^{\perp}=\{0\}$.
\end{proof}

\begin{corollary}
\label{corol_1}
Let $M$ be a $\wt{\C}$-submodule of $\G$. Then,  
\begin{itemize}
\item[(i)] $M$ is closed and edged if and only if $\G=M\oplus M^\perp$, i.e., $\G=M+M^\perp$ and $M\cap M^\perp=\{0\}$.
\end{itemize}
If $M$ is closed and edged, the following holds:
\begin{itemize}
\item[(ii)]if $M^\perp=\{0\}$ then $M=\G$;
\item[(iii)] $M^{\perp\perp}=M$;
\item[(iv)] the projection $P_M$ is a $\wt{\C}$-linear operator on $\G$;
\item[(v)] $\inner{P_M(u)}{u}\ge 0$ for all $u\in\G$;
\item[(vi)] $M^\perp$ is closed and edged and $P_{M^\perp}(u)=u-P_M(u)$.
\end{itemize}
\end{corollary}
\begin{proof}
$(i)$ $\Rightarrow:$ is clear from Corollary \ref{corol_1_new}. Let us assume that $\G=M+ M^\perp$ and that $u\in\overline{M}$. It follows that $u=u_1+u_2$ with $u_1\in M$ and $u_2\in M^\perp$. So, $u-u_1=u_2\in \overline{M}\cap\overline{M}^\perp=\{0\}$. Hence, $u\in M$ and $M$ is closed. From Corollary \ref{corol_1_new} we have that $M$ is edged. 

$(ii)$ Assume now that $M^\perp=\{0\}$. Then $M$ has to coincide with $\G$. This follows from the fact that any $u\in \G\setminus M$ can be written as $u_1+u_2$, where $u_1\in M$ and $u_2\neq 0$ belongs to $M^\perp$. 

$(iii)$ By construction $M\subseteq M^{\perp\perp}$. From the first assertion of this proposition we know that any $u\in M^{\perp\perp}$ can be written as $u_1+u_2$, where $u_1\in M\subseteq M^{\perp\perp}$ and $u_2\in M^\perp$. Hence, $u_2=u-u_1\in M^{\perp\perp}$ and since $u_2\in M^\perp$ we obtain that $\inner{u_2}{u_2}=0$. It follows that $u\in M$.

$(iv)$ The $\wt{\C}$-linearity of the operator $P_M$ is due to the uniqueness of the decomposition $u_1+u_2=P_M(u_1)+P_M(u_2)+(u_1+u_2-P_M(u_1)-P_M(u_2))$, where $P_M(u_1)+P_M(u_2)\in M$ and $u_1+u_2-P_M(u_1)-P_M(u_2)\in M^\perp$. Analogously, for all $\lambda\in\wt{\C}$ one has that $\lambda P_M(u)\in M$,  $\lambda(u-P_M(u))\in M^\perp$ and $\lambda u=\lambda P_M(u)+\lambda(u-P_M(u))$.
 
$(v)$ We write $u$ as the sum  of $u-P_M(u)\in M^\perp$ and $P_M(u)\in M$. It follows that 
\[
\inner{P_M(u)}{u}=\inner{P_M(u)}{u-P_M(u)}+\inner{P_M(u)}{P_M(u)}=\Vert P_M(u)\Vert^2.
\]

$(vi)$ It is clear that if $M$ is a closed $\wt{\C}$-submodule then $M^\perp$ is a closed $\wt{\C}$-submodule too. We want to prove that $M^\perp$ is edged, i.e. it is reachable from every element of $\G$. By Corollary \ref{corol_1_new} we know that every element $u$ of $\G$ can be uniquely wriiten as $P_M(u)+(u-P_M(u))$, where $u-P_M(u)\in M^\perp$. By assertion $(iii)$ we have that $P_M(u)\in M=M^{\perp\perp}$. So, again by Corollary \ref{corol_1_new} we conclude that $M^\perp$ is reachable from $u$.
\end{proof}
\begin{remark}
In \cite{Vernaeve:07a} it is shown that for $\G=\wt{\C}$ and $M$ a maximal ideal (in particular a closed submodule) of $\wt{\C}$, $M^\perp=\{0\}$ and thus $M^{\perp\perp}=\{0\}$. Hence, the condition that $M$ is edged can not be dropped in the statements $(ii)$ and $(iii)$ of the previous corollary.
\end{remark}

The proof of Corollary \ref{corol_2} makes use of the following lemma.

\begin{lemma}
\label{lemma_scratch}
Let $a,b,c\in\wt{\R}$ with $a,b,c\ge 0$. If $b\le a$ and $ab\le ac$ then $b\le c$.
\end{lemma}
\begin{proof}
Fix a representative $(a_\eps)_\eps$ of $a$. Let for each $n\in\N$, $S_n=\{\eps\in (0,1): |a_\eps|\ge \eps^n\}$. Since $a$ is invertible with respect to $S_n$, we have that $be_{S_n}\le ce_{S_n}\le c$. Further, $0\le b-be_{S_n}=b e_{\co S_n}\le a e_{\co S_n}\to 0$, so $b=\lim_n be_{S_n}\le c$.
\end{proof}
Note that Lemma \ref{lemma_scratch} allows us to deduce for positive real generalized numbers $a$ and $c$ that $a^2\le ac$ implies $a\le c$ without involving any invertibility assumption on $a$.

\begin{corollary}
\label{corol_2}
Let $C$ be a closed edged subset of the Hilbert $\wt{\C}$-module $\G$ such that $\lambda C+(1-\lambda)C\subseteq C$ for all real generalized numbers $\lambda\in\wt{[0,1]}$. Then,
\[
\Vert{P_C(u_1)-P_C(u_2)}\Vert\le \Vert{u_1-u_2}\Vert
\]
for all $u_1,u_2\in\G$.
\end{corollary}
\begin{proof}
From Proposition \ref{prop_charac} we have that the inequalities $\Re\inner{u_1-P_C(u_1)}{P_C(u_2)-P_C(u_1)}\le 0$ and $\Re\inner{u_2-P_C(u_2)}{P_C(u_1)-P_C(u_2)}\le 0$ hold for all $u_1,u_2\in \G$. Thus,
\begin{multline*}
-\Re\inner{u_1-P_C(u_1)}{P_C(u_1)-P_C(u_2)}+\Re\inner{u_2-P_C(u_2)}{P_C(u_1)-P_C(u_2)}\\
=\Re\inner{u_2-u_1+P_C(u_1)-P_C(u_2)}{P_C(u_1)-P_C(u_2)}\\
=\Re\inner{u_2-u_1}{P_C(u_1)-P_C(u_2)}+\Vert{P_C(u_1)-P_C(u_2)}\Vert^2\le 0.
\end{multline*}
By the Cauchy-Schwarz inequality, it follows that
\[
\Vert{P_C(u_1)-P_C(u_2)}\Vert^2\le\Re\inner{u_1-u_2}{P_C(u_1)-P_C(u_2)}\le \Vert{u_1-u_2}\Vert\Vert{P_C(u_1)-P_C(u_2)}\Vert.
\]
Lemma \ref{lemma_scratch} allows us to deduce that $\Vert{P_C(u_1)-P_C(u_2)}\Vert\le\Vert{u_1-u_2}\Vert$. 
\end{proof}
When we work on the Hilbert $\wt{\C}$-module $\G_H$ and the set $C\subseteq\G_H$ is internal, the projection operator $P_C$ and the set $C^\perp$ have the following expected properties.
\begin{proposition}
\label{prop_projec_int}
\leavevmode
\begin{itemize}
\item[(i)] Let $H$ be a Hilbert space, $(C_\eps)_\eps$ a net of nonempty convex subsets of $H$ and $C:=[(C_\eps)_\eps]$. If $C\neq \emptyset$ then it is closed and edged and $P_C(u)=[(P_{{\overline{C}_\eps}}(u_\eps))_\eps]$ for all $u\in\G_H$.
\item[(ii)] In particular if $(C_\eps)_\eps$ is a net of closed subspaces of $H$ then $C^{\perp}=[({C_\eps}^\perp)_\eps]$.
\end{itemize}
\end{proposition}
\begin{proof}
$(i)$ Let $u=[(u_\eps)_\eps]\in\G_H$. Working at the level of representatives we have that $\Vert u_\eps-P_{{\overline{C}}_\eps}(u_\eps)\Vert=\inf_{w\in{\overline{C}}_\eps}\Vert u_\eps-w\Vert$. Let $v$ be an arbitrary element of $C$. Then there exists a representative $(v_\eps)_\eps$ such that $v_\eps\in C_\eps$ for all $\eps$ and $\Vert u_\eps-P_{{\overline{C}}_\eps}(u_\eps)\Vert\le  \Vert u_\eps-v_\eps\Vert$. Since $\Vert P_{{\overline{C}}_\eps}(u_\eps)\Vert\le \Vert P_{{\overline{C}}_\eps}(u_\eps)-u_\eps\Vert +\Vert u_\eps\Vert\le \Vert u_\eps-v_\eps\Vert +\Vert u_\eps\Vert$, the net $(\Vert P_{{\overline{C}}_\eps}(u_\eps)\Vert)_\eps$ is moderate . It follows that $\Vert u-[(P_{{\overline{C}}_\eps}(u_\eps))_\eps]\Vert \le \Vert u-v\Vert$ for all $v\in C$. Since, as proved in \cite{OV:07}, the set $C$ is closed and edged, by Corollary \ref{theo_projec} we have that $[(P_{{\overline{C}}_\eps}(u_\eps))_\eps]$ coincides with $P_C(u)$.

$(ii)$ The inclusion $[({C_\eps}^\perp)_\eps]\subseteq C^\perp$ is clear. If $u\in C^\perp$ then $P_C(u)=0$ and from the first assertion of this proposition the net $(\Vert P_{{{C}}_\eps}(u_\eps)\Vert)_\eps$ is negligible. So, $(u_\eps-P_{{{C}}_\eps}(u_\eps))_\eps$ is another representative of $u$ and $u_\eps-P_{{{C}}_\eps}(u_\eps)$ belongs to $C_\eps^\perp$ for each $\eps$.
\end{proof} 

We conclude this section with a version of the Hahn-Banach theorem for operators acting on Hilbert $\wt{\C}$-modules.
\begin{theorem}
\label{theo_HB}
Let $\G$ be a Hilbert $\wt{\C}$-module, $M$ a closed and edged $\wt{\C}$-submodule of $\G$ and $\mH$ a topological $\wt{\C}$-module. Let $f:M\to\mH$ be a continuous $\wt{\C}$-linear map. Then, $f$ can be extended to a continuous $\wt{\C}$-linear map on $\G$.
\end{theorem}
\begin{proof}
Take the projection operator $P_M$. From Corollary \ref{corol_1} we know that $P_M:\G\to M$ is $\wt{\C}$-linear and continuous and that $P_M(u)=u$ when $u\in M$. Thus, $f\circ P_M:\G\to\mH$ is a continuous $\wt{\C}$-linear extension of $f$.
\end{proof}
Since there exists a (without loss of generality, closed) submodule $M$ of $\wt{\C}$ and a continuous $\wt{\C}$-linear functional $T:M\to\wt{\C}$ which can not be extended to the whole of $\wt{\C}$ \cite{Vernaeve:07a}, we see that the condition that $M$ is edged can not be dropped in the previous theorem.

\section{Edged submodules}
\label{sec_fin}
In this section, we have a closer look at edged submodules of a Hilbert $\wt\C$-module (cf.\ Definition \ref{def_edged}). In the case of finitely generated submodules, it appears that edged submodules can be characterized by a topological condition (Theorem \ref{thm_Gram-Schmidt}). Some of the results hold for more general $\wt\R$-normed $\wt\K$-modules (here $\wt\K$ denotes either $\wt\R$ or $\wt\C$) fulfilling the following \emph{normalization property}.
\begin{definition}
\label{def_normalization}
A $\wt{\R}$-normed $\wt{\K}$-module $\G$ fulfills the normalization property if for all $u\in\G$ there exists $v$ in $\G$ such that $v\Vert u\Vert = u$.
\end{definition}
\begin{proposition} 
Let $\G$ be an $\wt\R$-normed $\wt\K$-module. Then $\G$ has the normalization property iff for each $u\in\G$ and $\lambda\in\wt\K$, the following holds: if $\norm{u}\le C\abs{\lambda}$, for some $C\in\wt\R$, then there exists $v\in\G$ such that $u=\lambda v$.
\end{proposition}
\begin{proof}
$\Rightarrow$: By absolute convexity of ideals in $\wt\K$, there exists $\mu\in\wt\K$ such that $\norm{u}=\mu\lambda$. By the normalization property, there exists $v\in\G$ such that $\norm{u}v=u$. Hence $u=\lambda (\mu v)$.\\
$\Leftarrow$: choose $\lambda=\norm{u}$.
\end{proof}
We observe that the Colombeau space $\G_E$ of generalized functions based on the normed space $E$ fulfills the normalization property.

\begin{proposition}
\label{prop_normal}
Let $E$ be a normed space. The $\wt{\K}$-module $\G_E$ fulfills the normalization property.
\end{proposition}
\begin{proof}
Let $u\in\G_E$ with representative $(u_\eps)_\eps$. We define $v_\eps$ as $u_\eps/\Vert u_\eps\Vert$ when $\Vert u_\eps\Vert\neq 0$ and $0$ otherwise. The net $(v_\eps)_\eps$ is clearly moderate and $v_\eps \Vert u_\eps\Vert=u_\eps$ for all $\eps$. This defines an element $v\in\G_E$ such that $v\Vert u\Vert= u$.
\end{proof}
Note that Definition \ref{def_edged} can clearly be stated in the more general context of $\wt{\R}$-normed $\wt{\K}$-modules. We recall that a Banach $\wt\K$-module is a complete ultra-pseudo-normed $\wt\K$-module \cite{Garetto:05a, Garetto:06b}.
\begin{proposition}
\label{prop_inhe}
Let $M$ be a closed submodule of a Banach $\wt\K$-module $\G$ and let $\G/M$ be endowed with the usual quotient topology.
\begin{itemize}
\item[(i)] $\G/M$ is a Banach $\wt\K$-module.
\item[(ii)] If $\G$ is $\wt\R$-normed and $M$ is edged, then $\G/M$ is $\wt\R$-normed.
\item[(iii)] If $\G$ is a Hilbert $\wt\C$-module and $M$ is edged, then $\G/M$ is a Hilbert $\wt\C$-module.
\item[(iv)] If $\G$ is a Hilbert $\wt\C$-module satisfying the normalization property and $M$ is edged, then $M$ has the normalization property.
\end{itemize}
\end{proposition}
\begin{proof}
(i) By \cite[Example 1.12]{Garetto:05a}, the relative topology on $\G/M$ is generated by one ultra-pseudo-seminorm. It is easy to check that this ultra-pseudo-seminorm is an ultra-pseudo-norm if $M$ is closed. By Proposition 4.25 in \cite{Garetto:06b} we have that $\G/M$ is complete.\\
(ii) for $u\in\G$, let $\bar u:=u+M\in\G/M$. We define $\norm{.}$: $\G/M\to\wt\R$: $\norm{\bar u}=\inf_{w\in M}\norm{u-w}$. As $M$ is edged, the infimum exists. It is easy to see that $\norm{\bar u}$ does not depend on the representative $u\in\G$, that $\norm{\bar u}\ge 0$ and $\norm{\bar 0}=0$. If $\norm{\bar u}=\inf_{w\in M}\norm{u-w}=0$, then there exists a sequence $(w_n)_n$ with $w_n\in M$ and $u=\lim_n w_n$. Hence $u\in \overline M=M$ and $\bar u=0$. Let $u_1,u_2\in \G$ and $w_1,w_2\in M$. Then
\[
\norm{\bar u_1 + \bar u_2}=\inf_{w\in M}\norm{u_1+u_2-w}\le\norm{u_1 + u_2 - (w_1 + w_2)}\le\norm{u_1 - w_1} + \norm{u_2 - w_2}.\]
Taking the infimum over $w_1\in M$ and $w_2\in M$, we obtain $\norm{\bar u_1 + \bar u_2}\le\norm{\bar u_1} + \norm{\bar u_2}$. Now let $u\in\G$ and $\lambda\in\wt\K$. Then
\[\norm{\lambda \bar u}=\inf_{w\in M}\norm{\lambda u-w}\le \inf_{w\in M}\norm{\lambda u-\lambda w} =\inf_{w\in M}\abs\lambda\norm{u-w} =\abs\lambda\inf_{w\in M}\norm{u-w}=\abs\lambda\norm{\bar u}.\]
If $\lambda=0$, the converse inequality trivially holds. If $\lambda\ne 0$, let $S\subseteq(0,1]$ with $e_S\ne 0$ such that $\lambda$ is invertible w.r.t.\ $S$, say $\lambda\mu=e_S$, and let $w\in M$. Then
\[
\norm{\lambda e_S u - w}\ge \norm{\lambda e_S u - w} e_S
=\norm{\lambda e_S u - \lambda e_S (\mu w)}
\ge\inf_{w\in M}\norm{\lambda e_S u - \lambda e_S w}
=\abs{\lambda}e_S\norm{\bar u}.
\]
Fix a representative $(\lambda_\eps)_\eps$ of $\lambda$ and let $S_n=\{\eps\in(0,1]: \abs{\lambda_\eps}\ge\eps^n\}$ for each $n\in\N$. Then $e_{S_n}\ne 0$ and $\lambda$ is invertible w.r.t.\ $S_n$ for sufficiently large $n$. As $\lambda =\lim_n \lambda e_{S_n}$, by the continuity of the $\wt\R$-norm,
\[\norm{\lambda u - w}=\lim_n\norm{\lambda e_{S_n}u - w}\ge \lim_n\abs\lambda e_{S_n}\norm{\overline u}=\abs\lambda \norm{\overline u}.\]
Taking the infimum over $w\in M$, we obtain $\norm{\lambda \overline u}\ge\abs\lambda \norm{\overline u}$. So $\norm{.}$ is an $\wt\R$-norm on $\G/M$. By the continuity of the sharp norm $\abs{.}_\esp$ on $\wt\R$ and the fact that $\abs{.}_\esp$ is increasing on $\{x\in\wt\R: x\ge 0\}$, the corresponding ultra-pseudo-norm
\[
{\mathcal P}(\overline u)=\abs{\norm{\overline u}}_\esp
=\abs{\inf_{w\in M}\norm{u-w}}_\esp
=\inf_{w\in M}\abs{\norm{u-w}}_\esp
=\inf_{w\in M}{\mathcal P}(u-w)
\]
is the usual quotient ultra-pseudo-norm.\\
(iii) The map $f$: $\G/M\to M^\perp$: $u+M\mapsto P_{M^\perp}(u)$ is well-defined, since for $v\in\G$ with $u+M=v+M$, $P_{M^\perp}(u)-P_{M^\perp}(v)=P_{M^\perp} (u-v)=0$. Further, $f$ is $\wt\C$-linear and surjective and $\norm{u+M}=\inf_{w\in M}\norm{u-w}=\norm{u-P_M(u)}=\norm{P_{M^\perp}(u)}$, so $f$ is an algebraic and isometric isomorphism. Hence $\G/M$ is a Hilbert $\wt\C$-module for the scalar product $\inner{u+M}{v+M}_{\G/M}:=\inner{P_{M^\perp}(u)}{P_{M^\perp}(v)}_{\G}$.\\
(iv) Let $u\in M$. If there exists $v\in\G$ such that $\norm{u}v=u$, then $P_M(v)\in M$ and by the linearity of the projection operator, $\norm{u}P_M(v)=P_M(\norm{u}v)=P_M(u)=u$.
\end{proof}
\subsection{Cyclic submodules}
\begin{definition}
Let $\G$, $\mathcal{H}$ be $\wt\K$-modules with $\wt\R$-norm $\norm{.}$. A map $\phi$: $\G\to{\mathcal H}$ is an isometry iff 
$\norm{\phi(u)-\phi(v)}=\norm{u-v}$, for each $u,v\in\G$.
\end{definition}

The submodules considered in the sequel are always $\wt{\K}$-submodules. We recall that a $\wt\K$-module $M$ is called cyclic iff it is generated by one element, i.e., there exists $u\in M$ such that $M=u\wt\K$. An ideal $I$ of $\wt{\K}$ (in short $I\ideal\wt\K$) which is generated by one element is said to be principal. Before proving Proposition \ref{prop_cyclic_module} we collect some results concerning the ideals of $\wt{\K}$ which will be used later. Detailed proofs can be found in \cite{Vernaeve:07a}.

\begin{proposition}
\label{preli_prop_ideal}
Let $I\ideal\wt\K$. 
\begin{itemize}
\item[(i)] $I$ is absolutely order convex, i.e., if $x\in I$, $y\in\wt{\K}$ and $|y|\le|x|$ then $y\in I$.
\item[(ii)] If $x\in \overline{I}$ is invertible w.r.t. $S\subseteq(0,1]$ then $e_S\in I$.
\item[(iii)] A principal ideal $I$ of $\wt{\K}$ is closed if and only if there exists $S\subseteq(0,1]$ such that $I=e_S\wt{\K}$.
\end{itemize}
\end{proposition}
\begin{theorem}\label{thm_edged_ideals}
For an ideal $I$ of $\wt\K$, the next statements are equivalent:
\begin{itemize}
\item[(i)] $I$ is internal
\item[(ii)] $I$ is closed and edged
\item[(ii')] $I$ is edged
\item[(iii)] $I$ is a direct summand of $\wt\K$, i.e., there exists an ideal $J$ of $\wt\K$ such that $I+J=\wt\K$ and $I\cap J=\{0\}$
\item[(iv)] $(\exists S\subseteq (0,1])(I=e_S\wt\K)$.
\end{itemize}
\end{theorem}
\begin{proof}
$(i)\Rightarrow(ii)$: holds for any nonempty internal set of $\wt\K$ \cite{OV:07}.\\
$(ii)\Rightarrow(iii)$: by Corollary \ref{corol_1}, $I+I^\perp=\wt\K$ and $I\cap I^\perp=\{0\}$.\\
$(iii)\Rightarrow(iv)$: by hypothesis, $1=a+b$ with $a\in I$ and $b\in J$. As $ab\in I\cap J$, $ab=0$. Let $x\in I$. Then $xb\in I\cap J$, so $xb=0$ and $x=x(a+b)=xa$. So $I=a\wt\K$. As $a=a(a+b)=a^2$, $a$ is idempotent, hence $a=e_S$ for some $S\subseteq(0,1]$.\\
$(iv)\Rightarrow(i)$: let $I_\eps=\K$, if $\eps\in S$ and $I_\eps=\{0\}$, otherwise. Then $I=[(I_\eps)_\eps]$.\\
$(ii)\Leftrightarrow(ii')$: let $I$ be edged. As $\overline I$ is closed and edged, the previous equivalences show that $\overline I=e_S\wt\K$, for some $S\subseteq(0,1]$. But if $e_S\in \overline I$, then $e_S\in I$ by Proposition \ref{preli_prop_ideal}, so $I=\overline I$.
\end{proof}

\begin{proposition}\label{prop_cyclic_module}
Let $M=u\wt\K$ be a cyclic submodule of an $\wt\R$-normed $\wt\K$-module $\G$.
\begin{enumerate}
\item $M$ is isometrically isomorphic with the ideal $\norm{u}\wt\K\ideal\wt\K$.
\item $\overline M$ is isometrically isomorphic with an ideal $I\ideal\wt\K$ and $\norm{u}\wt\K\subseteq I\subseteq \overline{\norm{u}\wt\K}$. If $\G$ is a Banach $\wt\K$-module, then $I=\overline{\norm{u}\wt\K}$.
\item If $v\in\overline M$ and $\norm{v}$ is invertible w.r.t.\ $S$, then $ve_S\in M$.
\item If $v\in\overline M$ and $\norm{v}\le c\norm{u}$, for some $c\in\wt\R$, then $v\in M$.
\item If $v\in M$ and $\norm{u}\le c\norm{v}$, for some $c\in\wt\R$, then $M=v\wt\K$.
\item If there exists $w\in\G$ and $S\subseteq(0,1]$ such that $M=w\wt\K$ and $\norm{w}=e_S$ (or equivalently, $\norm{u}$ is invertible w.r.t.\ $S$ and zero w.r.t.\ $\co S$), then $M$ is closed.
\item If $\G$ is a Banach $\wt\K$-module, then $M$ is closed iff there exists $w\in\G$ and $S\subseteq(0,1]$ such that $M=w\wt\K$ and $\norm{w}=e_S$.
\item If $M$ is closed, then any edged submodule $N$ of $M$ is closed and cyclic.
\item If $\G$ has the normalization property, then $M$ is contained in a closed cyclic submodule of $\G$.
\item If $\G$ has the normalization property, and $M$ is edged, then $M$ is closed.
\item If $\G$ is a Hilbert $\wt\K$-module, $v\in \G$ and $\norm{v}\le c\norm{u}$, for some $c\in\wt\R$, then there exists $P_M(v)\in M$, which is both the unique element of $M$ such that $\norm{v-P_M(v)}= d(v, M)$, and the unique element of $M$ such that $\inner{v}{u}=\inner{P_M(v)}{u}$.
\item If $\G$ is a Hilbert $\wt\K$-module and $M$ is closed, then $M$ is edged. If $u$ is a generator of $M$ with idempotent norm, then for any $v\in\G$, $P_M(v)=\inner{v}{u}u$.
\item If $\G$ is a Hilbert $\wt\K$-module with the normalization property, then $M^{\perp\perp}=\overline M$.
\end{enumerate}
\end{proposition}
\begin{proof}
(1) Define $\phi$: $M\to \wt\K$: $\phi(\lambda u)=\lambda\norm{u}$ ($\lambda\in\wt\K$). Then the equality $\norm{\lambda u - \mu u}=\abs{\phi(\lambda u) - \phi(\mu u)}$ shows that $\phi$ is well-defined and isometric (hence also injective). It is easy to check that $\phi$ is $\wt\K$-linear and $\phi(M)=\norm{u}\wt\K$.\\
(2) We extend $\phi$: $M\to\wt\K$ to a map $\overline M\to\wt\K$ by defining $\phi(\lim_n \lambda_n u):=\lim_n \phi(\lambda_n u)$. Because $(\lambda_n u)_n$ is a Cauchy-sequence, $(\phi(\lambda_n u))_n$ is also a Cauchy-sequence in $\wt\K$, and hence convergent in $\wt\K$. To see that $\phi$ is well-defined, let $\lim_n \lambda_n u=\lim_n \mu_n u$. Then also the interlaced sequence $(\lambda_1 u, \mu_1 u,\dots, \lambda_n u,\mu_n u, \dots)$ is a Cauchy-sequence. Hence also $(\phi(\lambda_1 u), \phi(\mu_1 u),\dots, \phi(\lambda_n u),\phi(\mu_n u), \dots)$ is convergent to $\lim_n\phi(\lambda_n u)=\lim_n\phi(\mu_n u)$. It is easy to check that also the extended $\phi$ is linear and isometric and that $\phi(\overline{M})$ is an ideal of $\wt{\K}$ such that $\Vert u\Vert\wt{\K}\subseteq\phi(\overline{M})\subseteq\overline{\Vert u\Vert\wt{\K}}$. If $\G$ is complete, we find that the image under $\phi^{-1}$ of any convergent sequence in $\norm{u}\wt\K$ (say to $\lambda\in\overline{\norm{u}\wt\K}$) is a Cauchy sequence, and hence convergent to an element $v\in \overline M$. By definition of the extended map $\phi$, we have that $\phi(v)=\lambda$ and therefore $\phi(\overline{M})=\overline{\Vert u\Vert\wt{\K}}$.\\
(3) As $\phi$ is an isometry, $\abs{\phi(v)}$ is invertible w.r.t.\ $S$. As $\phi(v)\in\overline{\norm{u}\wt\K}$, by Proposition \ref{preli_prop_ideal}$(ii)$, $e_S\in\norm{u}\wt\K$. Hence also $\phi(ve_S)=\phi(v)e_S\in\phi(M)$. So $ve_S\in M$ by the injectivity of $\phi$.\\
(4) As $\phi$ is an isometry, $\abs{\phi(v)}=\norm{v}\le c\norm{u}$. By absolute order convexity of ideals in $\wt\K$, $\phi(v)\in\norm{u}\wt\K=\phi(M)$. So $v\in M$ by the injectivity of $\phi$.\\
(5) As $\norm{u}\le c\norm{v}= c\abs{\phi(v)}$, by absolute order convexity of ideals in $\wt\K$, $\norm{u}= \mu\phi(v)$, for some $\mu\in\wt\K$. So $\phi(u)=\norm{u}=\phi(\mu v)$, and $u\in v\wt\K$ by the injectivity of $\phi$. It follows that $M=v\wt{\K}$.\\
(6) Let us assume that $\norm{u}$ is invertible w.r.t.\ $S$ and zero w.r.t.\ $\co S$. Then, there exists $\lambda\in\wt{\R}$ such that $\lambda\Vert u\Vert=e_S=\Vert \lambda u\Vert$ and $\Vert u\Vert e_{\co S}=0$. It follows that $u e_{\co S}=0$ or equivalently $u=u e_S$. Hence, $u\wt{\K}=\lambda u\wt{\K}$ and we can choose $w=\lambda u$. Now, $\phi(M)=\norm{w}\wt\K=e_S\wt\K$ is closed in $\wt\K$, hence complete, so also $M$ is complete, hence closed.\\
(7) Let $M$ be closed. As a closed submodule of a Banach $\wt\K$-module, $M$ is complete. By part~1, also $\norm{u}\wt\K$ is a complete, hence closed principal ideal of $\wt\K$. By Proposition \ref{preli_prop_ideal}$(iii)$ this implies that $\norm{u}\wt\K=e_S\wt\K$ for some $S\subseteq(0,1]$.\\
(8) By part~7, we may assume that $\norm{u}=e_S$, for some $S\subseteq(0,1]$. If $N$ is edged, then it is reachable from any element of $M$. By the isometry, also $\phi(N)$ is reachable from every element of $\phi(M)=e_S\wt\K$. Hence, for each $\lambda\in \phi(N)\subseteq e_S\wt\K$ and $\mu\in \wt\K$,
\[
\abs{\mu-\lambda}
=\abs{\mu-\lambda}e_{\co S}+\abs{\mu-\lambda}e_S
=\abs{\mu}e_{\co S} + \abs{\mu e_S-\lambda},
\]
so $\phi(N)$ is also reachable from $\mu$. This implies that $\phi(N)$ is an edged ideal of $\wt{\K}$ and by Theorem \ref{thm_edged_ideals} that $\phi(N)$ is closed (hence complete) and principal. So $N$ is closed and cyclic.\\
(9) Consider $v\in\G$ with $\norm{u}v=u$. Then $\norm{u}\norm{v}=\norm{u}$, so $\norm{u}(1-\norm{v})=0$. By a characterization of zero divisors in $\wt\R$, there exists $S\subseteq (0,1]$ such that $\norm{u}e_{\co S}=0$ and $\norm{v}e_S= e_S$. Let $w=v e_S$. Then $\norm{w}=e_S$, hence the cyclic submodule $w\wt\K$ is closed by part~6. Further, $u=u e_S=\norm{u}w$, so $M=u\wt\K\subseteq w\wt\K$.\\
(10) By parts~8 and 9.\\
(11) By the Cauchy-Schwarz inequality (\ref{cs_ineq}), $\abs{\inner{v}{u}}\le \norm{u}\norm{v}\le c\norm{u}^2$, so by absolute order convexity of ideals in $\wt\K$, there exists $\lambda\in\wt\K$ such that $\inner{v}{u}=\lambda \norm{u}^2$. We show that $P_M(v)=\lambda u$. First, for $\mu\in\wt\K$, $\inner{\mu u}{u}=\inner{v}{u}$ iff $(\mu-\lambda)\norm{u}^2=0$. 
It follows that also $\abs{\mu-\lambda}^2\norm{u}^2=\norm{(\mu-\lambda)u}^2=0$, so $\mu u=\lambda u$, and $P_M(v)$ is the unique element in $M$ such that $\inner{v}{u}=\inner{P_M(v)}{u}$. From this equality, it follows that $\norm{v-P_M(v)+\mu u}^2=\norm{v-P_M(v)}^2+\norm{\mu u}^2$, which is only minimal if $\mu u=0$.\\
(12) By part~7, we can suppose that $\norm{u}=e_S$ for some $S\subseteq(0,1]$. In particular, $\norm{e_{\co S} u}=0$, so $u=e_S u$. Let $v\in\G$. Let $p=\inner{v}{u}u\in M$. Then $\inner{v-p}{u}=\inner{v}{u}-\inner{v}{u}\inner{u}{u} = \inner{v}{e_S u}-\inner{v}{u}e_S=0$. It follows from Corollary \ref{corol_1_new} that $M$ is reachable from $v$. So, $M$ is edged and $p=P_M(v)$.\\
(13) By parts 7, 9 and 12, $M\subseteq w\wt\K$ with $w\wt\K$ closed and edged, and $\norm w=e_S$ for some $S\subseteq(0,1]$.
Let $v\in M^{\perp\perp}$. If $\lambda\in\wt\K$ and $\inner{u}{\lambda w}=0$, then $\inner{v}{\lambda w}=0$. By Corollary \ref{corol_1}, $M^{\perp\perp}\subseteq (w\wt\K)^{\perp\perp}=w\wt\K$. Let $\phi$ be the isometric embedding $w\wt\K\to \wt\K$: $\phi(\lambda w)=\lambda e_S$. Since $\phi$ is a $\wt\K$-linear isometry, $\phi$ also preserves the scalar product. For $a\in e_S\wt\K$, $\inner{a}{\phi(\lambda w)}=a\bar{\lambda}e_S=a\bar{\lambda}$. So if $\lambda\in\wt\K$ and $\phi(u)\bar\lambda=0$, then $\phi(v)\bar\lambda=0$, i.e., $\phi(v)$ is orthogonal to any $\lambda\in\phi(M)^\perp$ (orthogonal complement in $\wt\K$). Hence $\phi(v)\in\phi(M)^{\perp\perp}=\overline{\phi(M)}$ since $\phi(M)$ is a principal ideal of $\wt\K$ \cite{Vernaeve:07a}. By part 2, $\overline{\phi(M)}=\phi(\overline M)$. By the injectivity of $\phi$, $v\in\overline M$. The converse inclusion holds for any submodule.
\end{proof}

The following example shows that the normalization-property can not be dropped in Proposition \ref{prop_cyclic_module} part~10.
\begin{example}\label{ex_edged_not_closed}
Let for each $m\in\N$, $S_m\subseteq(0,1]$ with $0\in\overline{S}_m$ and such that $S_n\cap S_m=\emptyset$ if $n\ne m$.\\
Let $\beta_\eps=\eps^m$, for each $\eps\in S_m$, and $\beta_\eps=0$ for $\eps\in(0,1]\setminus\bigcup S_n$. Let $\beta\in\wt\R$ be the element with representative $(\beta_\eps)_\eps$.\\
Then $\G=\overline{\beta\wt\K}$, the closure in $\wt\K$ of $\beta\wt\K$, is a Hilbert $\wt\K$-module and ${\beta\wt\K}\neq\overline{\beta\wt\K}$ as it is proven in \cite{Vernaeve:07a}. Then $M=\beta\wt\K$ is an edged cyclic submodule of $\G$, since for each $u\in\G$, $\overline{\inf}_{v\in M}|u-v|=0$. Yet $\beta\wt\K$ is not closed in $\G$.
\end{example}
The following example shows that Proposition \ref{prop_cyclic_module} part~1 does not hold for a general Banach $\wt\K$-module $\G$. In particular it provides an example of a Banach $\wt{\K}$-module which is not $\wt{\R}$-normed and proves that a quotient of a Hilbert $\wt{\K}$-module over a closed but not edged submodule is not necessarily a Hilbert $\wt{\K}$-module itself. We recall that for $\gamma\in\wt{\K}$, $\Ann(\gamma)$ denotes the set of all $x\in\wt{\K}$ such that $x\gamma=0$.
\begin{example}
Let $\beta\in\wt\R$ as in Example \ref{ex_edged_not_closed}. Then $\G=\wt\K/\overline{\beta\wt\K}$ is a cyclic Banach $\wt\K$-module. Yet $\G$ is not algebraically isomorphic with an ideal of $\wt\K$.
\end{example}
\begin{proof}
By Proposition \ref{prop_inhe}, $\G$ is a Banach $\wt\K$-module. For $x\in\wt\K$, we denote by $\bar x$ the class $x+\overline{\beta\wt\K} \in \G$. Then $\G$ is generated by the element $\bar 1\in\G$. Suppose that $\G\cong I$ (as a $\wt\K$-module), for some $I\ideal\wt\K$. Then there exists $a\in\wt\K$ such that $I=a\wt\K$. By the algebraic isomorphism, $x\bar 1=0$ iff $xa=0$, $\forall x\in\wt\K$. So the annihilator ideal $\Ann(a)=\Ann(\bar 1)=\overline{\beta\wt\K}$. But $\Ann(a)$ is either principal, either it is not the closure of a countably generated ideal, whereas $\overline{\beta\wt\K}$ is the closure of a countably generated ideal, but not principal \cite{Vernaeve:07a}.
\end{proof}
By means of Proposition \ref{prop_cyclic_module}, we are now able to prove that the formulation of convexity on $C$ given in Proposition \ref{prop_charac} automatically holds for all the values of $\lambda$ in $\wt{[0,1]}=\{x\in\wt\R: 0\le x\le 1\}$.
\begin{proposition}
\label{prop_new_convex}
Let $C$ be a closed edged subset of the Hilbert $\wt\C$-module $\G$ such that $\lambda C + (1-\lambda)C\subseteq C$ for all $\lambda\in\{\caninf^q\}_{q\in\N}\cup\{\frac{1}{2}\}$. Then $\lambda C + (1-\lambda)C\subseteq C$ for all $\lambda\in\wt{[0,1]}$.
\end{proposition}
\begin{proof}
Let $u,u'\in C$ and $\lambda\in\wt{[0,1]}$. We show that $v=\lambda u + (1-\lambda)u'\in C$. As the properties of $C$ are translation invariant, we may suppose that $u'=0$ (so $v=\lambda u$). If $\norm{u}=0$, then trivially $v=0\in C$. So, without loss of generality $\Vert u\Vert\neq 0$. Let $S\subseteq (0,1]$ with $e_S\ne 0$ such that $\norm{u}$ is invertible w.r.t.\ $S$. Then $\norm{ue_S}=\norm{u}e_S$ is invertible w.r.t.\ $S$ and zero w.r.t.\ $\co S$, so $M=ue_S\wt\C$ is a closed, edged submodule by Proposition \ref{prop_cyclic_module}. Let $P_M(P_C(v)e_S) = (\mu + i\kappa) u e_S$ for some $\mu,\kappa\in\wt\R$. Then $P_C(v)e_S = (\mu + i\kappa) u e_S + w$, with $\inner{u}{w}=\inner{u e_S}{w}=0$. Fix representatives $(\lambda_\eps)_\eps$ of $\lambda$ and $(\mu_\eps)_\eps$ of $\mu$. Let $T=\{\eps\in S: \lambda_\eps\le\mu_\eps\}$. Then $0\le \lambda e_T\le\mu e_T$. By Proposition \ref{prop_charac},
\begin{multline*}
0\ge \Re \inner{v-P_C(v)}{-P_C(v)}e_T 
=-\Re\inner{v}{P_C(v)e_T} + \norm{P_C(v)}^2e_T\\
=-\lambda\mu\norm{u}^2e_T + (\mu^2 + \kappa^2)\norm{u}^2e_T + \norm{w}^2e_T
\end{multline*}
so
\[
0\le \norm{w}^2 e_T + \kappa^2\norm{u}^2 e_T\le (\lambda-\mu)\mu\norm{u}^2e_T\le 0.
\]
By the invertibility of $\norm{u}$ w.r.t.\ $S$, $(\lambda-\mu)\mu e_T=w e_T=\kappa e_T =0$. Then also $0\le(\lambda - \mu)^2e_T = (\lambda-\mu)\lambda e_T\le 0$, and $\lambda e_T = \mu e_T$.\\
Denoting $U=S\setminus T$, we have $\mu e_U\le \lambda e_U\le e_U$. Again by Proposition \ref{prop_charac},
\begin{multline*}
0\ge \Re \inner{v-P_C(v)}{u-P_C(v)} e_U
=\Re \inner{v}{u}e_U - \Re\inner{v}{P_C(v)e_U}- \Re\inner{P_C(v)}{u}e_U + \norm{P_C(v)}^2 e_U\\
=\lambda \norm{u}^2e_U - \lambda\mu\norm{u}^2e_U - \mu\norm{u}^2e_U + (\mu^2 + \kappa^2)\norm{u}^2e_U + \norm{w}^2e_U
\end{multline*}
so
\[
0\le \norm{w}^2 e_U + \kappa^2\norm{u}^2 e_U\le (\lambda - \mu)(\mu - 1)\norm{u}^2e_U\le 0,
\]
hence, as before, $(\lambda - \mu)(\mu - 1) e_U= w e_U=\kappa e_U=0$. Then also $0\le(\lambda-\mu)^2e_U \le (\lambda-\mu)(1-\mu)e_U= 0$, and $\lambda e_U=\mu e_U$.\\
Together, this yields $w=w e_S=0$, $\kappa e_S=0$ and $\lambda e_S = \mu e_S$. It follows that $P_C(v)e_S=\mu u e_S = \lambda u e_S = v e_S$. Now fix a representative $(\norm{u}_\eps)_\eps$ of $\norm{u}$ and consider $S_n=\{\eps\in (0,1]: \norm{u}_\eps\ge \eps^n\}$, for $n\in\N$. Since $\norm{u}\neq 0$, $e_{S_n}\ne 0$ for sufficiently large $n$. As $\norm{u}$ is invertible w.r.t.\ $S_n$, $P_C(v)e_{S_n}=v e_{S_n}$ for sufficiently large $n$. Further, as $0\in C$, $\norm{P_C(v)}\le\norm{P_C(v)-v}+\norm{v}\le 2\norm{v}\le 2\abs{\lambda}\norm{u}$. As $\lim_n\norm{u}e_{\co S_n}=0$, also $\lim_n\norm{P_C(v)}e_{\co S_n}=\lim_n\norm{v}e_{\co S_n}=0$, so $v=\lim_n v e_{S_n}=\lim_n P_C(v) e_{S_n}=P_C(v)\in C$.
\end{proof}

\begin{theorem}
\label{theo_hans_new}
\leavevmode
\begin{itemize}
\item[(i)] Let $\G$ be a Hilbert $\wt\K$-module with the normalization property. Then a cyclic submodule is edged iff it is closed iff it is generated by an element with idempotent $\wt{\R}$-norm.
\item[(ii)] Let $\G_E$ be a Banach $\wt\K$-module constructed by means of a Banach space $E$. Then a cyclic submodule is edged iff it is closed iff it is generated by an element with idempotent norm iff it is internal.
\end{itemize}
\end{theorem}
\begin{proof}
(i) Follows by proposition \ref{prop_cyclic_module}, assertions 7, 10, 12.\\
(ii) By Proposition \ref{prop_normal}, $\G_E$ has the normalization property. So, by Proposition \ref{prop_cyclic_module} we already have the implications edged $\implies$ closed $\implies$ generated by an element with idempotent $\wt{\R}$-norm.\\
Let $M=u\wt\K$ be a submodule of $\G_E$, and suppose that $\norm{u}=e_S$, for some $S\subseteq(0,1]$. We show that $M$ is internal. Let $(u_\eps)_\eps$ be a representative of $u$. As $u e_{\co S}=0$, we may suppose that $u_\eps=0$, for each $\eps\in\co S$. Let $A_\eps=u_\eps\K$, for each $\eps\in (0,1]$. If $v=\lambda u$, for some $\lambda\in\wt\K$, then there exist representatives such that $v_\eps=\lambda_\eps u_\eps$, so $v\in[(A_\eps)_\eps]$. For the converse inclusion, if $v\in[(A_\eps)_\eps]$, we find $\lambda_\eps\in\K$ such that, on representatives, $v_\eps=\lambda_\eps u_\eps$. We may assume that $\lambda_\eps=0$ for $\eps\in\co S$. Then, denoting by $\chi_S$ the characteristic function of $S$, the net $(\abs{\lambda_\eps})_\eps = (\frac{\norm{v_\eps}}{\norm{u_\eps}}\chi_S(\eps))_\eps$ is moderate (since $\norm{u}$ is invertible w.r.t.\ $S$). So $(\lambda_\eps)_\eps$ represents $\lambda\in\wt\K$ and $v=\lambda u$.\\
Finally, any internal set in $\G_E$ is edged \cite{OV:07}.
\end{proof}

In spite of the obtained results, some elementary operations on cyclic modules appear not to preserve the property of being edged. Even in $\wt\R^2$, neither intersections nor projections, nor sums preserve this property.
\begin{example}
Let $\beta\in\wt\R$ as in Example~\ref{ex_edged_not_closed}. Then $(1,\beta)\wt\R\cap (1,0)\wt\R=\Ann(\beta)\times\{0\}$ is not edged in $\wt\R^2$ (since $\Ann(\beta)$ is not edged in $\wt\R$ \cite{Vernaeve:07a}). Since $\norm{(1,\beta)}$ is invertible by Theorem \ref{theo_hans_new} we have that $(1,\beta)\wt\R$ is edged.
\end{example}

\begin{example}
Let $\beta\in\wt\R$ as in Example~\ref{ex_edged_not_closed}. Let $M=(1,0)\wt\R\subseteq\wt{\R}^2$. Then $P_M((\beta,1)\wt\R)=(\beta,0)\wt\R$ is not edged in $\wt\R^2$ (since $\beta\wt\R$ is not generated by an idempotent \cite{AJOS:06}).
\end{example}
This gives also an example of a projection of a closed submodule on a closed submodule which is not closed.

\begin{example}
Let $\beta\in\wt\R$ as in Example~\ref{ex_edged_not_closed}. Let $M=(1,\beta)\wt\R+(1,0)\wt\R\subset\wt\R^2$. As $\norm{(1,\beta)}$ and $\norm{(1,0)}$ are invertible, $M$ is the sum of cyclic edged submodules. Yet $M$ is not edged, since $M=(0,\beta)\wt\R+(1,0)\wt\R$, so $\overline{\inf}_{v\in M}\Vert(0,a)-v\Vert=\overline{\inf}_{\lambda,\mu\in\wt\R}(\abs{\mu}^2 + \abs{a-\lambda\beta}^2)^{1/2} = \overline{\inf}_{\lambda\in\wt\R}\abs{a-\lambda\beta}$ does not exist for some $a\in\wt{\R}$, since $\beta\wt\R$ is not edged.
\end{example}
This gives also an example of two submodules $M$, $N$ with $\overline M +\overline N\ne\overline{M+N}$.

Concerning \emph{direct} sums of edged submodules, see however Theorem \ref{thm_edged_direct_sum} below.

\subsection{Submodules generated by $m\ge 1$ elements}
\begin{theorem}\label{thm_Gram-Schmidt}
Let $\G$ be a Hilbert $\wt\K$-module and $M$ a submodule of $\G$ generated by $m$ elements. Then
\begin{enumerate}
\item $M$ is a direct sum of $m$ mutually orthogonal cyclic modules (`interleaved Gram-Schmidt').
\item $M$ is isometrically isomorphic with a submodule $M'$ of $\wt\K^m$.
\item $\overline M$ is isometrically isomorphic with $\overline{M'}$ (closure in $\wt\K^m$).
\item $M$ is closed iff $M$ is a direct sum of $m$ mutually orthogonal closed cyclic modules.
\item If $M$ is closed, then $M$ is edged.
\item If $M$ is closed, any edged submodule $N$ of $M$ is closed and finitely generated.
\item If $\G$ has the normalization property and $M$ is a direct sum of mutually orthogonal cyclic modules $M_1$, \dots, $M_m$, then there exist mutually orthogonal closed cyclic modules $N_j$ such that $M_j\subseteq N_j$, for $j=1$,\dots, $m$.
\item If $\G$ has the normalization property and $M$ is edged, then $M$ is closed.
\end{enumerate}
\end{theorem}
\begin{proof}
(1) We proceed by induction on $m$. The case $m=1$ is trivial.\\
Let $M=u_1\wt\K+\cdots + u_m\wt\K$. Fix representatives $(\norm{u_j}_\eps)_\eps$ of $\norm{u_j}$ and define recursively for $j=1$, \dots, $m$
\[S_j=\big\{\eps\in (0,1]: \norm{u_j}_\eps\ge\max_{k\ne j}\norm{u_k}_\eps\big\}\setminus \{S_1,\dots,S_{j-1}\}.\]
Then $e_{S_j}e_{S_k}=0$ if $j\ne k$, and $e_{S_1}+\cdots+e_{S_m}=1$.
By Proposition~\ref{prop_cyclic_module}, we can project $u_j e_{S_1}$ on $N=u_1 e_{S_1}\wt\K$, obtaining $\tilde u_j=u_j e_{S_1}-P_{N}(u_j e_{S_1})$ with $\inner{u_1}{\tilde u_j}=0$ and $\tilde u_j=\tilde u_j e_{S_1}$ ($j>1$). With $N'=\tilde u_2 \wt\K +\cdots +\tilde u_m \wt\K$, we also have $Me_{S_1}=u_1 e_{S_1}\wt{\K}+ N'$ and $u_1\in N'^{\perp}$. By induction, there exist mutually orthogonal generators $v_2^{(1)}$, \dots, $v_m^{(1)}$ of $N'$. Since $N'e_{S_1}=N'$, $v_j^{(1)}=v_j^{(1)}e_{S_1}$, for all $j$. With $v_1^{(1)}=u_1 e_{S_1}$, we obtain $m$ mutually orthogonal generators of $Me_{S_1}$. Similarly, we obtain $m$ mutually orthogonal generators $v_1^{(k)}$, \dots, $v_m^{(k)}$ of $Me_{S_k}$ ($k=1$, \dots, $m$) (in particular satisfying $v_j^{(k)}=v_j^{(k)}e_{S_k}$). Then $v_j=\sum_{k=1}^m v_j^{(k)}$ ($j=1$, \dots, $m$) are mutually orthogonal generators of $M$. By orthogonality, it follows that the sum is a direct sum: if $\sum_j\lambda_jv_j=0$, for some $\lambda_j\in\wt\K$, then $0= \inner[]{\sum_j\lambda_jv_j}{\sum_j\lambda_jv_j} = \sum_j\norm{\lambda_j v_j}^2$, so each $\lambda_j v_j=0$.\\
(2) By part~(1), $M=v_1\wt\K +\cdots + v_m\wt\K$, with $v_j$ mutually orthogonal. Define $\phi$: $M\to \wt\K^m$: $\phi(\sum_j\lambda_j v_j)=(\lambda_1\norm{v_1},\dots,\lambda_m\norm{v_m})$ ($\lambda_j\in\wt\K$). Then, by the orthogonality,
\[
\Big\Vert\phi\Big(\sum_j\lambda_j v_j\Big)-\phi\Big(\sum_j\mu_j v_j\Big)\Big\Vert^2
=\sum_j\norm{(\lambda_j-\mu_j)v_j}^2
=\Big\Vert{\sum_j(\lambda_j-\mu_j)v_j}\Big\Vert^2
=\Big\Vert{\sum_j\lambda_j v_j-\sum_j\mu_j v_j}\Big\Vert^2\,,
\]
which shows that $\phi$ is well-defined and isometric (hence also injective). It is easy to check that $\phi$ is $\wt\K$-linear.\\
(3) We extend $\phi$: $M\to\wt\K^m$ to a map $\overline M\to\wt\K^m$ by defining $\phi(\lim_n w_n):=\lim_n \phi(w_n)$ ($w_n\in M$). Because $(w_n)_n$ is a Cauchy-sequence, $(\phi(w_n))_n$ is also a Cauchy-sequence in $\wt\K^m$, and hence convergent in $\wt\K^m$. To see that $\phi$ is well-defined, let $\lim_n w_n =\lim_n w'_n$. Then also the interlaced sequence $(w_1, w'_1,\dots, w_n,w'_n, \dots)$ is a Cauchy-sequence. Hence also $(\phi(w_1), \phi(w'_1),\dots, \phi(w_n),\phi(w'_n), \dots)$ is convergent to $\lim_n\phi(w_n)=\lim_n\phi(w'_n)$. It is easy to check that also the extended $\phi$ is linear and isometric. As $\G$ is complete, we find that the image under ${\phi}^{-1}$ of any convergent sequence in $\phi(M)$ (say to $\xi\in\overline{\phi(M)}$) is a Cauchy sequence, and hence convergent to an element $w\in \overline M$. By definition of the extended map $\phi$, $\phi(w)=\xi$. So $\overline {\phi(M)}\subseteq \phi(\overline M)$. The converse inclusion holds by continuity of $\phi$.\\
(4), (5) Let $M$ be closed. By part~1, $M=u_1\wt\K+\cdots+u_m\wt\K$ with $u_j$ mutually orthogonal. Let $w\in\overline{u_1\wt\K}$, so $w=\lim_n \lambda_n u_1$, for some $\lambda_n\in\wt\K$. As $M$ is closed, $w=\sum_j\mu_j u_j$, for some $\mu_j\in\wt\K$. By the continuity of the scalar product, $\inner{w}{u_j}=\lim_n\lambda_n \inner{u_1}{u_j}=0$ for $j>1$. So $0=\inner{w}{u_j}=\sum_k\mu_k\inner{u_k}{u_j}=\mu_j\norm{u_j}^2$, for $j>1$. So also $\norm{\mu_j u_j}^2=0$, for $j>1$ and $w=\mu_1u_1\in u_1\wt\K$. Similarly, $u_j\wt\K$ is closed ($j=1$, \dots, $m$).\\
Conversely, let $M=u_1\wt\K+\cdots+u_m\wt\K$, with $u_j\wt\K$ closed and $u_j$ mutually orthogonal.
By Proposition~\ref{prop_cyclic_module} part 12, we know that $M_j=u_j\wt\K$ are edged. 
Let $v\in\G$. Let $p=\sum P_{M_j}(v)\in M$. Then by orthogonality, $\inner{v-p}{u_j}=\inner{v-P_{M_j}(v)}{u_j}=0$. So, by Corollary \ref{corol_1}$(i)$ it follows that $M$ is closed and edged.\\ 
(6) Let $M= u_1\wt\K +\cdots + u_m\wt\K$. As $\overline N$ is edged and closed, by the linearity of the projection operator $P_{\overline N}$, $\overline N=P_{\overline N}(M)=v_1\wt\K +\cdots + v_m\wt\K$, with $v_j=P_{\overline N}u_j$. By part~4, we may suppose that $\norm{v_j}=e_{S_j}$ for some $S_j\subseteq(0,1]$ and that $v_j$ are mutually orthogonal. So
\[
0=\overline{\inf}_{\sum_j\mu_j v_j\in N}\Big\Vert{v_1-\sum_j\mu_j v_j}\Big\Vert
=\overline{\inf}_{\sum_j\mu_j v_j\in N}\big(\abs{1-\mu_1}^2 e_{S_1}+\sum_{j>1}\norm{\mu_j v_j}^2\big)^{1/2},
\]
so for each $m\in\N$, there exist $\sum_j\mu_jv_j\in N$ with $\abs{1-\mu_1}e_{S_1}\le [(\eps^m)_\eps]$, $\norm{\mu_j v_j}\le[(\eps^m)_\eps]$. For sufficiently large $m$, this implies that $\mu_1$ is invertible w.r.t.\ $S_1$. Let $\lambda_1\mu_1=e_{S_1}$ with $\lambda_1=\lambda_1 e_{S_1}$. Then, $|\mu_1|e_{S_1}\ge (1-[(\eps^m)_\eps])e_{S_1}\ge \frac{1}{2} e_{S_1}$. Hence, $|\lambda_1|\le 2 e_{S_1}$. So for each $m$, there exist $v_1+\sum_{j\ne 1} \mu_j v_j\in N$ with $\norm{\mu_j v_j}\le[(\eps^m)_\eps]$. Similarly, for each $m$, there exist $v_k+\sum_{j\ne k} \mu_j^{(k)} v_j\in N$ with $\norm{\mu_j^{(k)} v_j}\le[(\eps^m)_\eps]$ ($k=1$, \dots, $m$). Then also linear combinations \[\Big(v_1+\sum_{j\ne 1} \mu_j^{(1)} v_j\Big) -\mu_2^{(1)}e_{S_2} \Big(v_2+\sum_{j\ne 2} \mu_j^{(2)} v_j\Big)=(1-\mu_2^{(1)}e_{S_2}\mu_1^{(2)}e_{S_1})v_1+\sum_{j\ne 1,2} \mu'_j v_j\in N,\]
with $\norm{\mu'_j v_j}$ arbitrarily small. As $1-\mu_2^{(1)}e_{S_2}\mu_1^{(2)}e_{S_1}$ is invertible (for $m$ sufficiently large), this also implies that there exist $v_1+\sum_{j\ne 1,2} \mu_j v_j\in N$ with $\norm{\mu_j v_j}$ arbitrarily small, and so on. We conclude that $v_1$, \dots, $v_m$ $\in N$. So $N=\overline N$ is closed and finitely generated.\\
(7) Let $M=u_1\wt\K+\cdots+u_m\wt\K$ with $u_j$ mutually orthogonal. By the normalization property, there exists $v_1\in\G$ with $\norm{u_1}v_1=u_1$. As in Proposition~\ref{prop_cyclic_module} part~9, there exists $S_1\subseteq(0,1]$ such that $\norm{u_1}e_{\co S_1}=0$ and $\norm{v_1}e_{S_1}=e_{S_1}$. As $\inner{u_1}{u_j}=0$ for $j>1$, also $\norm{u_1}\inner{v_1}{u_j}=0$, so by a characterization of zero divisors in $\wt\K$, there exist $S_j\subseteq (0,1]$ such that $\norm{u_1}e_{\co S_j}= 0$ and $\inner{v_1}{u_j}e_{S_j}=0$. Let $w_1=v_1e_{S_1}\cdots e_{S_m}$. Then $\norm{w_1}=e_{S_1}\cdots e_{S_m}$ is idempotent, hence $w_1\wt\K$ is closed by proposition~\ref{prop_cyclic_module} part 6. Further, $u_1=u_1e_{S_1}\cdots e_{S_m}=\norm{u_1}w_1$, so $u_1\wt\K\subseteq w_1\wt\K$. Finally, $\inner{w_1}{u_j}=\inner{v_1}{u_j}e_{S_1}\cdots e_{S_m}=0$, for $j>1$.
Similarly, we find $w_2\in\G$ such that $w_2\wt\K$ is closed, $u_2\wt\K\subseteq w_2\wt\K$ and $\inner{w_2}{w_1}=\inner{w_2}{u_3}=\cdots=\inner{w_2}{u_m}=0$, and so on.\\
(8) By parts~1 and 7, $M=u_1\wt\K +\cdots+ u_m\wt\K$ and there exist $w_j$ with $w_j\wt\K$ closed, $u_j\in w_j\wt\K$ and $w_j$ mutually orthogonal. As $\mathcal H=w_1\wt\K$ is closed, it is itself a Hilbert $\wt\K$-module. 
We show that $u_1\wt\K$ is an edged submodule of $\mathcal H$.\\
So let $\lambda\in\wt\K$. Since $M$ is edged in $\G$, $\overline{\inf}_{v\in M}\Vert{\lambda w_1-v}\Vert$ exists. So by orthogonality,
\[
\overline{\inf}_{\mu_j\in\wt\K}\Big\Vert{\lambda w_1-\sum_j\mu_ju_j}\Big\Vert
=\overline{\inf}_{\mu_j\in\wt\K}\Big(\norm{\lambda w_1 - \mu_1 u_1}^2 + \sum_{j>1}\norm{\mu_ju_j}^2\Big)^{1/2}\\
=\overline{\inf}_{\mu_1\in\wt\K}\norm{\lambda w_1 - \mu_1 u_1}
\]
hence $u_1\wt{\K}$ is edged in $\mathcal{H}$. By Proposition~\ref{prop_cyclic_module} part~8, $u_1\wt\K$ is closed in $\mathcal H$ and by completeness, also in $\G$. Similarly, $u_j\wt\K$ is closed ($j=1$, \dots, $m$). By the fourth assertion of this theorem, $M$ is closed.
\end{proof}

\begin{theorem}
Let $\G_H$ be a Hilbert $\wt\K$-module constructed by means of a Hilbert space $H$. Then a finitely generated submodule $M$ of $\G_H$ is edged iff $M$ is closed iff $M$ is a finite direct sum of mutually orthogonal closed cyclic modules iff $M$ is internal.
\end{theorem}
\begin{proof}
Let $M$ be a finite direct sum of mutually orthogonal closed cyclic modules, so $M=u_1\wt\K + \cdots + u_m\wt\K$ with $u_j$ mutually orthogonal and $\norm{u_j}=e_{S_j}$, for some $S_j\subseteq(0,1]$. We show that $M$ is internal.\\
Fix representatives $(u_{j,\eps})_\eps$ of $(u_j)$. By interleaved Gram-Schmidt at the level of representatives, we may suppose that $\inner{u_{j,\eps}}{u_{k,\eps}}=0$, for $j\ne k$. As $u_j e_{\co S_j}=0$, we may also suppose that $u_{j,\eps}=0$, for each $\eps\in \co S_j$. Let $A_\eps=u_{1,\eps}\K +\cdots + u_{m,\eps}\K$ for each $\eps\in(0,1]$. If $v\in M$, looking at representatives, $v\in [(A_\eps)_\eps]$. Conversely, if $v\in [(A_\eps)_\eps]$, we find $\lambda_{j,\eps}\in\K$ such that, on representatives, $v_\eps=\sum_j\lambda_{j,\eps} u_{j,\eps}$. We may assume that $\lambda_{j,\eps}=0$ for $\eps\in\co S_j$. Then $\inner{v_\eps}{u_{j,\eps}} =\lambda_{j,\eps}\inner{u_{j,\eps}}{u_{j,\eps}}$, so $(\lambda_{j,\eps})_\eps$ are moderate (since $\norm{u_j}$ are invertible w.r.t.\ $S_j$). So $(\lambda_{j,\eps})_\eps$ represent $\lambda_j\in\wt\K$ and $v=\sum_j\lambda_j u_j\in M$.\\
Further, any nonempty internal set in $\G_H$ is edged \cite{OV:07}.\\
Since $\G_H$ has the normalization property, the other equivalences follow by the previous theorem.
\end{proof}

\begin{theorem}
\leavevmode
\begin{itemize}
\item[(i)] Let $M$ be a finitely generated submodule of $\wt\K^d$. Then $M$ is generated by $d$ elements.
\item[(ii)] Let $M$ be a submodule of a Hilbert $\wt\K$-module $\G$ that is generated by $m$ elements. Then any finitely generated submodule of $M$ is generated by $m$ elements.
\end{itemize}
\end{theorem}
\begin{proof}
$(i)$ Let $M=u_1\wt\K+\cdots+u_m\wt\K$ with $m>d$. Applying interleaved Gram-Schmidt at the level of representatives, we can obtain representatives $(u_{j,\eps})_\eps$ of $u_j$ such that for each $\eps$, $\inner{u_{j,\eps}}{u_{k,\eps}}=0$ if $j\ne k$. Define recursively for $j=1$, \dots, $m$
\[S_j=\big\{\eps\in (0,1]: \norm{u_j}_\eps\le\max_{k\ne j}\norm{u_k}_\eps\big\}\setminus \{S_1,\dots, S_{j-1}\}.\]
Then $e_{S_j}e_{S_k}=0$ if $j\ne k$ and $e_{S_1}+\cdots+e_{S_m}=1$. Let $\eps\in S_1$. Should $u_{1,\eps}\ne 0$, then also $u_{j,\eps}\ne 0$, for all $j$. So we would obtain $m>d$ orthogonal (hence linearly independent) elements of $\K^d$, a contradiction. So $u_1e_{S_1}=0$, and $Me_{S_1}=v_1^{(1)}\wt\K+\cdots+v_{m-1}^{(1)}\wt\K$, for $v_j^{(1)}=u_{j+1} e_{S_1}$. Similarly, $Me_{S_k}=v_1^{(k)}\wt\K+\cdots+v_{m-1}^{(k)}\wt\K$, for some $v_j^{(k)}\in\wt\K^d$ satisfying $v_j^{(k)}=v_j^{(k)} e_{S_k}$ ($k=1$, \dots, $m$). Then $v_j=\sum_{k=1}^m v_j^{(k)}$ ($j=1$, \dots, $m-1$) are $m-1$ generators of $M$.\\
$(ii)$ Follows from part~1 and Theorem \ref{thm_Gram-Schmidt}.
\end{proof}

\begin{theorem}\label{thm_edged_orthogonal_sum}
Let $M$, $N$ be edged submodules of a Hilbert $\wt\K$-module $\G$. If $M\perp N$, then $M+N$ is edged and $\overline M + \overline N= \overline{M+N}$.
\end{theorem}
\begin{proof}
First, by the continuity of the scalar product in $\G$, if $M\perp N$ then also $\overline M\perp \overline N$.\\
Let $v\in\G$. For each $u\in \overline M$, \[\inner{v-(P_{\overline M}(v)+P_{\overline N}(v))}{u}=\inner{(v-P_{\overline M}(v))+P_{\overline N}(v)}{u}=0\] by the properties of the $P_{\overline M}$ and the fact that $\overline M\perp \overline N$. Switching roles of $M$ and $N$, we obtain that $v-(P_{\overline M}(v)+P_{\overline N}(v))\in (\overline M+\overline N)^\perp$. As also $P_{\overline M}(v)+P_{\overline N}(v)\in \overline M+\overline N$, it follows from Corollary \ref{corol_1} that $\overline M+ \overline N$ is closed and edged.
As $M+N\subseteq \overline M +\overline N \subseteq \overline{M + N}$ and $\overline M + \overline N$ is closed, $\overline M + \overline N=\overline{M+N}$.  
\end{proof}

\begin{theorem}\label{thm_edged_direct_sum}
Let $M$, $N$ be submodules of a Hilbert $\wt\K$-module $\G$. Let $M$ be closed and finitely generated, and $N$ closed and edged. If $M\cap N=\{0\}$, then $M+N$ is closed and edged.
\end{theorem}
\begin{proof}
We proceed by induction on the number $m$ of generators of $M$.\\
First, let $M=u\wt\K$ be cyclic. By Proposition~\ref{prop_cyclic_module}, we may suppose that $\norm{u}=e_S$, for some $S\subseteq(0,1]$. Now suppose that $\norm{u-P_N(u)}$ is not invertible w.r.t.\ $S$. Then there exists $T\subseteq S$ with $0\in \overline T$ such that $\norm{ue_T-P_N(u)e_T}=\norm{u-P_N(u)}e_T=0$, so $ue_T=P_N(u)e_T\in M\cap N$, and $\norm{ue_T}=e_Se_T=e_T\ne 0$, which contradicts $M\cap N=\{0\}$. As $0\le\norm{P_N(u)}\le \norm{u}$, also $\norm{u-P_N(u)}e_{\co S}=0$, and $M'=(u-P_N(u))\wt\K$ is also closed, hence edged by Proposition~\ref{prop_cyclic_module} parts 6 and 12. Since $M'\perp N$, by Theorem \ref{thm_edged_orthogonal_sum}, $M+N=M'+N$ is closed and edged.\\
Now let $M$ be generated by $m$ elements. By Theorem~\ref{thm_Gram-Schmidt}, $M$ is a direct sum of a closed cyclic module $M_1$ and a closed module $M_2$ generated by $m-1$ elements. By induction, as $M_2\cap N=\{0\}$, $M_2+N$ is closed and edged. As also $M_1+(M_2+N)$ is a direct sum, then by the first part of the proof $M+N=M_1+(M_2+N)$ is closed and edged.
\end{proof}

\section{A Riesz-representation theorem for continuous $\wt{\C}$-linear functionals on $\G$}
\label{section_riesz}
In this section we consider Hilbert $\wt{\C}$-modules with the normalization property and we prove a Riesz representation theorem for the corresponding continuous $\wt{\C}$-linear functionals.

\begin{theorem}
\label{theo_riesz}
Let $\G$ be a Hilbert $\wt{\C}$-module with the normalization property and $T$ a continuous $\wt{\C}$-linear functional on $\G$. The following assertions are equivalent:
\begin{itemize}
\item[(i)] there exists a closed edged $\wt{\C}$-submodule $M$ of $\Ker T$ and a subset $S$ of $(0,1]$ such that 
\begin{itemize}
  \item[(a)] there exists $u_1\in M^{\perp}$ with $\Vert u_1\Vert=e_S$;
  \item[(b)] $\Vert u\Vert= e_S\Vert u\Vert$ for all $u\in M^\perp$;
  \item[(c)] $T(u)v-T(v)u\in M$ for all $u,v\in M^\perp$;
\end{itemize}
\item[(ii)] there exists a closed, cyclic (and hence edged) $\wt{\C}$-submodule $N$ of $\G$ such that $N^\perp\subseteq \Ker T$;
\item[(iii)] there exists a unique $c\in\G$ such that $T(u)=\inner{u}{c}$.
\end{itemize}
\end{theorem}
\begin{proof}
$(i)\Rightarrow (iii)$ Let $u_1\in M^\perp$ satisfying the condition $(a)$ and $u\in M^\perp$. From $(c)$ we get that $T(u)u_1-T(u_1)u\in M$ and thus $T(u)e_S=T(u)\Vert u_1\Vert^2=T(u_1)\inner{u}{u_1}$. Since $T$ is continuous there exists $C>0$ such that $\vert T(u)\vert\le C\Vert u\Vert$. It follows that $\vert T(u)\vert e_{\compl{S}}\le C\Vert u\Vert e_{\compl{S}}=0$ because of property $(b)$. So, $T(u)=T(u)e_S=\inner{u}{c}$, where $c=\overline{T(u_1)}u_1$. 

Now let $u\in \G$. By Corollary \ref{corol_1} we know that $u=(u-P_M(u))+P_M(u)$, where $u-P_M(u)\in M^\perp$ and $P_M(u)\in M$. Since $T(u)=T(u-P_M(u))$ by the previous case we have that $T(u)=\inner{u-P_M(u)}{c}=\inner{u}{c}$.

$(iii)\Rightarrow (ii)$ By the normalization property and the assertions 9. and 12. of Proposition \ref{prop_cyclic_module} we know that there exists a closed cyclic and edged $\wt{\C}$-module $N$ such that $c\wt{\C}\subseteq N$. Thus, $N^\perp\subseteq (c\wt{\C})^\perp=\Ker T$.  

$(ii)\Rightarrow (i)$ From the assertion 7. of Proposition \ref{prop_cyclic_module} we have that $N$ is generated by an element $w\in\G$ such that $\Vert w\Vert=e_S$, for some $S\subseteq(0,1]$. Let us define $M=N^\perp$. By Corollary \ref{corol_1} $M$ is closed and edged. $(a)$ and $(b)$ are straightforward and $(c)$ follows from the fact that $M^\perp$ is cyclic. 
\end{proof}
\begin{remark}
Note that $(i)\Rightarrow(iii)$ is valid without assuming the normalization property on $\G$. In the first assertion of Theorem \ref{theo_riesz} we assume the existence of an edged and closed $\wt{\C}$-submodule $M$ contained in $\Ker T$ because in general the kernel of a continuous $\wt{\C}$-linear functional is not edged. Indeed, taking $\beta$ as in \ref{ex_edged_not_closed} and the functional $T:\wt{\C}\to\wt{\C}:z\to \beta z$ we have that $\Ker{T}=\Ann(\beta)$ is not edged \cite{Vernaeve:07a}.
\end{remark}

The following example shows that (at least under some set-theoretic assumption) there are continuous $\wt\C$-linear functionals on Hilbert $\wt\C$-modules for which the Riesz representation theorem does not hold.
\begin{example}
Under the assumption that $2^{\aleph_0}<2^{\aleph_1}$ (e.g., if one assumes the continuum hypothesis), by \cite{Vernaeve:07a}, there exists a submodule (=ideal) $M$ of $\wt\C$ and a continuous $\wt\C$-linear map $T$: $M\to\wt\C$ that cannot be extended to a $\wt\C$-linear map $\wt\C\to\wt\C$. Let $\G=\overline M$ (the topological closure of $M$ in $\wt\C$). Then $M$ is a Hilbert $\wt\C$-module as a closed submodule of a Hilbert $\wt\C$-module. By continuity, $T$ can be uniquely extended to a continuous $\wt{\C}$-linear map $\wt T$: $\G\to\wt\C$. Suppose that there exists $c\in\G$ such that $\wt T(u)=\inner{u}{c}$. Then the $\wt\C$-linear map $\wt\C\to\wt\C$: $u\mapsto \inner{u}{c}$ would be an extension of $T$, a contradiction.
\end{example}

\begin{proposition}
\label{riesz_int}
Let $H$ be a Hilbert space and $\G_H$ the corresponding Hilbert $\wt{\C}$-module. A continuous $\wt{\C}$-linear functional $T$ on $\G_H$ is basic if and only if it fulfills the equivalent properties of the previous theorem. 
\end{proposition} 
\begin{proof}
Apply the Riesz theorem at the level of representatives, noting that $T_\eps(u)=\inner{u}{c_\eps}$ with $\Vert c_\eps\Vert=\Vert T_\eps\Vert$.
\end{proof}
Conjecture: there exists a Hilbert space $H$ (necessarily infinitely dimensional) and a continuous $\wt{\C}$-linear functional that it is not basic.

 
We now investigate the structural properties of continuous $\wt{\C}$-sesquilinear forms on Hilbert $\wt{\C}$-modules by making use of the previous representation theorem.
\begin{theorem}
\label{theo_riesz_b}
Let $\G$ and $\mH$ be Hilbert $\wt{\C}$-modules with $\mH$ satisfying the normalization property and $a:\G\times\mH\to\wt{\C}$ be a continuous $\wt{\C}$-sesquilinear form. The following assertions are equivalent:
\begin{itemize}
\item[(i)] for all $u\in\G$ there exists a closed and cyclic $\wt{\C}$-submodule $N_u$ of $\mH$ such that $N_u^\perp\subseteq\{v\in \mH:\, a(u,v)=0\}$;
\item[(ii)] there exists a unique continuous $\wt{\C}$-linear map $T:\G\to\mH$ such that $a(u,v)=\inner{Tu}{v}$ for all $u\in\G$ and $v\in\mH$.
\end{itemize}
\end{theorem} 
\begin{proof}
$(i)\Rightarrow (ii)$ Let $u\in\G$. We consider the continuous $\wt{\C}$-linear functional $a_u:\mH\to\wt{\C}:v\to\overline{a(u,v)}$. Since $\Ker\, a_u=\{v\in \mH:\, a(u,v)=0\}$ contains the orthogonal complement of a closed and cyclic $\wt{\C}$-submodule $N_u$ of $\mH$, by Theorem \ref{theo_riesz} there exists a unique $c\in\mH$ such that $\overline{a(u,v)}=\inner{v}{c}$. We define $T:\G\to\mH:u\to c$. By construction, $a(u,v)=\inner{Tu}{v}$. We leave to the reader to check that the map $T$ is $\wt{\C}$-linear. By definition of the operator $T$ we have that 
\beq
\label{cont_T}
\Vert T(u)\Vert^2=\inner{Tu}{Tu}=a(u,Tu)\le C\Vert u\Vert \Vert Tu\Vert,
\eeq
where the constant $C\in\wt{\R}$ comes from the continuity of $a$. Applying Lemma \ref{lemma_scratch} to \eqref{cont_T} we have that $\Vert T(u)\Vert\le C\Vert u\Vert$ for all $u$. This shows that $T$ is continuous.\\
$(ii)\Rightarrow (i)$ Let us fix $u\in\G$. Since $v\mapsto \overline{a(u,v)}=\inner{v}{Tu}$ is a continuous $\wt{\C}$-linear functional on $\mH$ satisying the assertion $(iii)$ of Theorem \ref{theo_riesz} we find a subset $N_u$ as desired.
\end{proof}
\begin{proposition}
\label{prop_riesz_b_int}
Let $H$ and $K$ be Hilbert spaces and $a$ be a basic $\wt{\C}$-sesquilinear form on $\G_H\times\G_K$. Then, the hypotheses of Theorem \ref{theo_riesz_b} are satisfied. Moreover, the map $T:\G_H\to\G_K$ such that $a(u,v)=\inner{Tu}{v}$ is basic. 
\end{proposition}
\begin{proof}
By Proposition \ref{prop_normal} the $\wt{\C}$-module $\G_K$ has the normalization property. If $a$ is basic then for any fixed $u\in\G_H$ the $\wt{\C}$-linear functional $\G_K\to\wt{\C}:v\to\overline{a(u,v)}$ is basic too. Hence, from Theorem \ref{theo_riesz} there exists a closed and cyclic $\wt{\C}$-submodule $N_u$ of $\G_K$ such that $N_u^\perp\subseteq\{v\in \G_K:\, a(u,v)=0\}$. It remains to prove that the continuous $\wt{\C}$-linear map $T:\G_H\to\G_K$, that we know to exist from Theorem \ref{theo_riesz_b}, has a basic structure. Let us take a net $(a_\eps)_\eps$ representing the $\wt{\C}$-sesquilinear form $a$. By fixing $u\in H$ we obtain from the continuity of $a_\eps$ that there exist a net $(c_\eps)_\eps$ of elements of $K$ and a net $t_\eps(u)=c_\eps$ of linear maps from $H$ to $K$ such that $\overline{a_\eps(u,v)}=\inner{v}{c_\eps}$ for all $v\in K$. Since for some $N\in\N$ and $\eta\in(0,1]$ the inequality
\[
\Vert t_\eps(u)\Vert^2= a_\eps(u,t_\eps(u))\le \eps^{-N}\Vert u\Vert\, \Vert t_\eps(u)\Vert
\]
holds for all $u\in H$ and $\eps\in(0,\eta]$, we obtain that $(t_\eps)_\eps$ defines a basic map $T':\G_H\to \G_K$ such that $a(u,v)=\inner{T'u}{v}$. By Theorem \ref{theo_riesz_b} there exists a unique continuous $\wt{\C}$-linear map from $\G_H$ to $\G_K$ having this property. It follows that $T'=T$ and that $T$ is basic.
\end{proof}

\section{Continuous $\wt{\C}$-linear operators on a Hilbert $\wt{\C}$-module}
\label{section_cont}
In this section we focus on continuous $\wt{\C}$-linear operators acting on a Hilbert $\wt{\C}$-module. In particular we deal with isometric, unitary, self-adjoint and projection operators obtaining an interesting characterization for the projection operators.
\subsection{Adjoint}
\begin{definition}
\label{def_adjoints}
Let $\G$ and $\mH$ be Hilbert $\wt{\C}$-modules and $T:\G\to\mH$ a continuous $\wt{\C}$-linear map. A continuous $\wt{\C}$-linear operator $T^\ast:\mH\to\G$ is called adjoint of $T$ if 
\beq
\label{eq_T_ast}
\inner{Tu}{v}=\inner{u}{T^\ast v}
\eeq
for all $u\in\G$ and $v\in\mH$.
\end{definition}
Note that if there exists an operator $T^\ast$ satisfying \eqref{eq_T_ast} then it is unique. 

The following proposition characterizes the existence of the adjoint $T^\ast$ under suitable hypotheses on the spaces $\G$ and $\mH$.
\begin{proposition}
\label{prop_adjoint}
Let $\G$ and $\mH$ be Hilbert $\wt{\C}$-modules with $\G$ satisfying the normalization property and $T:\G\to\mH$ be a continuous $\wt{\C}$-linear map. The adjoint $T^\ast:\mH \to \G$ exists if and only if for all $v\in\mH$ there exists a closed and cyclic $\wt{\C}$-submodule $N_v$ of $\G$ such that $N_v^\perp\subseteq\{u\in \G:\, \inner{v}{Tu}=0\}$.  
\end{proposition}
\begin{proof}
The proof is clear by applying Theorem \ref{theo_riesz_b} to the continuous $\wt{\C}$-sesquilinear form $a:\mH\times\G\to\wt{\C}:(v,u)\to\inner{v}{Tu}$.  
\end{proof}

\begin{proposition}
\label{prop_adjoint_int}
If $H$ and $K$ are Hilbert spaces and $T$ is a basic $\wt{\C}$-linear map from $\G_H$ to $\G_K$ then the hypotheses of the previous proposition are fulfilled. In particular the operator $T^\ast:\G_K\to\G_H$ is basic. 
\end{proposition}
\begin{proof}
It suffices to observe that the $\wt{\C}$-sesquilinear form $\inner{v}{Tu}$ is basic and to apply Proposition \ref{prop_riesz_b_int}.
\end{proof}


\begin{proposition}
\label{prop_adjoint_1}
Let $\G$ and $\mH$ be Hilbert $\wt{\C}$-modules and $S,T:\G\to\mH$ continuous $\wt{\C}$-linear maps having an adjoint.  The following properties hold:
\begin{itemize}
\item[(i)] $(S+T)^\ast=S^\ast+T^\ast$;
\item[(ii)] $(\lambda T)^\ast=\overline{\lambda}\,T^\ast$ for all $\lambda\in\wt{\C}$;
\item[(iii)] $\inner{T^\ast(v)}{u}=\inner{v}{Tu}$ for all $u\in\G$ and $v\in\mH$;
\item[(iv)] $T^{\ast\ast}=T$;
\item[(v)] $T^\ast T=0$ if and only if $T=0$;
\item[(vi)] $(ST)^\ast=T^\ast S^\ast$;
\item[(vii)] if $M\subseteq\G$, $N\subseteq\mH$ and $T(M)\subseteq N$ then $T^\ast(N^\perp)\subseteq M^\perp$;
\item[(viii)] if $M\subseteq\G$ and $N$ is a closed and edged $\wt{\C}$-submodule of $\mH$, then $T(M)\subseteq N$ if and only if $T^\ast(N^\perp)\subseteq M^\perp$.
\end{itemize} 
\end{proposition}
\begin{proof}
We omit the proof of the first seven assertions of the proposition because they are elementary.\\
$(viii)$ From assertion $(vii)$ we have that $T(M)\subseteq N$ implies $T^\ast(N^\perp)\subseteq M^\perp$. Conversely, assume that $T^\ast(N^\perp)\subseteq M^\perp$ and apply $(vii)$ to $T^\ast$. It follows that $(T^\ast)^\ast (M^{\perp\perp})\subseteq N^{\perp\perp}$. By $(iv)$ we can write $T(M)\subseteq T(M^{\perp\perp})\subseteq N^{\perp\perp}$ and from Corollary \ref{corol_1}$(iii)$ we have that $T(M)\subseteq N$.
\end{proof}


\begin{proposition}
\label{prop_adjoint_2}
Let $\G$ and $\mH$ be Hilbert $\wt{\C}$-modules and $T:\G\to\mH$ be a continuous $\wt{\C}$-linear map. Assume that the adjoint of $T$ exists. The following equalities hold:
\begin{itemize}
\item[(i)] $\Ker{T}=(T^\ast (\mH))^\perp$;
\item[(ii)] $\Ker{T^\ast}=(T(\G))^\perp$;
\item[(iii)] if $T^\ast(\mH)$ is a closed and edged $\wt{\C}$-submodule of $\G$ then $(Ker{T})^\perp=T^\ast(\mH)$;
\item[(iv)] if $T(\G)$ is a closed and edged $\wt{\C}$-submodule of $\mH$ then $(Ker{T^\ast})^\perp=T(\G)$.
\end{itemize}
\end{proposition} 
\begin{proof}
An application of Proposition \ref{prop_adjoint_1}$(vii)$ and $(iii)$ to $T$ and $T^\ast$ yields
\beq
\label{i}
T^\ast(\mH)\subseteq(\Ker{T})^\perp,
\eeq
\beq
\label{i_ast}
T(\G)\subseteq(\Ker{T^\ast})^\perp,
\eeq
\beq
\label{ii}
(T(\G))^\perp\subseteq\Ker{T^\ast},
\eeq
\beq
\label{ii_ast}
(T^\ast(\mH))^\perp\subseteq\Ker{T}.
\eeq
\eqref{i} combined with \eqref{ii_ast} entails the first assertion while \eqref{ii} combined with \eqref{i_ast} entails the second assertion. The assertions $(iii)$ and $(iv)$ are obtained from $(i)$ and $(ii)$ respectively making use of Corollary \ref{corol_1}$(iii)$.
\end{proof}

\subsection{Isometric, unitary, self-adjoint and projection operators}

\begin{definition}
\label{def_isom}
Let $\G$ and $\mH$ be Hilbert $\wt{\C}$-modules. A continuous $\wt{\C}$-linear operator $T:\G\to\mH$ is said to be isometric if $\Vert Tu\Vert=\Vert u\Vert$ for all $u\in\G$.
\end{definition}
\begin{lemma}
\label{lemma_seq}
Any $\wt{\C}$-sesquilinear form $a:\G\times\G\to\wt{\C}$ is determined by its values on the diagonal, in the sense that \[
a(u,v)=\frac{1}{4}\big[a(u+v,u+v)-a(u-v,u-v)+ia(u+iv,u+iv)-ia(u-iv,u-iv)\big]
\]
for all $u,v\in\G$.	
\end{lemma}
\begin{proposition}
\label{prop_isom}
Let $\G$ and $\mH$ be Hilbert $\wt{\C}$-modules and $T:\G\to\mH$ a continuous $\wt{\C}$-linear operator with an adjoint. The following assertions are equivalent:
\begin{itemize}
\item[(i)] $T$ is isometric;
\item[(ii)] $T^\ast T=I$;
\item[(iii)] $\inner{Tu}{Tv}=\inner{u}{v}$ for all $u,v\in\G$.
\end{itemize}
\end{proposition}
\begin{proof}
$(i)\Rightarrow(ii)$ By definition of isometric operator and adjoint operator we have that $\inner{u}{u}=\inner{Tu}{Tu}=\inner{T^\ast Tu}{u}$. Hence, $\inner{T^\ast Tu-Iu}{u}=0$. Since the form $(u,v)\to\inner{T^\ast Tu-Iu}{v}$ is $\wt{\C}$-sesquilinear, from Lemma \ref{lemma_seq} we conclude that $\inner{T^\ast Tu-Iu}{v}=0$ for all $u,v$, that is $T^\ast T=I$. The implications $(ii)\Rightarrow (iii)$ and $(iii)\Rightarrow (i)$ are immediate.
\end{proof}
More generally, from Lemma \ref{lemma_seq} we have that $(i)$ is equivalent to $(iii)$ for any continuous $\wt{\C}$-linear operator $T:\G\to\mH$ even when the adjoint $T^\ast$ does not exist.
\begin{proposition}
\label{prop_range}
The range of an isometric operator $T:\G\to\mH$ between Hilbert $\wt{\C}$-modules is a closed $\wt{\C}$-submodule of $\mH$.
\end{proposition}
\begin{proof}
Let $v\in\overline{T(\G)}$. There exists a sequence $(u_n)_n$ of elements of $\G$ such that $Tu_n\to v$ in $\mH$. By definition of isometric operator we obtain that $(u_n)_n$ is a Cauchy sequence in $\G$ and therefore it is convergent to some $u\in\G$. It follows that $v=Tu$.
\end{proof}

\begin{proposition}
\label{prop_ad_Hans}
Let $T$: $\G\to\mH$ be a continuous $\wt\C$-linear operator between Hilbert $\wt{\C}$-modules. If $T^*$ exists and there exists a continuous $\wt\C$-linear operator $S$: $\mH\to\G$ such that $T^* T S= T^*$, then $T(\G)$ is closed and edged. Moreover, $P_{T(\G)}=TS$.
\end{proposition}
\begin{proof}
Let $u\in \mH$. Then $T^*(u- T S u)=T^*u-T^*u=0$, so by Proposition \ref{prop_adjoint_2}, $u= (u- T S u) + T S u\in \Ker{T^*} + T(\G)=T(\G)^\perp + T(\G)$. By Corollary \ref{corol_1}, $T(\G)$ is closed and edged and $P_{T(\G)}=TS$.
\end{proof}
\begin{corollary}
\label{corol_isom_edged}
Let $\G$ and $\mH$ be Hilbert $\wt{\C}$-modules and $T$ an isometric operator with adjoint. Then, $T(\G)$ is closed and edged.
\end{corollary}
\begin{proof}
Apply Proposition \ref{prop_ad_Hans} to $T$ with $S=T^\ast$.
\end{proof}

\begin{example}
\label{ex_isom}
A basic operator $T:\G_H\to\G_H$ given by a net of isometric operators $(T_\eps)_\eps$ on $H$ is clearly isometric on $\G_H$. In particular by the corollary above, the range $T(\G_H)$ is a closed and edged $\wt{\C}$-submodule of $\G_H$.  
\end{example}

\begin{definition}
\label{def_unitary}
Let $\G$ be a Hilbert $\wt{\C}$-module and $T:\G\to\G$ a continuous $\wt{\C}$-linear operator with an adjoint. $T$ is unitary if and only if $T^\ast T=TT^\ast=I$.
\end{definition}
\begin{proposition}
\label{prop_unitary}
Let $\G$ be a Hilbert $\wt{\C}$-module and $T:\G\to\G$ a continuous $\wt{\C}$-linear operator with an adjoint. The following conditions are equivalent:
\begin{itemize}
\item[(i)] $T$ is unitary;
\item[(ii)] $T^\ast$ is unitary;
\item[(iii)] $T$ and $T^\ast$ are isometric;
\item[(iv)] $T$ is isometric and $T^\ast$ is injective;
\item[(v)]  $T$ is isometric and surjective;
\item[(vi)] $T$ is bijective and $T^{-1}=T^\ast$.
\end{itemize}
\end{proposition}
\begin{proof}
By Proposition \ref{prop_isom} it is clear that $(i)$, $(ii)$ and $(iii)$ are equivalent. Since any isometric operator is injective we have that $(iii)$ implies $(iv)$. 

$(iv)\Rightarrow (v)$ By Corollary \ref{corol_isom_edged} we know that $T(\G)$ is a closed and edged $\wt{\C}$-submodule of $\G$ and that $\Ker{T^\ast}=\{0\}$. Hence, by Proposition \ref{prop_adjoint_2}$(iv)$ we have that $\{0\}^\perp=(\Ker{T^\ast})^\perp=T(\G)$, which means that $\G=T(\G)$.

$(v)\Rightarrow(vi)$ $T$ is isometric and surjective. Thus, it is bijective. Moreover, $T^\ast T=I=TT^{-1}$. Thus, $T^\ast=T^\ast(TT^{-1})=(T^\ast T)T^{-1}=T^{-1}$. The fact that $(vi)$ implies $(i)$ is clear.
\end{proof}

\begin{definition}
\label{def_self_adjoint}
\label{def_self_adj}
Let $\G$ be a Hilbert $\wt{\C}$-module and $T:\G\to\G$ a continuous and $\wt{\C}$-linear operator. $T$ is said to be self-adjoint if $\inner{Tu}{v}=\inner{u}{Tv}$ for all $u,v\in\G$.
\end{definition}
If $T$ is self-adjoint then the adjoint operator $T^\ast$ exists and coincides with $T$.
\begin{proposition}
\label{prop_self_adj_1}
The following conditions are equivalent:
\begin{itemize}
\item[(i)] $T$ is self-adjoint;
\item[(ii)] $\inner{Tu}{u}=\inner{u}{Tu}$ for all $u\in\G$;
\item[(iii)] $\inner{Tu}{u}\in\wt{\R}$ for all $u\in\G$.
\end{itemize}
\end{proposition}
\begin{proof}
We prove that $(iii)$ implies $(i)$. By Lemma \ref{lemma_seq} we can write $\inner{Tu}{v}$ as 
\[
\frac{1}{4}\big[\inner{T(u+v)}{u+v}-\inner{T(u-v)}{u-v}\big]+i\frac{1}{4}\big[\inner{T(u+iv)}{u+iv}-\inner{T(u-iv)}{u-iv}\big].
\]
Since each scalar product belongs to $\wt{\R}$ and therefore $\inner{Tw}{w}=\inner{w}{Tw}$ for all $w\in\G$, we obtain that $\inner{Tu}{v}=\inner{u}{Tv}$.
\end{proof}
We leave to the reader to prove the following proposition.
\begin{proposition}
\label{prop_self_adj_2}
Let $\G$ be a Hilbert $\wt{\C}$-module and let $S,T:\G\to\G$ be continuous $\wt{\C}$-linear operators. 
\begin{itemize}
\item[(i)] If $S,T$ are self-adjoint then $S+T$ is self-adjoint;
\item[(ii)] if $T$ is self-adjoint and $\alpha\in\wt{\R}$ then $\alpha T$ is self-adjoint;
\item[(iii)] if $T^\ast$ exists then $T^\ast T$ and $T+T^\ast$ are self-adjoint; 
\item[(iv)] if $S$ and $T$ are self-adjoint then $ST$ is self-adjoint if and only if $ST=TS$.
\end{itemize}
\end{proposition}
Note that Proposition \ref{prop_adjoint_2} can be stated for self-adjoint operators on a Hilbert $\wt{\C}$-module $\G$ by replacing $T^\ast$ with $T$.
\begin{example}
There are self-adjoint operators whose range is not edged. Indeed, let $\beta\in\wt{\R}$ be as in Example \ref{ex_edged_not_closed} and $T:\wt{\C}\to\wt{\C}:u\to\beta u$. $T$ is self-adjoint but $T(\wt{\C})=\beta\wt{\C}$ is not edged \cite{Vernaeve:07a}.
\end{example}
\begin{definition}
\label{def_projection}
A continuous $\wt{\C}$-linear operator $T:\G\to\G$ on a Hilbert $\wt{\C}$-module $\G$ is called a projection if it is self-adjoint and $T=TT$.
\end{definition}
Note that when $M$ is a closed and edged $\wt{\C}$-submodule of $\G$ then the corresponding $P_M$ is a projection in the sense of Definition \ref{def_projection}. Indeed, by Proposition \ref{prop_cont}$(iii)$ $P_M$ is idempotent and combining Corollary \ref{corol_1}$(v)$ with Proposition \ref{prop_self_adj_1} we have that $P_M$ is self-adjoint.
We prove the converse in the following proposition.  
\begin{proposition}
If $T$ is a projection then $T(\G)$ is an edged and closed $\wt{\C}$-submodule of $\G$ and $T=P_{T(\G)}$.
\end{proposition}
\begin{proof}
We apply Proposition \ref{prop_ad_Hans} with $S=I$.
\end{proof}


\section{Lax-Milgram theorem for Hilbert $\wt{\C}$-modules}
\label{section_lax}
As in the classical theory of Hilbert spaces we prove that for any $f\in\G$ the problem
\[
a(u,v)=\inner{v}{f},\qquad\qquad\text{for all $v\in\G$},
\]
can be uniquely solved in $\G$ under suitable hypotheses on the $\wt{\C}$-sesquilinear form $a$. In this way, we obtain a Lax-Milgram theorem for Hilbert $\wt{\C}$-modules.
\begin{definition}
\label{def_coercivity}
A $\wt{\C}$-sesquilinear form $a$ on a Hilbert $\wt{\C}$-module $\G$ is coercive if there exists an invertible $\alpha\in\wt{\R}$ with $\alpha\ge 0$ such that
\beq
\label{first_ineq}
a(u,u)\ge \alpha \Vert{u}\Vert^2
\eeq
for all $u\in\G$. 
\end{definition}
\begin{theorem}
\label{LM_theo}
Let $\G$ be a Hilbert $\wt{\C}$-module and $g$ a $\wt{\C}$-linear continuous map on $\G$ such that $g(\G)$ is edged. Let $a$ be the $\wt{\C}$-sesquilinear form on $\G$ defined by $a(u,v)=\inner{u}{g(v)}$. If $a$ is coercive then for all $f\in\G$ there exists a unique $u\in\G$ such that $$a(v,u)=\inner{v}{f}$$ for all $v\in\G$.
\end{theorem}
\begin{proof}
We want to prove that the map $g$ is an isomorphism on $\G$. We begin by observing that the coercivity of $a$ combined with the Cauchy-Schwarz inequality yields for all $u\in \G$
\beq
\label{g_ineq_1}
\alpha \Vert u\Vert^2\le |a(u,u)|=|\inner{u}{g(u)}|\le \Vert u\Vert \Vert g(u)\Vert.
\eeq
By applying Lemma \ref{lemma_scratch} it follows that 
\beq
\label{g_ineq}
\alpha\Vert u\Vert\le\Vert g(u)\Vert.
\eeq
This means that $g$ is an isomorphism of $\G$ onto $g(\G)$. It remains to prove that $g$ is surjective. The $\wt{\C}$-submodule $g(\G)$ is closed. Indeed, if $g(u_n)\to v\in\G$ then from \eqref{g_ineq} we have that $(u_n)_n$ is a Cauchy sequence in $\G$ converging to some $u\in\G$. Since $g$ is continuous we conclude that $v=g(u)$. In addition, $g(\G)$ is edged by assumption and \eqref{g_ineq_1} entails $g(\G)^\perp=\{0\}$. Hence, by Corollary \ref{corol_1}, $g(\G)$ coincides with $\G$. Let now $f\in\G$. We have proved that there exists a unique $u\in\G$ such that $f=g(u)$. Thus, $a(v,u)=\inner{v}{g(u)}=\inner{v}{f}$ for all $v\in\G$. 
\end{proof}
Note that when $C$ is a subspace of $H$ then the corresponding space $\G_C$ of generalized functions based on $C$ is canonically embedded into $\G_H$.
\begin{lemma}
\label{lemma_coer}
Let $H$ be a Hilbert space, $C$ a subspace of $H$, $\alpha\in\wt{\R}$ positive and invertible, and $a$ be a basic $\wt{\C}$-sesquilinear form on $\G_H$. The following assertions are equivalent:
\begin{itemize}
\item[(i)] $a(u,u)\ge\alpha\Vert u\Vert^2$ for all $u\in\G_C$;
\item[(ii)] for all representatives $(a_\eps)_\eps$ of $a$ and $(\alpha_\eps)_\eps$ of $\alpha$ and for all $q\in\N$ there exists $\eta\in(0,1]$ such that $$a_\eps(u,u)\ge(\alpha_\eps-\eps^q)\Vert u\Vert^2$$ for all $u\in C$ and $\eps\in(0,\eta]$.
\item[(iii)] for all representatives $(a_\eps)_\eps$ of $a$ there exists a representative $(\alpha_\eps)_\eps$ of $\alpha$ and a constant $\eta\in(0,1]$ such that $$a_\eps(u,u)\ge\alpha_\eps\Vert u\Vert^2$$ for all $u\in C$ and $\eps\in(0,\eta]$.
\end{itemize}
\end{lemma}
\begin{proof}
It is clear that $(iii)$ implies $(i)$. We begin by proving that $(ii)$ implies $(iii)$. Let $(\alpha'_\eps)_\eps$ be a representative of $\alpha$. Assume that there exists a decreasing sequence $(\eta_q)_q$ tending to $0$ such that $a_\eps(u,u)\ge(\alpha'_\eps-\eps^q)\Vert u\Vert^2$ for all $u\in C$ and $\eps\in(0,\eta_q]$. The net $n_\eps=\eps^q$ for $\eps\in(\eta^{q+1},\eta^q]$ is negligible and therefore $\alpha_\eps=\alpha'_\eps-n_\eps$ satisfies the inequality of the assertion $(iii)$ on the interval $(0,\eta_0]$.\\
Note that the first assertion is equivalent to claim that $e_S a(u,u)\ge \alpha\Vert u\Vert^2 e_S$ for all $S\subseteq(0,1]$. We want now to prove that if 
\beq
\label{neg_2}
\exists (a_\eps)_\eps\, \exists(\alpha_\eps)_\eps\, \exists q\in\N\, \forall\eta\in(0,1]\, \exists\eps\in(0,\eta]\, \exists u\in C\qquad a_\eps(u,u)<(\alpha_\eps-\eps^q)\Vert u\Vert^2
\eeq
then we can find $S\subseteq(0,1]$ and $u\in\G_C$ such that $e_S a(u,u)< \alpha\Vert u\Vert^2 e_S$. From \eqref{neg_2} it follows that there exists a decreasing sequence $(\eps_k)_k\subseteq(0,1]$ converging to $0$ and a sequence $(u_{\eps_k})_k$ of elements of $C$ with norm $1$ such that 
\[
a_{\eps_k}(u_{\eps_k},u_{\eps_k})<\alpha_{\eps_k}-\eps_k^q.
\]
Let us fix $x\in C$ with $\Vert x\Vert=1$. The net $v_\eps=u_{\eps_k}$ when $\eps=\eps_k$ and $v_\eps=x$ otherwise generates an element $v=[(v_\eps)_\eps]$ of $\G_C$ with $\wt{\R}$-norm $1$. Let now $S=\{\eps_k:\, k\in\N\}$. By construction we have that
\[
e_S a(v,v)=[(\chi_S a_{\eps_k}(u_{\eps_k},u_{\eps_k}))_\eps]\le e_S(\alpha -[(\eps^q)_\eps])< e_S\alpha\Vert v\Vert^2.
\]
This contradicts assertion $(i)$.
\end{proof}

\begin{proposition}
\label{LM_int}
Let $H$ be a Hilbert space, $a$ be a basic coercive $\wt{\C}$-sesquilinear form on $\G_H$ and $f$ a basic functional on $\G_H$. Then there exists a unique $u\in\G_H$ such that $a(v,u)=f(v)$ for all $v\in\G_H$.
\end{proposition}
\begin{proof}
By applying Proposition \ref{prop_riesz_b_int} to the $\wt{\C}$-sesquilinear form $b(u,v):=\overline{a(v,u)}$, there exists a basic map $g:\G_H\to\G_H$ such that $a(u,v)=\inner{u}{g(v)}$. In order to apply Theorem \ref{LM_theo} it remains to prove that $g(\G_H)$ is edged. By the continuity of $g$ and the inequality \eqref{g_ineq}, we find by Proposition \ref{prop_cont_basic} and Lemma \ref{lemma_coer} $C=[(C_\eps)_\eps]\in\wt\R$ and an invertible $\alpha=[(\alpha_\eps)_\eps]\in\wt\R$, $\alpha\ge 0$ for which
\beq\label{eq_equiv_norm}
\alpha_\eps\norm{u}\le \norm{g_\eps(u)}\le C_\eps\norm{u}, \forall u\in H, \forall\eps\le\eta.
\eeq
Let us call $H_\eps$ the Hilbert space $H$ provided with the scalar product $\inner{u}{v}_\eps:=\inner{g_\eps(u)}{g_\eps(v)}$ and consider the Hilbert $\wt\C$-module $\G=\M_{(H_\eps)_\eps}/\Neg_{(H_\eps)_\eps}$ as in Proposition \ref{prop_ex_3}. By equation \eqref{eq_equiv_norm}, a net $(u_\eps)_\eps$ of elements of $H$ is moderate (resp.\ negligible) in $\G_H$ iff it is moderate (resp.\ negligible) in $\G$. Hence the map $\tilde g$: $\G\to\G_H$: $\tilde g([(u_\eps)_\eps])=[g_\eps(u_\eps)_\eps]$ is a well-defined isometric $\wt\C$-linear operator with $\tilde g(\G)=g(\G_H)$. Let $\tilde g_\eps$: $H_\eps\to H$: $\tilde g_\eps(u)=g_\eps(u)$. As $\tilde g_\eps$ is a continuous linear map, there exists $\tilde g_\eps^*$: $H\to H_\eps$ such that $\inner{\tilde g_\eps(u)}{v}=\inner{u}{\tilde g_\eps^*(v)}_\eps$, $\forall u\in H_\eps$, $\forall v\in H$ and with $\norm{\tilde g_\eps^*}=\norm{\tilde g_\eps}$. Hence the map $\tilde g^*$: $\G_H\to\G$: $\tilde g^*([(u_\eps)_\eps])=[\tilde g_\eps^*(u_\eps)_\eps]$ is a well-defined continuous $\wt\C$-linear map and is the adjoint of $\tilde g$. By Corollary \ref{corol_isom_edged}, $\tilde g(\G)=g(\G_H)$ is edged in $\G_H$.
\end{proof}

\section{Variational inequalities in Hilbert $\wt{\R}$-modules}
\label{section_varia}
In the framework of Hilbert $\wt{\R}$-modules we now study variational inequalities involving a continuous and $\wt{\R}$-bilinear form. We will make use of the results proved in the previous sections in the context of Hilbert $\wt{\C}$-modules which can be easily seen to be valid for Hilbert $\wt{\R}$-modules. We begin with a general formulation in Theorem \ref{theo_varia} and then we concentrate on some internal versions in Proposition \ref{prop_suff} and Theorem \ref{theo_contrac}.

\begin{theorem}
\label{theo_varia}
Let $a(u,v)$ be a symmetric, coercive and continuous $\wt{\R}$-bilinear form on a Hilbert $\wt{\R}$-module $\G$. Let $C$ be a nonempty closed subset of $\G$ such that $\lambda C+(1-\lambda)C\subseteq C$ for all real generalized numbers $\lambda\in\{[(\eps^q)_\eps]\}_{q\in\N}\cup\{\frac{1}{2}\}$. For all $f\in\G$ such that the functional $I(u)=a(u,u)-2\inner{f}{u}$ has a close infimum on $C$ in $\wt{\R}$, there exists  a unique solution $u\in C$ of the following problem:
\beq
\label{prob_ineq}
a(u,v-u)\ge\inner{f}{v-u}, \qquad\qquad\text{for all $v\in C$.}
\eeq
\end{theorem}
\begin{proof}
Let $d$ be the close infimum of the functional $I$ on $C$ and $(u_n)_n\subseteq C$ a minimizing sequence such that $d\le I(u_n)\le d+[(\eps^n)_\eps]$. By means of the parallelogram law and the assumptions on $C$ we obtain that
\[
\begin{split}
\alpha\Vert{u_n-u_m}\Vert^2\le a(u_n-u_m,u_n-u_m)&=2a(u_n,u_n)+2a(u_m,u_m)-4a(\frac{u_n+u_m}{2},\frac{u_n+u_m}{2})\\
&=2I(u_n)+2I(u_m)-4I(\frac{u_n+u_m}{2})\\
&\le 2(d+[(\eps^n)_\eps]+d+[(\eps^m)_\eps]-2d)\\
&\le 2[(\eps^{\min{(m,n)}})_\eps]
\end{split}
\]
Since $\alpha$ is invertible it follows that $(u_n)_n$ is a Cauchy sequence and therefore it is convergent to some $u\in C$ such that $I(u)=\lim_{n\to \infty}I(u_n)=d$. 

For any $v\in C$ let us take $w=u+\lambda(v-u)$ with $\lambda=[(\eps^q)_\eps]$. By the properties of $C$ we know that $w\in C$ and $I(w)\ge I(u)$. It follows that
\begin{multline*}
a(u+\lambda(v-u),u+\lambda(v-u))-2\inner{f}{u+\lambda(v-u)}-a(u,u)+2\inner{f}{u}\\
=\lambda a(u,v-u)+\lambda a(v-u,u)+\lambda^2 a(v-u,v-u)-2\lambda\inner{f}{v-u}\ge 0
\end{multline*}
and since $\lambda$ is invertible,
\[
a(u,v-u)\ge\inner{f}{v-u}-\frac{1}{2}\lambda a(v-u,v-u).
\]
Letting $\lambda=[(\eps^q)_\eps]$ tend to $0$ in $\wt{\R}$ we conclude that $a(u,v-u)\ge\inner{f}{v-u}$ for all $v\in C$, or in other words that $u$ is a solution of our problem.

Finally, assume that $u_1,u_2$ are both solution in $C$ of the variational inequality problem \eqref{prob_ineq}. Then,
$a(u_1,u_1-u_2)\le\inner{f}{u_1-u_2}$, $-a(u_2,u_1-u_2)\le -\inner{f}{u_1-u_2}$ and  
\[
\alpha\Vert{u_1-u_2}\Vert^2\le a(u_1-u_2,u_1-u_2)\le 0.
\]
This means that $u_1=u_2$ and that the problem \eqref{prob_ineq} is uniquely solvable in $C$.
\end{proof}

\begin{corollary}
\label{corollary_eq}
Let $a(u,v)$ be a symmetric, coercive and continuous $\wt{\R}$-bilinear form on a Hilbert $\wt{\R}$-module $\G$. For all $f\in\G$ such that the functional $I(u)=a(u,u)-2\inner{f}{u}$ has a close infimum in $\wt{\R}$, there exists  a unique solution $u\in\G$ of the problem
\[
a(u,v)=\inner{f}{v}, \qquad\qquad\text{for all $v\in\G$.}
\]
\end{corollary}
\begin{proof}
Since Theorem \ref{theo_varia} applies to the case of $C=\G$ we have that there exists a unique $u\in \G$ such that
$a(u,v-u)\ge\inner{f}{v-u}$ for all $v\in\G$. This implies that $a(u,v)=\inner{f}{v}$ for all $v\in\G$. 
\end{proof}
Note that differently from Theorem \ref{LM_theo}, Corollary \ref{corollary_eq} does not require the particular structure $\inner{u}{g(v)}$ for the symmnetric $\wt{\R}$-bilinear form $a(u,v)$.

As a particular case of Theorem \ref{theo_varia} we obtain the following result for basic, symmetric and coercive $\wt{\R}$-bilinear forms on $\G_H$.

\begin{proposition}
\label{prop_suff}
Let $H$ be a real Hilbert space, $(C_\eps)_\eps$ a net of convex subsets of $H$, $C=[(C_\eps)_\eps]$ and $a$ be a basic, symmetric and coercive $\wt{\R}$-bilinear form on $\G_H$. If $C\neq\emptyset$ then, for all basic functionals $f$ on $\G_H$ there exists a unique solution $u\in C$ of the problem:
\[
a(u,v-u)\ge f(v-u), \qquad\qquad\text{for all $v\in C$.}
\]    
\end{proposition}
\begin{proof}
By Proposition \ref{riesz_int} and Theorem \ref{theo_riesz}, there exists $b\in\G_H$ such that $f(v)=\inner{b}{v}$, $\forall v\in\G_H$. Since $C$ is a closed and edged subset of $\G_H$ such that $\lambda C+(1-\lambda)C\subseteq C$ for all real generalized numbers $\lambda\in\{[(\eps^q)_\eps]\}_{q\in\N}\cup\{\frac{1}{2}\}$, in order to apply Theorem \ref{theo_varia} it suffices to prove that the functional $I(u)=a(u,u)-2\inner{b}{u}$ has a close infimum on $C$ in $\wt{\R}$. We fix representatives $(a_\eps)_\eps$ and $(b_\eps)_\eps$ of $a$ and $b$ respectively and we denote the corresponding net of functionals by $(I_\eps)_\eps$. From the coercivity of $a$ and Lemma \ref{lemma_coer} it follows that for each sufficiently small $\eps$ the inequality
\beq
\label{ineq_bound_2}
I_\eps(w)\ge\alpha'_\eps\Vert w\Vert-2c_\eps\Vert w\Vert=\big(\sqrt{\alpha'_\eps}\Vert{w}\Vert-\frac{c_\eps}{\sqrt{\alpha'_\eps}}\big)^2-\frac{1}{\alpha'_\eps}{c_\eps}^2\ge -\frac{1}{\alpha'_\eps}{c_\eps}^2 
\eeq
holds on $H$, where $\alpha'_\eps=\alpha_\eps-\eps^q$, $(c_\eps)_\eps$ is a representative of $c=\Vert b\Vert$ and $q\in\N$. Hence, $I_\eps$ has an infimum $d_\eps$ on $C_\eps$ such that $-{c_\eps}^2/{\alpha'_\eps}\le d_\eps$. Let $v_\eps\in C_\eps$ be such that $I_\eps(v_\eps)\le d_\eps + \eps^{1/\eps}$. From \eqref{ineq_bound_2} we see that for every moderate net of real numbers $(\lambda_\eps)_\eps$ there exists a moderate net $(\mu_\eps)_\eps$ such that $I_\eps(u_\eps)\ge\lambda_\eps$ as soon as $u_\eps\in H$ and $\Vert{u_\eps}\Vert^2\ge\mu_\eps$. Applying this to $\lambda_\eps=1+d_\eps$, we conclude that the net $(\Vert{v_\eps}\Vert^2)_\eps$ is moderate, and $v=[(v_\eps)_\eps]\in C$. It follows that the functional $I$ reaches its minimum $d=[(d_\eps)_\eps]$ on $C$ in $v$. The uniqueness of the solution follows as in the proof of Theorem \ref{theo_varia}.
\end{proof}

\begin{remark}
\label{rem_for_appl}
Note that Proposition \ref{prop_suff} makes use of the completeness of $\G_H$ which holds even if $H$ is not complete (see \cite[Proposition 3.4]{Garetto:05a}).
\end{remark}
We extend now Proposition \ref{prop_suff} to $\wt{\R}$-bilinear forms which are not necessarily symmetric by making use of some contraction techniques.   

\begin{theorem}
\label{theo_contrac}
Let $H$ be a real Hilbert space, $(C_\eps)_\eps$ a net of convex subsets of $H$, $C=[(C_\eps)_\eps]$ and $a$ be a basic and coercive $\wt{\R}$-bilinear form on $\G_H$. If $C\ne\emptyset$, then, for all basic functionals $f$ on $\G_H$ there exists a unique solution $u\in C$ of the problem:
\[
a(u,v-u)\ge f(v-u), \qquad\qquad\text{for all $v\in C$.}
\]
\end{theorem}
\begin{proof}
By Proposition \ref{riesz_int} and Theorem \ref{theo_riesz}, there exists $c\in\G_H$ such that $f(v)=\inner{c}{v}$, $\forall v\in\G_H$. By Proposition \ref{prop_riesz_b_int} and Theorem \ref{theo_riesz_b}, there exists a basic $\wt{\R}$-linear map $T$: $\G_H\to\G_H$ such that $a(u,v)=\inner{Tu}{v}$, $\forall u,v \in\G_H$.
We look for $u\in C$ satisfying the inequality
\[\inner{T u}{v-u}\ge \inner{c}{v-u},\qquad\forall v\in C.\]
For any $\rho\in\wt{\R}$, with $\rho\ge 0$ invertible,
the inequality is equivalent with
\[
\inner{(\rho c-\rho T u + u) - u}{v-u}\le 0,\qquad\forall v\in C.
\]
By Proposition \ref{prop_projec_int}$(i)$, $C$ is closed and edged; further, $\lambda C + (1-\lambda) C\subseteq C$, $\forall\lambda\in \wt{\R}$ with $0\le\lambda\le 1$. So by Proposition \ref{prop_charac}, we search $u\in \G_H$ with $u=P_{C}(\rho c-\rho T u + u)$ (for a suitable $\rho$ that will be determined below).

By the basic structure of $T$ we know that there exists a moderate net $(M_\eps)_\eps$ and $\eta_1\in(0,1]$ such that $\Vert{T_\eps u}\Vert\le M_\eps\Vert{u}\Vert$, $\forall u\in H$, $\forall\eps\in(0,\eta_1]$. By coercivity of $a$, there exists a moderate net $(\alpha_\eps)_\eps$ and $m\in\N$ with $\alpha_\eps\ge\eps^m$, $\forall\eps$ and there exists $\eta_2\in (0,1]$, such that $\inner{T_\eps u}{u}\ge\alpha_\eps \Vert u\Vert^2$, $\forall u\in H$, $\forall\eps\in(0,\eta_2]$ (Lemma \ref{lemma_coer}). Let $\eta=\min(\eta_1,\eta_2)$. 
Fix $\eps\in(0,\eta]$. Let $\rho_\eps=\frac{\alpha_\eps}{M_\eps^2}$ and
\[S_\eps: \overline C_\eps\to \overline C_\eps: S_\eps(v)=P_{\overline C_\eps}(\rho_\eps c_\eps-\rho_\eps T_\eps v + v).\]
For $v_1$, $v_2\in \overline C_\eps$, by the properties of $P_{\overline C_\eps}$,
\[\Vert S_\eps(v_1)-S_\eps(v_2)\Vert\le\Vert(v_1-v_2) - \rho_\eps(T_\eps v_1- T_\eps v_2)\Vert,\]
so
\begin{align*}
\Vert S_\eps(v_1)-S_\eps(v_2)\Vert^2 &\le \Vert v_1 - v_2\Vert^2-2\rho_\eps\inner{T_\eps v_1 - T_\eps v_2}{v_1-v_2}+ \rho_\eps^2\Vert T_\eps v_1 -T_\eps v_2\Vert^2\\
&\le (1-2\rho_\eps\alpha_\eps + \rho_\eps^2M_\eps^2)\Vert v_1 - v_2\Vert^2=\left(1-\frac{\alpha_\eps^2}{M_\eps^2}\right)\Vert v_1 - v_2\Vert^2.
\end{align*}
So $S_\eps$ is a contraction. Let $w\in C$ with representative $(w_\eps)_\eps$, $w_\eps\in C_\eps$, $\forall\eps$. Denoting the contraction constant by $k_\eps=\big(1-\frac{\alpha_\eps^2}{M_\eps^2}\big)^{1/2}$, by the properties of a contraction, $\Vert{S_\eps^n(w_\eps)-w_\eps}\Vert \le \frac{1}{1-k_\eps}\Vert S_\eps(w_\eps)-w_\eps\Vert$, $\forall n\in\N$. In particular, for the fixed point $u_\eps$ of $S_\eps$ in $\overline C_\eps$, also $\Vert u_\eps - w_\eps\Vert\le\frac{1}{1-k_\eps}\Vert S_\eps(w_\eps)-w_\eps\Vert$. Hence
\[
\Vert u_\eps\Vert\le\Vert w_\eps\Vert + \frac{1}{1-k_\eps} \Vert S_\eps(w_\eps)-w_\eps\Vert.
\]
Now there exists $m\in\N$ such that $k_\eps^2\le 1- \eps^m$, $\forall\eps$, hence $k_\eps\le \sqrt{1-\eps^{m}}\le 1-\frac{\eps^{m}}{2}$, $\forall\eps$. Further, if $\rho=[(\rho_\eps)_\eps]$, then $[(S_\eps(w_\eps))_\eps] = P_C(\rho c-\rho T w + w)$ by Proposition \ref{prop_projec_int}. So $(\Vert u_\eps\Vert)_\eps$ is a moderate net, hence
$(u_\eps)_\eps$ represents some $u\in C$. Similarly, as $S_\eps(u_\eps)=u_\eps$, $\forall\eps\le\eta$, we have $u=P_C(\rho c - \rho T u + u)$ for $\rho=\big[\big(\frac{\alpha_\eps}{M_\eps^2}\big)_\eps\big]\in\wt{\R}$ with $\rho\ge 0$ and invertible, as required.

The uniqueness of the solution follows as in the proof of Theorem \ref{theo_varia}.
\end{proof}

\section{Applications}
\label{sec_appl}
We conclude the paper by applying the theorems on variational equalities and inequalities of Section \ref{section_varia} to some concrete problems coming from partial differential operators with highly singular coefficients. The generalized framework within which we work allows us to approach problems which are not solvable classically and to get results consistent with the classical ones when the latter exist.
\subsection{The generalized obstacle problem}
In the sequel $\Om\subseteq\R^n$ is assumed to be open, bounded and connected with smooth boundary $\partial\Om$. We consider a net $(\psi_\eps)_\eps\in H^1(\Om)^{(0,1]}$ such that $\psi_\eps\le 0$ a.e. on $\partial\Om$ for all $\eps$, and we define the set 
\[
C_\psi=\{[(u_\eps)_\eps]\in\G_{H^1_0(\Om)}:\,  \forall\eps\quad u_\eps\ge\psi_\eps\ \text{a.e. on}\, \Om\}.
\]
One can easily see that $C_\psi$ is a nonempty internal subset of the Hilbert $\wt{\C}$-module $\G_{H^1_0(\Om)}$ given by a net of convex subsets of $H^1_0(\Om)$. Note that the net $(\psi_\eps)_\eps$ can be generated by a highly singular obstacle $\psi$ regularized via convolution with the mollifier $\varphi_\eps$, where $\varphi\in\Cinfc(\Om)$, $\int\varphi\, dx=1$ and $\varphi_\eps(x)=\eps^{-n}\varphi(x/\eps)$. For instance, on $\Om=\{x\in\R^n:\, |x| < 1\}$ one can take an arbitrary $\psi\in \E'(\Om)$. From the structure theorem for distributions with compact support we obtain that there exists some $\eta\in(0,1]$ such that $(\psi\ast\varphi_\eps)_{\eps\le\eta}$ is a $H^1(\Om)$-moderate net. 


Let $(a_{i,j,\eps})_\eps$ be moderate nets of $L^\infty$-functions on $\Om$ such that
\beq
\label{def_a_repr}
\lambda_\eps^{-1}\xi^2\le \sum_{i,j=1}^n a_{i,j,\eps}(x)\xi_j\xi_i\le\lambda_\eps \xi^2
\eeq
holds for some positive and invertible $[(\lambda_\eps)_\eps]\in\wt{\R}$ and for all $(x,\xi)\in\Om\times\R^n$. From \eqref{def_a_repr} it follows that
\beq
\label{def_a_ob}
a(u,v)=\biggl[\biggl(\int_\Om\sum_{i,j=1}^n a_{i,j,\eps}(x)\partial_{x_j}u_\eps(x)\partial_{x_i}v_\eps(x)\, dx\biggr)_\eps\biggr]
\eeq
is a well-defined basic $\wt{\R}$-bilinear form on $\G_{H^1_0(\Om)}$. Before proceeding we recall that from \cite{Garetto:05a}[Proposition 3.22] the space $\G_{H^{-1}(\Om)}$ coincides with the set of basic functionals in $\L(\G_{H^1_0(\Om)},\wt{\C})$.

We are now ready to state the following theorem.
\begin{theorem}
\label{theo_obstacle}
Let $a$ as in \eqref{def_a_ob}. For any $f\in\G_{H^{-1}(\Om)}$ there exists a unique solution $u\in C_\psi$ of the problem:
\[
a(u,v-u)\ge f(v-u), \qquad\qquad\text{for all $v\in C_\psi$.}
\] 
\end{theorem}
\begin{proof}
In order to apply Theorem \ref{theo_contrac} we have to prove that the $\wt{\R}$-bilinear form $a$ is coercive, in the sense of Definition \ref{def_coercivity}. The condition \eqref{def_a_repr} on the coefficients of $a$ and the Poincar\'e inequality yield that
\[
a(v,v)\ge \lambda^{-1}\Vert v\Vert^2_{H^1_0(\Om)},
\]
is valid for all $v\in\G_{H^1_0(\Om)}$. This completes the proof.
\end{proof}

\begin{remark}
\label{rem_classic}
When the obstacle $\psi$ and the coefficients $a_{i,j}$ are classical, for any $f\in H^{-1}(\Om)$ the problem $a(u,v-u)\ge f(v-u)$ can be classically settled in $H^1_0(\Om)$ by looking for a solution $u$ in $C^{\rm{cl}}_\psi:=\{u\in {H^1_0(\Om)}:\,  u\ge\psi\ \text{a.e. on}\, \Om\}$. Let $u_0\in C^{\rm{cl}}_\psi$ such that 
\beq
\label{class_solva}
a(u_0,v-u_0)\ge f(v-u_0)
\eeq
for all $v\in C^{\rm{cl}}_\psi$. Note that by embedding $H^{-1}(\Om)$ into $\G_{H^{-1}(\Om)}$ by means of $f\to [(f)_\eps]$, we can study the previous obstacle problem in the generalized context of $\G_{H^{1}_0(\Om)}$. By Theorem \ref{theo_obstacle} we know that there exists a unique $u\in C_\psi:=\{[(v_\eps)_\eps]\in\G_{H^1_0(\Om)}:\,  \forall\eps\quad v_\eps\ge\psi\ \text{a.e. on}\, \Om\}$ such that 
\beq
\label{gen_solva}
a(u,v-u)\ge[(f)_\eps](v-u),
\eeq
for all $v\in C_\psi$. By the fact that \eqref{gen_solva} is uniquely solvable it follows that $u$ coincides with the classical solution, i.e., $u=[(u_0)_\eps]$. Indeed, since for any $v\in C_\psi$ we can find a representative $(v_\eps)_\eps$ such that $v_\eps\in C^{\rm{cl}}_\psi$ for all $\eps$, from \eqref{class_solva} we have that
\[
a(u_0,v_\eps-u_0)\ge f(v_\eps-u_0)
\]
for all $v=[(v_\eps)_\eps]\in C_\psi$.
\end{remark}

\begin{example}
\label{ex_conv}
When the coefficients $a_{i,j}$ are not bounded, the obstacle problem is in general not solvable in the Sobolev space $H^1_0(\Om)$. In this case one can think of regularizing the coefficients by convolution with a mollifier $\varphi_\eps$ and looking for a generalized solution in some subset of $\G_{H^{1}_0(\Om)}$. For instance, let $\mu_i$ be finite measures on $\R^n$ with $\mu_{i}\ge c\chi_V$, for $i=1,\dots, n$, where $V$ is a neighbourhood of $\overline\Omega$, $\chi_V$ denotes the characteristic function of $V$ and $c\in\R$, $c >0$. Let us take $a_{i,j}=0$ when $i\neq j$ and $a_{i,i}=\mu_i$, 
for $i,j=1,\dots,n$. If $\varphi$ is a nonnegative function in $\Cinfc(\R^n)$ such that $\int\varphi=1$ and $\varphi_\eps(x)=\eps^{-n}\varphi(x/\eps)$, we obtain for sufficiently small $\eps$ and $x\in\Omega$ that 
\[
c \le \mu_{i}\ast\varphi_\eps (x)
\le \mu_{i}(\R^n)
\Vert \varphi_\eps\Vert_{L^\infty}\le c'\eps^{-n},
\] 
for some constant $c'\in\R$ depending on $\varphi$ and the measures $\mu_{i}$. It follows that setting $a_{i,i,\eps}(x)=\mu_{i}\ast\varphi_\eps(x)$, the net 
\[
\int_\Om\sum_{i=1}^n a_{i,i,\eps}(x)\partial_{x_i}u(x)\partial_{x_i}v(x)\, dx
\]
of bilinear forms on $H^{1}_0(\Om)$ generates a basic and coercive $\wt{\R}$-bilinear form on $\G_{H^1_0(\Om)}$. For a generalized $C_\psi$ as at the beginning of this subsection and any $f\in \G_{H^{-1}(\Om)}$, the corresponding obstacle problem is uniquely solvable. 
\end{example}

\subsection{A generalized Dirichlet problem}
We want to study the homogeneous Dirichlet problem 
\beq
\label{Dirich}
-\nabla\cdot( A\nabla u)+a_0 u=f\quad \text{in $\Om$},\qquad\qquad u=0\quad \text{on $\partial\Om$},
\eeq
where $A=[(A_\eps)_\eps]$ and $a_0=[(a_{0,\eps})_\eps]$ are $\G_{L^\infty(\Om)}$ generalized functions satisfying the following conditions: there exist some positive moderate nets $(\lambda_\eps)_\eps$, with moderate inverse $(\lambda_\eps^{-1})_\eps$, and $(\mu_\eps)_\eps$ such that for all $x\in\Om$ and $\eps\in(0,1]$,
\beq
\label{cond_A}
\lambda_\eps^{-1}\le A_\eps(x)\le\lambda_\eps,
\eeq
and
\beq
\label{cond_a_0}
0\le a_{0,\eps}(x)\le \mu_\eps.
\eeq
Let $f\in\G_{H^{-1}(\Om)}$. We formulate the problem \eqref{Dirich} in $\G_{H^{-1}(\Om)}$. Its variational formulation is given within the Hilbert $\wt{\R}$-module $\G_{H^1_0(\Om)}$ in terms of the equation 
\beq
\label{vario_Dirich}
a(u,v)=f(v),\qquad\qquad \text{for all $v\in\G_{H^1_0(\Om)}$},
\eeq
where 
\[
a(u,v)=\biggl[\biggl(\int_\Om A_\eps(x)\nabla u_\eps(x)\cdot\nabla v_\eps(x)\, dx+\int_\Om a_{0,\eps}(x)u_\eps(x)v_\eps(x)\, dx\biggr)_\eps\biggr].
\]
From \eqref{cond_A}, \eqref{cond_a_0} and the Poincar\'e inequality it follows that $a$ is a basic and coercive $\wt{\R}$-bilinear form on $\G_{H^1_0(\Om)}$. Recalling that $f$ is a basic functional on $\G_{H^1_0(\Om)}$, an application of Proposition \ref{LM_int} yields the desired solvability in $\G_{H^1_0(\Om)}$.
\begin{theorem}
\label{theo_Dirich}
For any $f\in\G_{H^{-1}(\Om)}$ the variational problem \eqref{vario_Dirich} is uniquely solvable in $\G_{H^1_0(\Om)}$.
\end{theorem}

\begin{example}
\label{exam_Hd}
Let us consider the one-dimensional Dirichlet problem given by the equation
\beq
\label{dirich_sing}
-(H(x)u')'+\delta u=f,
\eeq
in some interval $I=(-a,a)$. One can think of approximating the singular coefficients which appear in \eqref{dirich_sing} by means of moderate nets of $L^\infty$ functions which satisfy the conditions \eqref{cond_A} and \eqref{cond_a_0}. Let $(\nu_\eps)_\eps\in\R^{(0,1]}$ be a positive moderate net with moderate inverse $(\nu_\eps^{-1})_\eps$ such that $\nu_\eps\to 0$ if $\eps\to 0$. We easily see that $A_\eps(x)$ equal to $1$ for $x\in(0,a)$ and to $\nu_\eps$ for $x\in(-a,0]$, fulfills \eqref{cond_A}, while $a_{0,\eps}(x)=\varphi_\eps(x)$, $x\in I$, with $\varphi\in\Cinfc(\R)$ as in Example \ref{ex_conv}, has the property \eqref{cond_a_0}. From Theorem \ref{theo_Dirich} we have that for any $f\in\G_{H^{-1}(I)}$ the variational problem associated to \eqref{dirich_sing} is uniquely solvable in $\G_{H^1_0(I)}$. It is clear that a similar result can be obtained for other functions $h\ge 0$ with zeroes instead of the Heaviside-function $H$.
\end{example}
In a similar way we can deal with the inhomogeneous problem
\beq
\label{Dirich_2}
-\nabla\cdot( A\nabla u)+a_0 u=f\quad \text{in $\Om$},\qquad\qquad u=g\quad \text{on $\partial\Om$},
\eeq
where we assume that $g\in\G_{\mathcal{C}^1(\partial{\Omega})}$. Let $(g_\eps)_\eps$ be a representative of $g$ and  $(\wt{g_\eps})_\eps$ a net in $\M_{H^1(\Om)\cap\mathcal{C}(\overline{\Omega})}$ such that $g_\eps=\wt{g_\eps}$ on $\partial\Omega$. Defining,   
\[
C=\{v=[(v_\eps)_\eps]\in\G_{H^1(\Om)}:\ v_\eps-\wt{g_\eps}\in H^1_0(\Om)\},
\]
a variational formulation of \eqref{Dirich_2} is  
\beq
\label{vario_Dirich_2}
a(u,v)=f(v),\qquad\qquad \text{for all $v\in C$}. 
\eeq
Under the assumptions \eqref{cond_A} and \eqref{cond_a_0} for $A$ and $a_0$ we obtain the following result.
\begin{theorem}
\label{theo_Dirich_2}
For any $f\in\G_{H^{-1}(\Om)}$ the variational problem \eqref{vario_Dirich_2} is uniquely solvable in $C$.
\end{theorem}

\newcommand{\SortNoop}[1]{}

\end{document}